\documentclass[11pt]{amsart}

\usepackage[normalem]{ulem}
\usepackage{booktabs} 
\usepackage{amsmath}
\usepackage{textcomp} 
\usepackage[font=footnotesize]{subcaption}
\usepackage[]{algorithm2e}
\usepackage{scalerel}
\usepackage{stackengine,wasysym} 
\usepackage{mathdots} 

\usepackage{algorithmic}
\newlength\bshft 
\def\fakebold#1{\ThisStyle{\ooalign{$\SavedStyle#1$\cr%
  \kern-\bshft$\SavedStyle#1$\cr%
  \kern\bshft$\SavedStyle#1$}}}

\usepackage[usenames,dvipsnames]{xcolor}
\usepackage[mathscr]{euscript}
\usepackage{array, amsmath, amsfonts, amssymb, amsthm, geometry, graphics, graphicx, multicol, multirow, enumerate, float, floatflt, fullpage, tikz, outlines, url}
\usepackage[utf8]{inputenc}
\usetikzlibrary{shapes.multipart, matrix, positioning, arrows}
\newcolumntype{P}[1]{>{\centering\arraybackslash}p{#1}}

\usepackage{wrapfig}

\usepackage{tikz}

\usepackage{lineno}
\usepackage{diagbox} 
\usetikzlibrary{shapes,decorations}

\newsavebox\crossing
\begin{lrbox}{\crossing}
\begin{tikzpicture}[scale=1]
\useasboundingbox (0,0) rectangle (2,1);
\draw (1,0)--(0,1);
\draw (0,0)--(.4,.4);
\draw (.6,.6)--(1,1);
\end{tikzpicture}
\end{lrbox}

\newsavebox\poscrossing
\begin{lrbox}{\poscrossing}
\begin{tikzpicture}[scale=1]
\useasboundingbox (0,0) rectangle (2,1);
\draw[->] (1,0)--(0,1);
\draw (0,0)--(.4,.4);
\draw[->] (.6,.6)--(1,1);
\end{tikzpicture}
\end{lrbox}

\newsavebox\negcrossing
\begin{lrbox}{\negcrossing}
\begin{tikzpicture}[scale=1]
\draw (1,0)--(.6,.4);
\draw[->] (.4,.6)--(0,1);
\draw[->] (0,0)--(1,1);
\end{tikzpicture}
\end{lrbox}

\newsavebox\zerosmoothing
\begin{lrbox}{\zerosmoothing}
\begin{tikzpicture}[scale=1]
\draw (0,1) to [out=-45,in=-135](1,1);
\draw (0,0) to [out=45,in=135](1,0);
\end{tikzpicture}
\end{lrbox}

\newsavebox\onesmoothing
\begin{lrbox}{\onesmoothing}
\begin{tikzpicture}[scale=1]
\draw (0,0) to [out=45,in=-45](0,1);
\draw (1,0) to [out=135,in=-135] (1,1);
\end{tikzpicture}
\end{lrbox}

\newsavebox\posnegconv
\begin{lrbox}{\posnegconv}
\begin{tikzpicture}[scale=.5]
\node at (3,3){\usebox\poscrossing};

\node at (6.5,3){\usebox\negcrossing};

\node at (1.8,1.2){$-$};
\node at (6.5,1.2){$+$};
\end{tikzpicture}
\end{lrbox}

\newsavebox\crossingconv
\begin{lrbox}{\crossingconv}
\begin{tikzpicture}[scale=.5]
\node at (2,3){\usebox\zerosmoothing};
\node at (12.4,1){$L_1$};

\draw[->, thick] (5.5,3)--(3.5,3);
\node at (4.5,3.4){$0$};

\node at (7.3,3){\usebox\negcrossing};
\node at (7.3,1){$L_+$};

\draw[->, thick] (8.8,3)--(10.8,3);
\node at (9.8,3.4){$1$};

\node at (12.4,3){\usebox\onesmoothing};
\node at (2,1){$L_0$};

\end{tikzpicture}
\end{lrbox}

\usetikzlibrary{decorations.pathreplacing, hobby, calc, knots, shapes, external}

\usepackage{environ}
\makeatletter
\newsavebox{\measure@tikzpicture}
\NewEnviron{scaletikzpicturetowidth}[1]{%
  \def\tikz@width{#1}%
  \begin{lrbox}{\measure@tikzpicture}%
  \BODY
  \end{lrbox}%
  \pgfmathparse{#1/\wd\measure@tikzpicture}%
  \BODY
}
\makeatother

\usepackage{xcolor}
\usepackage[utf8]{inputenc}
\usepackage{hyperref}

\usepackage[foot]{amsaddr}

\graphicspath{ {figures/}}

\title{Mapper--type algorithms for complex data and relations }

\author{Pawe{\l} D{\l}otko$^1$, Davide Gurnari$^1$ and Radmila Sazdanovic$^{1,2}$}

\email{pdlotko@impan.pl, dgurnari@impan.pl, rsazdanovic@math.ncsu.edu}

\address{$^1$Dioscuri Centre in Topological Data Analysis,
Mathematical Institute PAN,
Warsaw, PL}
\address{$^2$Department of Mathematics,
North Carolina State University,
Raleigh, NC}

\thanks{PD and DG acknowledge support by Dioscuri program initiated by the Max Planck Society, jointly managed with the National Science Centre (Poland), and mutually funded by the Polish Ministry of Science and Higher Education and the German Federal Ministry of Education and Research. RS is partially supported by NSF grant DMS-1854705.}

\date{\today}

\newcommand{\NEW}[1]{#1}
\newcommand{\OLD}[1]{{\leavevmode\color{WildStrawberry} \sout{#1}}}

\usepackage[numbers]{natbib}
\usepackage{graphicx}

\newcommand*{\ifNotPictures}{}% 

\begin{document}

\begin{abstract}
Mapper and Ball Mapper are Topological Data Analysis tools used for exploring  high dimensional point clouds and visualizing scalar--valued functions on those point clouds. 
Inspired by open questions in knot theory, new features are added to Ball Mapper that enable encoding \NEW{of the structure, internal relations and symmetries of the point cloud.} 
Moreover, the strengths of Mapper and Ball Mapper constructions are combined to create a tool for comparing high dimensional data descriptors of a single dataset. \NEW{
This new hybrid algorithm,  Mapper on Ball Mapper, is applicable to high dimensional lens functions.}
 As a proof of concept we include applications to knot and game theory, as well as material science \NEW{and cancer research}.
\end{abstract}

\maketitle

\section{Introduction}\label{sec:intro}
Mapper~\cite{singh2007topological} and Ball Mapper~\cite{ball_mapper} algorithms are Topological Data Analysis (TDA) tools for exploring and visualizing data. 
The Mapper algorithm was introduced in 2007 and since then it has been used in a number of different applications, most notably medical data analysis~\cite{nicolau},~\cite{diabetes} and  material science~\cite{zeolites}. Ball Mapper~\cite{ball_mapper} was introduced in 2019 as an easy-to-use effective alternative to the Mapper algorithm, see  applications to data in economics~\cite{altman}.

Both techniques take a point cloud $X$, possibly with a scalar--valued function $f : X \rightarrow \mathbb{R}$, as an input, and return an abstract graph $G = (V,E)$, with a induced function $\hat{f} : V \rightarrow \mathbb{R}$. The ''layout'' of $G$ resembles the multidimensional layout of the input set $X$. In addition, visualising $\hat{f}$ over the vertices of $G$, typically using an appropriate color scale, provides insights into properties of the input function $f$. One of the first and most famous examples\NEW{, further analyzed in this paper,} is the work in~\cite{nicolau}, where the input point cloud consists of  gene expression data for breast cancer patients, and the function $f$ determine the survival rate. Mapper associated rare cancer 
with subtypes with $100\%$ survival  rate and allowed the authors to characterize its genetic profile. 

\NEW{This research introduces novel techniques that broaden and improve the scope and applicability of mapper-type algorithms to utilize additional structure of data and visualize the maps between datasets. }
Main contributions include: 
\begin{enumerate}
\item \emph{Equivariant Ball Mapper EqBM}: the most natural choice of mapper-type algorithms when the input data admits an action of a group of isometries or distance preserving bijection. The structure of the resulting Ball Mapper graph reflects the action of the input isometries. 
\item \emph{MappingMappers}: Mapper and Ball Mapper--based representations of high dimensional datasets $X$ and $Y$ can be used to visualize a given relation $f : X \rightarrow Y$. For example, if  $X$ and $Y$ are values of different descriptors of the \NEW{ object of interest, such as different knots, materials, or cancer patients. } MappingMappers can be used to compare  descriptors' relevance 
and discover potential correlations or dependencies. 

\item \emph{Mapper on Ball Mapper MoBM}: Extension of the Mapper algorithm that allows the use of high dimensional lens functions $f : X \rightarrow \mathbb{R}^n$, with  $n >> 1$. In the proposed technique, a Ball Mapper graph of the range of $f$ is built to obtain an \emph{adaptive cover} which is then used to construct the final Mapper graph.
\end{enumerate}

The main running and motivating example in this paper is an analysis of data coming from knot theory which has recently opened up to big data techniques such as machine learning \cite{hughes2016neural,jejjala2019deep,ward2018using,levitt2019big, davies2021advancing}. Knot theory point clouds created from polynomial knot invariants such as the Alexander, Jones and HOMFLYPT, are perfect for showcasing the strengths of the Equivariant, MappingMappers and Mapper on Ball Mapper algorithms.  Additional examples related to game theory 
\NEW{material science and cancer research are discussed to highlight the presented algorithms. however the  detailed interpretations, as well as in depth analysis of knot theory results, are omitted as they do not fit in the scope of this paper}.

Techniques developed in this work are accompanied with sample public-domain implementations and have wide applicability in  different areas of science. 
The software implementing the discussed techniques and the interactive visualizations of all the plots in this paper are  available at the webpage \url{https://dioscuri-tda.org/BallMapperKnots.html} and in~\cite{zenodo}.

The paper is organized as follows.  Section~\ref{sec:Background} provides the necessary background on the Mapper and Ball Mapper algorithms. Section~\ref{sec:new_developement_mapper} focuses on adaptations and new developments of Mapper algorithms that are applicable to any point cloud.
In Section~\ref{sec:equivarinat_ball_mapper}  we develop a version of Ball Mapper that takes into account symmetries (global or partial) of the data.  Section~\ref{sec:BonBM} describes a way to  construct Mapper graphs for lens functions with high dimensional domains and codomains. In Section~\ref{sec:relation_between_mapper} we discuss how to combine the strengths of the  Mapper and Ball Mapper  algorithms to analyze data, relations between high dimensional point clouds, and visualize maps between different datasets.
Section~\ref{sec:KnotIntro}
provides informal minimal background information about relevant knot invariants while Section \ref{sec:Data} describes a way to turn knot invariants into point clouds so they can  be analyzed by TDA techniques. Section~\ref{sec:KnotDataBM} present the results of application of the the presented Mapper--type techniques to the space of knots and their invariants, while Section~\ref{sec:comparisons} focuses on the comparison of knot invariants using techniques from Section~\ref{sec:relation_between_mapper} \NEW{and shows how the presented technique can benefit theoretical  mathematics}. 
Finally, Section~\ref{sec:further_examples} provides additional examples of the use of the proposed techniques in game theory, materials science, \NEW{ and cancer research.}
\ifdefined\showOldText 
\OLD{and materials science}.
\fi

\section{Background}\label{sec:Background}

 In this section we present the Mapper algorithm~\cite{singh2007topological}, a standard tool of Topological Data Analysis, as well as the recently developed Ball Mapper~\cite{ball_mapper}.

\ifdefined\ifNotPictures
\begin{figure}[h]
    \centering
    \includegraphics[width=0.6\textwidth]{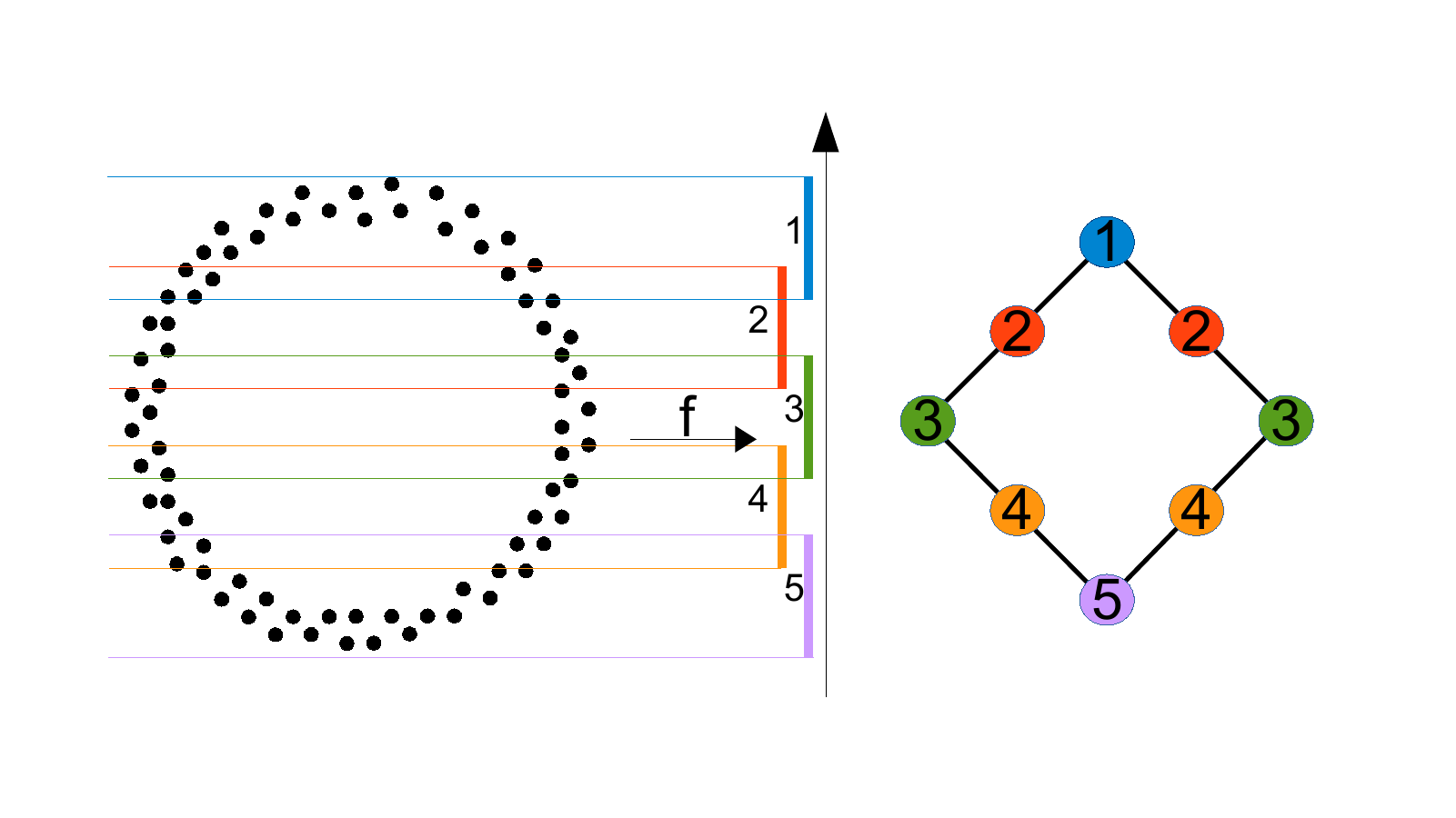}
    \caption{The Mapper construction illustration. The input is the 2-dimensional point cloud $X$ in a shape of a circle  shown on the left. The function $f : X \rightarrow \mathbb{R}^1$ is a projection of a point to its second coordinate. We cover the range of $f$ with five intervals, enumerated from $1$ to $5$, then compute the inverse image of each interval. In this example, the inverse images of intervals $1$ and $5$ contain one cluster, while inverse images of intervals $2,3$ and $4$ contain two clusters each. Moreover there are obvious connections between the clusters in the inverse image of intervals $i$ and $i+1$ for $i \in {1,2,3,4}$. They give rise to edges in the Mapper graph presented on the right. 
    Note that the (non-unique) enumeration of vertices comes from the enumeration of elements of the cover $\mathcal{C}$.}
    \label{fig:idea_of_conventional_mapper}
\end{figure}
\fi 
\subsection{Mapper}
\label{sec:conventional_mapper}
The \emph{Mapper algorithm}, introduced in 2007 in~\cite{singh2007topological}, is one of the most recognizable tools of Topological Data Analysis. It can be considered a discrete approximation of a \emph{Reeb graph}~\cite{reeb, Mapper2reeb}.
An input for Mapper is %Let us have 
a collection of points $X$, typically embedded in some high dimensional space,  and  a function $f : X \rightarrow \mathbb{R}^n$, referred to as the \emph{lens function}. Typically, the range of the lens function is 1-dimensional.
Next, construct the cover of the range $f(X) \subset \mathbb{R}^n$ by a \emph{collection of overlapping cubes of intervals} $\mathcal{C}$. For $n=1$, $\mathcal{C}$ usually consists of $k$ intervals covering $f(X)$ 
such that two constitutive intervals overlap on $p$ percent of their length.
The number of intervals (k) and the overlapping percentage (p),
often refereed to as \emph{resolution} and \emph{gain},
are additional input parameters of the Mapper algorithm. 

For each element $I \in \mathcal{C}$, consider  $f^{-1}(I) \subset X$ and search for the clusters therein. The clustering algorithm used for that purpose is yet another parameter of the Mapper construction. For $I \in \mathcal{C}$ let $C_I$ indicate the collection of clusters found in $f^{-1}(I)$. Each cluster in $\bigcup_{I \in \mathcal{C}} C_{I}$ corresponds to a vertex of an abstract graph $G = (V,E)$. Given a cluster $A \in C_{I}$, the vertex corresponding to $A$ is denoted by $v(A)$. 
\NEW{An undirected graph $G$ is constructed using the following rule: for any two vertices $v(A)$ and $v(B)$ corresponding to clusters $A \in C_{I}$ and $B \in C_{J}$, an undirected edge is placed between $v(A)$ and $v(B)$ if and only if $A \cap B \neq \emptyset$. }
The resulting graph $G = (V,E)$ is called a \emph{Mapper graph}. An illustration of the Mapper construction on a simple 2-dimensional input is shown in Figure~\ref{fig:idea_of_conventional_mapper}.

\ifdefined\showOldText 
\OLD{Mapper graph $G$ can be used to further analyze the original point cloud by coloring it by a function of interests that presents some other property of the data. Given a function $g : X \rightarrow \mathbb{R}$ define a \emph{induced function} $\hat{g} : G \rightarrow \mathbb{R}$, on the the Mapper graph $G$ defined above. Since each vertex $v(A)$ in $G$ corresponds to a cluster $A \in \bigcup_{I \in \mathcal{C}} C_{I}$, we define $\hat{g}( v(A) ) = \frac{ \sum_{x \in A} g(x) }{ |A| }$, i.e. take an average value of a function $g$ on the cluster $A$ as the value of the induced function on the vertex corresponding to $A$. This induced function is subsequently visualized using a color scale. This way, intuitively, the shape of Mapper graph $G$ is a reflection of the shape of the input point cloud $X$ and the coloring of $G$ provides an insight into the variation of a scalar--valued function on the original point cloud. }
\fi

\NEW{While Mapper is a well established tool for data analysis its stability with respect to perturbation of input parameters is still open for explorations. Main results in this direction focus on convergence of Mapper graphs to the Reeb graph~\cite{reeb} of the manifold from which the points are sampled from, when the number of sampled points goes to infinity \cite{DBLP:journals/jact/BrownBMW21, JMLR:v19:17-291}. However, the practical scope of those results is limited, especially when dealing with finite noisy samples. 

Mapper and Ball Mapper are data visualization tools that provide the initial understanding of the data. Due to the stability issues, the Mapper and Ball Mapper graphs should be computed for a range of parameters and permutation of the input data,  for example see Figure~\ref{fig:stability_Jones_construction} and~\ref{fig:stability_Jones_construction_epsilon} and compared using MappingMappers or Mapper on Ball Mapper construction. The validity of any findings should be confirmed using other methods.}
\subsection{Ball Mapper}
\label{sec:ball_mapper}
The \emph{Ball Mapper Algorithm} introduced in~\cite{ball_mapper}~\NEW{\cite{altman}}, provides a conceptually different and simpler way to obtain a cover of the input cloud $X$ than the original Mapper. 
Starting from a point cloud $X$, and a constant $\epsilon > 0$,
a subset $L \subset X$ is selected having the property that for every $x \in X$, there exists $l \in L$ such that $d(x,l) \leq \epsilon$. Such a subset $L$ is called an \emph{$\epsilon$--net} of $X$ \NEW{and  its points are referred to as "landmarks"}. Algorithm~\ref{alg:greedy_epsilon_net} is an example of a greedy algorithm to compute an $\epsilon$--net. Other algorithms can also be used to compute $\epsilon$--nets~\cite{mustafa:LIPIcs:2019:10663}.

\begin{algorithm}[H]
\KwData{Point cloud $X$, $\epsilon > 0$}
\KwResult{$L \subset X$, an $\epsilon$-net}
$L = \emptyset$ ;\\
\While{There exists a point $x \in X$ farther than $\epsilon$ from every point in $L$}{
  $L = L \cup x$ ;
}
\Return L
\label{alg:greedy_epsilon_net}
\caption{Greedy $\epsilon$--net ~\cite{ball_mapper}.}
\end{algorithm}
\medskip

Note that $X \subset \bigcup_{l \in L} B( l , \epsilon )$. The \emph{Ball Mapper graph} is obtained by assigning each ball $B(l,\epsilon)$ to  a vertex $v(l)$ of the graph, and by placing an edge between any two vertices whose corresponding balls jointly cover points from $X$. An illustration of the Ball Mapper algorithm, for the dataset from Figure~\ref{fig:idea_of_conventional_mapper}, is shown in  Figure~\ref{fig:Ball_mapper_example}.

\ifdefined\ifNotPictures
\begin{figure}[h]
 \centering
 \includegraphics[width=0.6\textwidth]{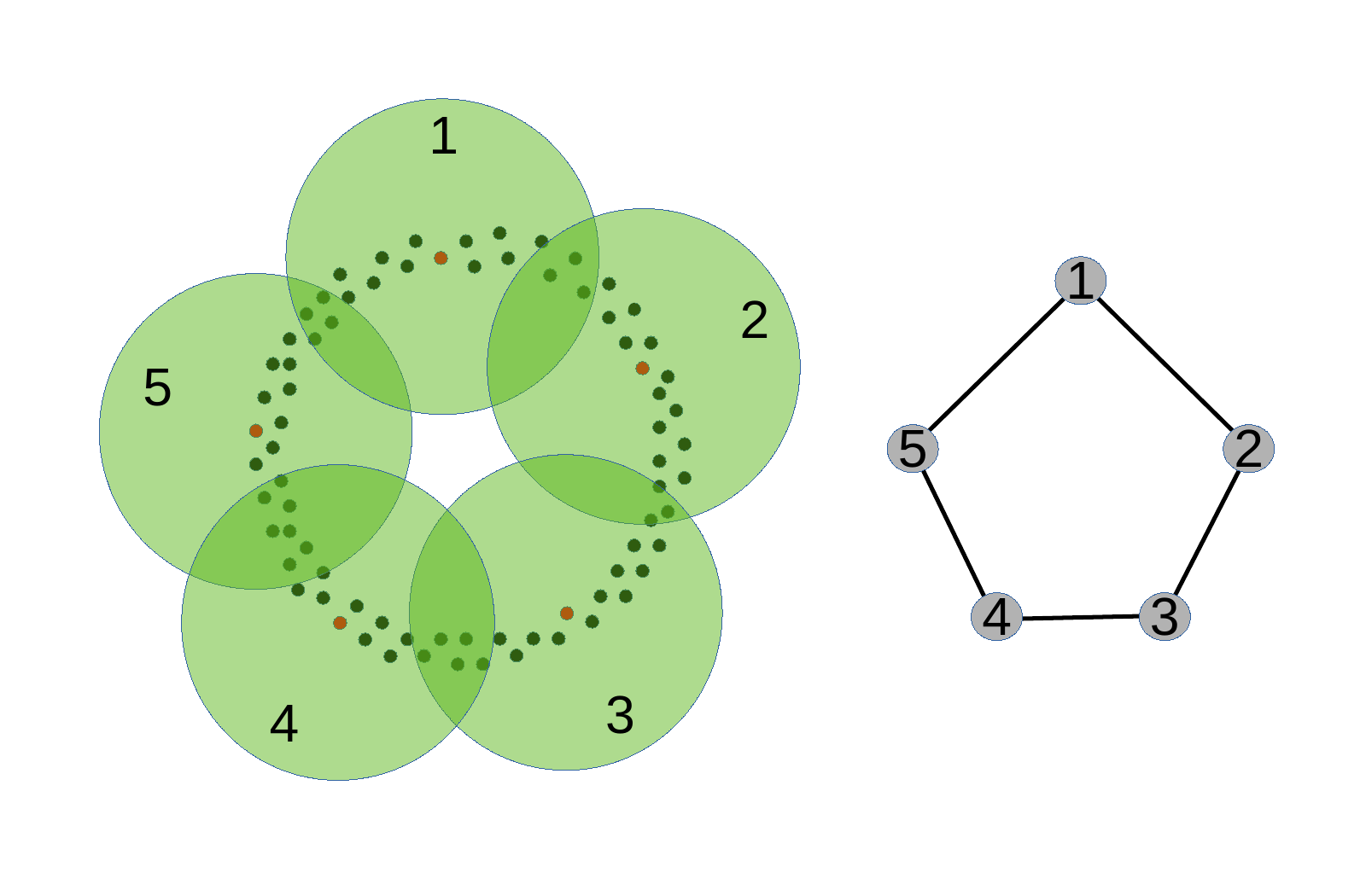}
\caption{Ball Mapper construction illustration. Using the point cloud $X$ in the Figure~\ref{fig:idea_of_conventional_mapper} we first select an $\epsilon$--net (red points in the left figure). Construct balls of radius $\epsilon$ centered in all the points of the net (green). The union of those balls provides a cover of our data space. Ball Mapper gives a one dimensional nerve of the obtained cover. In more detail, to each ball $B(A,\epsilon)$ we assign a vertex $v(A)$ in the abstract graph $G = (V,E)$. Two balls $B(N,\epsilon)$ and $B(K,\epsilon)$ that jointly cover some points from $X$ give rise to an edge between $v(N)$ and $v(K)$. The Ball Mapper graph of this cover is shown on the right with the corresponding balls and vertices labeled with the same number.}
\label{fig:Ball_mapper_example}
\end{figure}
\fi

\ifdefined\showOldText 
\OLD{
Ball Mapper graphs, analogous to Mapper graphs, can also be used to explore relations between different data descriptors. Given a function on the point cloud $g : X \rightarrow \mathbb{R}$, we define an \emph{induced function} $\hat{g} : G \rightarrow \mathbb{R}$ on the Ball Mapper graph $G$ in the following way;  for every landmark point $l$ the corresponding vertex $v(l)$ in the Ball Mapper graph $G$ gets the value determined by }
$$\hat{g}( v(l) ) =  \frac{ \displaystyle\sum_{x \in B(l,\epsilon) \cap X} g(x) }{|B(l,\epsilon) \cap X|}$$ 
\OLD{
The induced function $\hat{g}$ is subsequently visualized using a color scale. 
To confirm that the induced function $\hat{g}$ is a good approximation of the original function $g$ one can consider the \emph{variation of the values of $g$} on each element of the cover and compare it to the value of the induced function. If the variation is small relative to the function value $\hat{g}$ can be used for further data analysis both in case of Mapper and Ball Mapper. }
\fi
\NEW{
\subsection{Mapper and Ball Mapper: plotting scalar-valued functions}
Let $G$ be a Mapper or Ball Mapper graph. Its vertices cover collections of points of the input point cloud $X$. Therefore, given a function $g:X\rightarrow \mathbb{R}$  we define an \emph{induced function} $\hat{g}: G \rightarrow \mathbb{R}$, on $G$.
The value of the induced function $\hat{g}$ on a vertex $v(A) \in G$ covering $A \subset X$ is an average value of $g$ over all points in $A$: $\hat{g}( v(A) ) = \frac{ \sum_{x \in A} g(x) }{ |A| }$.}
Since the structure of Mapper or Ball Mapper graph $G$ reflects the shape of the input data, visualizing the induced function on $G$ using a color scale provides an insight into the variation of a scalar--valued function $g$ on $X$. 
\NEW{This procedure can be seen as a generalization to high dimensional samples of standards techniques, such as heatmaps, which allow for visualization of a scalar function on a compact subset of $\mathbb{R}^2$ by means of a color scale. 
To confirm that the induced function $\hat{g}$ is a good approximation of the original function $g$ one can consider the \emph{variation of the values of $g$} on each element of the cover and compare it to the value of the induced function. If the variation is small relative to the function value $\hat{g}$ can be used for further data analysis both in case of Mapper and Ball Mapper.
}

\section{New developments of Mapper algorithms}
\label{sec:new_developement_mapper}

\subsection{Equivariant Ball Mapper}
\label{sec:equivarinat_ball_mapper}

Suppose our data belongs to a metric space $(X,d)$ with \NEW{a finite} automorphism group of isometries $H$  acting on it. Although it is not necessary for the construction, in this paper we assume that the data is not noisy. Moreover, we assume that \NEW{for every sample point $x \in X$, the data contains all elements that are in its orbit,} i.e. for every $g \in H$, $g(x) \in X$.
\NEW{If that is not the case, the data should be properly augmented to ensure that the sample has all the symmetries imposed by the group action.} 
For example, see the integer--valued knot polynomials considered in  Section~\ref{sec:Data} and the Tic-Tac-Toe endgame configurations data in Section~\ref{sec:additional_equivariant_bm}.  In case that noise is present or that data sample might not accurately reflect expected isometries, one can include the orbit of each sampled point or points in the cloud closest to the points in orbits. 

Since $H$ is an automorphism group, for every $x,y \in X$ and every $g \in H$, $d(x,y) = d( g(x),g(y) )$. For every point $x \in X$, and isometry $g\in H$ the \emph{orbit} $\Omega_g(x)$ of $x$  contains the sequence of points $x,g(x),g^2(x),\ldots,g^{n}(x) \in X$, where $g^{n+1}(x) = x$ for some $n \in \mathbb{N}$. The relation \eqref{MirJ} for Jones polynomial data provides a dataset with a simple automorphism group generated by a permutation of the coordinates given by the exchange matrix.

Given a point cloud $X$ and an automorphism group $H$ acting on it, we modify the Ball Mapper algorithm in such a way that there is 
\ifdefined\showOldText \OLD{\emph{lift} or} 
\fi
an induced action of $H$ in the Ball Mapper graph $G$. 
The induced action is described in the following way: given a vertex $v(l) \in G$ consider $B(l,\epsilon) \cap X$, i.e. all the points in $X$ covered by a ball corresponding to $v(l) \in G$. 
For every isometry $g \in H$ we require the "covering" condition:
\begin{itemize}
    \item all the points in $g( B(l,\epsilon) \cap X )$ are covered by a ball $B( g(l) , \epsilon )$ and
    \item there are no other points in this ball.
\end{itemize}  
The vertex in $G$ corresponding to the ball $B( g(l) , \epsilon )$ is therefore denoted by $\hat{g}( v(l) )$. Note that such $\hat{g}$, induced by
$g$, is acting on an abstract graph, and therefore certain properties of $g$ will not be reflected in $\hat{g}$. An example of this construction is given in Figure~\ref{fig:eq_bm_8_rotated}.

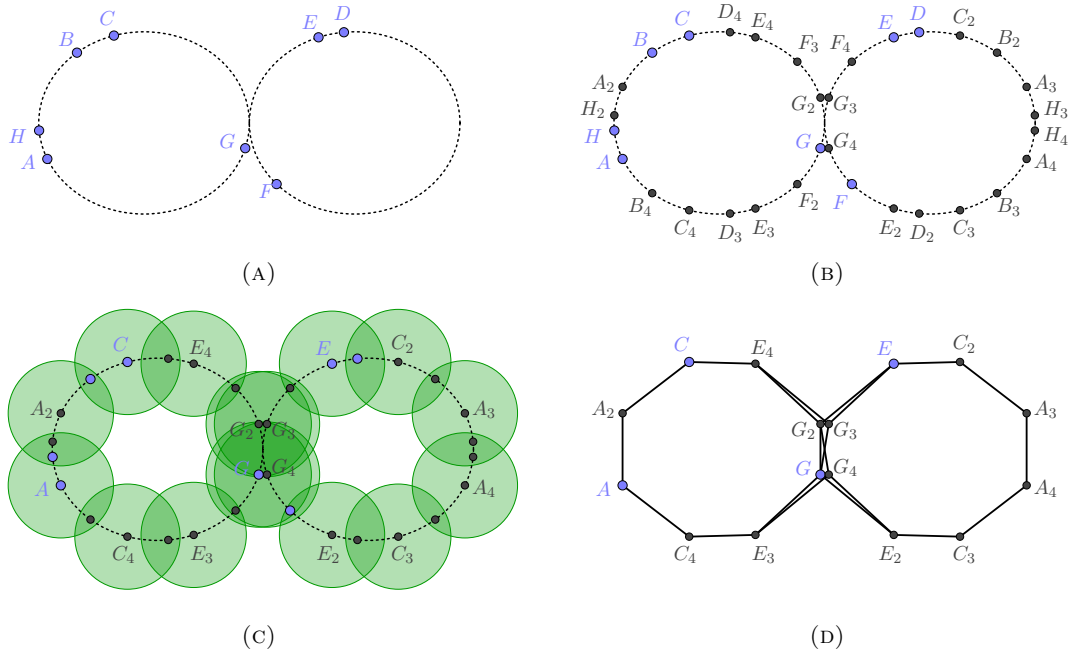
\begin{figure}
     \begin{subfigure}[]{0.45\textwidth}
         \scalebox{0.7}{\definecolor{xdxdff}{rgb}{0.49019607843137253,0.49019607843137253,1}
\begin{tikzpicture}[line cap=round,line join=round,>=triangle 45,x=1cm,y=1cm]

\clip(-5,-2.3) rectangle (5,2.3);

\draw [rotate around={0:(-2,0)},line width=0.8pt,dash pattern=on 1pt off 2pt] (-2,0) ellipse (2cm and 1.7320508075688772cm);
\draw [rotate around={0:(2,0)},line width=0.8pt,dash pattern=on 1pt off 2pt] (2,0) ellipse (2cm and 1.7320508075688772cm);

\draw [fill=xdxdff] (-3.836930294054438,-0.685029430819441) circle (2.5pt);
\draw[color=xdxdff] (-4.2,-0.75) node {$A$};
\draw [fill=xdxdff] (-3.273900626945473,1.3352463796996221) circle (2.5pt);
\draw[color=xdxdff] (-3.473379890246376,1.5922058739273686) node {$B$};
\draw [fill=xdxdff] (-2.57156380627509,1.6598150835310685) circle (2.5pt);
\draw[color=xdxdff] (-2.697599021325328,2) node {$C$};
\draw [fill=xdxdff] (1.7961821901392556,1.723033291984713) circle (2.5pt);
\draw[color=xdxdff] (1.7713358433043718,2.08389797394775) node {$D$};
\draw [fill=xdxdff] (1.3136826487660311,1.626876261442889) circle (2.5pt);
\draw[color=xdxdff] (1.1485258499452204,1.9200006072742897) node {$E$};
\draw [fill=xdxdff] (0.5178133080004113,-1.1629023852174727) circle (2.5pt);
\draw[color=xdxdff] (0.329039016577916,-1.3033142706371015) node {$F$};
\draw [fill=xdxdff] (-0.07850064402872259,-0.4804999154541802) circle (2.5pt);
\draw[color=xdxdff] (-0.41396237900843996,-0.3417830528194662) node {$G$};
\draw [fill=xdxdff] (-3.9929980876658053,-0.1448066535775655) circle (2.5pt);
\draw[color=xdxdff] (-4.4,-0.25) node {$H$};
\end{tikzpicture}}
         \caption{}
     \end{subfigure}
     \begin{subfigure}[]{0.45\textwidth}
         \scalebox{0.7}{\definecolor{uuuuuu}{rgb}{0.26666666666666666,0.26666666666666666,0.26666666666666666}
\definecolor{xdxdff}{rgb}{0.49019607843137253,0.49019607843137253,1}
\begin{tikzpicture}[line cap=round,line join=round,>=triangle 45,x=1cm,y=1cm]

\clip(-5,-2.3) rectangle (5,2.3);

\draw [rotate around={0:(-2,0)},line width=0.8pt,dash pattern=on 1pt off 2pt] (-2,0) ellipse (2cm and 1.7320508075688772cm);
\draw [rotate around={0:(2,0)},line width=0.8pt,dash pattern=on 1pt off 2pt] (2,0) ellipse (2cm and 1.7320508075688772cm);

\draw [fill=xdxdff] (-3.836930294054438,-0.685029430819441) circle (2.5pt);
\draw[color=xdxdff] (-4.2,-0.75) node {$A$};
\draw [fill=uuuuuu] (-3.836930294054438,0.685029430819441) circle (2pt);
\draw[color=uuuuuu] (-4.2,0.75) node {$A_{2}$};
\draw [fill=uuuuuu] (3.836930294054438,0.685029430819441) circle (2pt);
\draw[color=uuuuuu] (4.2,0.75) node {$A_{3}$};
\draw [fill=uuuuuu] (3.836930294054438,-0.685029430819441) circle (2pt);
\draw[color=uuuuuu] (4.2,-0.75) node {$A_{4}$};

\draw [fill=xdxdff] (-3.273900626945473,1.3352463796996221) circle (2.5pt);
\draw[color=xdxdff] (-3.5,1.6) node {$B$};
\draw [fill=uuuuuu] (-3.273900626945473,-1.3352463796996221) circle (2pt);
\draw[color=uuuuuu] (-3.5,-1.6) node {$B_{4}$};
\draw [fill=uuuuuu] (3.273900626945473,1.3352463796996221) circle (2pt);
\draw[color=uuuuuu] (3.5,1.6) node {$B_{2}$};
\draw [fill=xdxdff] (-2.57156380627509,1.6598150835310685) circle (2.5pt);
\draw[color=uuuuuu] (3.5,-1.6) node {$B_{3}$};
\draw [fill=xdxdff] (1.7961821901392556,1.723033291984713) circle (2.5pt);

\draw[color=xdxdff] (-2.697599021325328,2) node {$C$};
\draw [fill=uuuuuu] (2.57156380627509,1.6598150835310685) circle (2pt);
\draw[color=uuuuuu] (2.65638162334106,2) node {$C_{2}$};
\draw [fill=uuuuuu] (-2.57156380627509,-1.6598150835310685) circle (2pt);
\draw[color=uuuuuu] (-2.642966565767508,-2) node {$C_{4}$};
\draw [fill=uuuuuu] (2.57156380627509,-1.6598150835310685) circle (2pt);
\draw[color=uuuuuu] (2.65638162334106,-2) node {$C_{3}$};
\draw [fill=uuuuuu] (3.273900626945473,-1.3352463796996221) circle (2pt);

\draw[color=xdxdff] (1.7713358433043718,2.08389797394775) node {$D$};
\draw [fill=uuuuuu] (1.7961821901392556,-1.723033291984713) circle (2pt);
\draw[color=uuuuuu] (1.85,-2.1) node {$D_{2}$};
\draw [fill=uuuuuu] (-1.7961821901392556,-1.723033291984713) circle (2pt);
\draw[color=uuuuuu] (-1.7907002590655114,-2.1) node {$D_{3}$};
\draw [fill=uuuuuu] (-1.7961821901392556,1.723033291984713) circle (2pt);
\draw[color=uuuuuu] (-1.8016267501770753,2.094824465059314) node {$D_{4}$};
\draw [fill=xdxdff] (1.3136826487660311,1.626876261442889) circle (2.5pt);
\draw[color=xdxdff] (1.1485258499452204,1.92) node {$E$};
\draw [fill=uuuuuu] (1.3136826487660311,-1.626876261442889) circle (2pt);
\draw[color=uuuuuu] (1.246864269949297,-2) node {$E_{2}$};
\draw [fill=uuuuuu] (-1.3136826487660311,-1.626876261442889) circle (2pt);
\draw[color=uuuuuu] (-1.156963774594796,-2) node {$E_{3}$};
\draw [fill=uuuuuu] (-1.3136826487660311,1.626876261442889) circle (2pt);
\draw[color=uuuuuu] (-1.1788167568179242,1.92) node {$E_{4}$};
\draw [fill=xdxdff] (0.5178133080004113,-1.1629023852174727) circle (2.5pt);

\draw[color=xdxdff] (0.3,-1.5) node {$F$};
\draw [fill=uuuuuu] (-0.5178133080004113,-1.1629023852174727) circle (2pt);
\draw[color=uuuuuu] (-0.3,-1.5) node {$F_{2}$};
\draw [fill=uuuuuu] (-0.5178133080004113,1.1629023852174727) circle (2pt);
\draw[color=uuuuuu] (-0.3,1.5) node {$F_{3}$};
\draw [fill=uuuuuu] (0.5178133080004113,1.1629023852174727) circle (2pt);
\draw[color=uuuuuu] (0.3,1.5) node {$F_{4}$};
\draw [fill=xdxdff] (-0.07850064402872259,-0.4804999154541802) circle (2.5pt);

\draw[color=xdxdff] (-0.4,-0.35) node {$G$};
\draw [fill=uuuuuu] (-0.07850064402872259,0.4804999154541802) circle (2pt);
\draw[color=uuuuuu] (-0.4,0.35) node {$G_{2}$};
\draw [fill=uuuuuu] (0.07850064402872259,0.4804999154541802) circle (2pt);
\draw[color=uuuuuu] (0.4,0.35) node {$G_{3}$};
\draw [fill=uuuuuu] (0.07850064402872259,-0.4804999154541802) circle (2pt);
\draw[color=uuuuuu] (0.4,-0.35) node {$G_{4}$};

\draw [fill=xdxdff] (-3.9929980876658053,-0.1448066535775655) circle (2.5pt);
\draw[color=xdxdff] (-4.4,-0.25) node {$H$};
\draw [fill=uuuuuu] (3.9929980876658053,-0.1448066535775655) circle (2pt);
\draw[color=uuuuuu] (4.4,-0.25) node {$H_{4}$};
\draw [fill=uuuuuu] (-3.9929980876658053,0.1448066535775655) circle (2pt);
\draw[color=uuuuuu] (-4.4,0.25) node {$H_{2}$};
\draw [fill=uuuuuu] (3.9929980876658053,0.1448066535775655) circle (2pt);
\draw[color=uuuuuu] (4.4,0.25) node {$H_{3}$};
\end{tikzpicture}}
         \caption{}
     \end{subfigure}
     \\
     \begin{subfigure}[]{0.45\textwidth}
         \scalebox{0.7}{\definecolor{qqzzqq}{rgb}{0,0.6,0}
\definecolor{uuuuuu}{rgb}{0.26666666666666666,0.26666666666666666,0.26666666666666666}
\definecolor{xdxdff}{rgb}{0.49019607843137253,0.49019607843137253,1}
\begin{tikzpicture}[line cap=round,line join=round,>=triangle 45,x=1cm,y=1cm]

\clip(-5,-3) rectangle (5,3);

\draw [line width=0pt,color=qqzzqq,fill=qqzzqq,fill opacity=0.3] (-3.836930294054438,-0.685029430819441) circle (1cm);
\draw [line width=0pt,color=qqzzqq,fill=qqzzqq,fill opacity=0.3] (-3.836930294054438,0.685029430819441) circle (1cm);
\draw [line width=0pt,color=qqzzqq,fill=qqzzqq,fill opacity=0.3] (3.836930294054438,0.685029430819441) circle (1cm);
\draw [line width=0pt,color=qqzzqq,fill=qqzzqq,fill opacity=0.3] (3.836930294054438,-0.685029430819441) circle (1cm);
\draw [line width=0pt,color=qqzzqq,fill=qqzzqq,fill opacity=0.3] (-2.57156380627509,1.6598150835310685) circle (1cm);
\draw [line width=0pt,color=qqzzqq,fill=qqzzqq,fill opacity=0.3] (2.57156380627509,1.6598150835310685) circle (1cm);
\draw [line width=0pt,color=qqzzqq,fill=qqzzqq,fill opacity=0.3] (2.57156380627509,-1.6598150835310685) circle (1cm);
\draw [line width=0pt,color=qqzzqq,fill=qqzzqq,fill opacity=0.3] (-2.57156380627509,-1.6598150835310685) circle (1cm);
\draw [line width=0pt,color=qqzzqq,fill=qqzzqq,fill opacity=0.3] (1.3136826487660311,1.626876261442889) circle (1cm);
\draw [line width=0pt,color=qqzzqq,fill=qqzzqq,fill opacity=0.3] (1.3136826487660311,-1.626876261442889) circle (1cm);
\draw [line width=0pt,color=qqzzqq,fill=qqzzqq,fill opacity=0.3] (-1.3136826487660311,-1.626876261442889) circle (1cm);
\draw [line width=0pt,color=qqzzqq,fill=qqzzqq,fill opacity=0.3] (-1.3136826487660311,1.626876261442889) circle (1cm);
\draw [line width=0pt,color=qqzzqq,fill=qqzzqq,fill opacity=0.3] (-0.07850064402872259,-0.4804999154541802) circle (1cm);
\draw [line width=0pt,color=qqzzqq,fill=qqzzqq,fill opacity=0.3] (-0.07850064402872259,0.4804999154541802) circle (1cm);
\draw [line width=0pt,color=qqzzqq,fill=qqzzqq,fill opacity=0.3] (0.07850064402872259,0.4804999154541802) circle (1cm);
\draw [line width=0pt,color=qqzzqq,fill=qqzzqq,fill opacity=0.3] (0.07850064402872259,-0.4804999154541802) circle (1cm);

\draw [rotate around={0:(-2,0)},line width=0.8pt,dash pattern=on 1pt off 2pt] (-2,0) ellipse (2cm and 1.7320508075688772cm);
\draw [rotate around={0:(2,0)},line width=0.8pt,dash pattern=on 1pt off 2pt] (2,0) ellipse (2cm and 1.7320508075688772cm);

\draw [fill=xdxdff] (-3.836930294054438,-0.685029430819441) circle (2.5pt);
\draw[color=xdxdff] (-4.2,-0.75) node {$A$};
\draw [fill=uuuuuu] (-3.836930294054438,0.685029430819441) circle (2pt);
\draw[color=uuuuuu] (-4.2,0.75) node {$A_{2}$};
\draw [fill=uuuuuu] (3.836930294054438,0.685029430819441) circle (2pt);
\draw[color=uuuuuu] (4.2,0.75) node {$A_{3}$};
\draw [fill=uuuuuu] (3.836930294054438,-0.685029430819441) circle (2pt);
\draw[color=uuuuuu] (4.2,-0.75) node {$A_{4}$};

\draw [fill=xdxdff] (-3.273900626945473,1.3352463796996221) circle (2.5pt);
\draw [fill=uuuuuu] (-3.273900626945473,-1.3352463796996221) circle (2pt);
\draw [fill=uuuuuu] (3.273900626945473,1.3352463796996221) circle (2pt);
\draw [fill=xdxdff] (-2.57156380627509,1.6598150835310685) circle (2.5pt);
\draw[color=xdxdff] (-2.697599021325328,2) node {$C$};
\draw [fill=uuuuuu] (2.57156380627509,1.6598150835310685) circle (2pt);
\draw[color=uuuuuu] (2.65638162334106,2) node {$C_{2}$};
\draw [fill=uuuuuu] (-2.57156380627509,-1.6598150835310685) circle (2pt);
\draw[color=uuuuuu] (-2.642966565767508,-2) node {$C_{4}$};
\draw [fill=uuuuuu] (2.57156380627509,-1.6598150835310685) circle (2pt);
\draw[color=uuuuuu] (2.65638162334106,-2) node {$C_{3}$};
\draw [fill=uuuuuu] (3.273900626945473,-1.3352463796996221) circle (2pt);
\draw [fill=xdxdff] (1.7961821901392556,1.723033291984713) circle (2.5pt);
\draw [fill=uuuuuu] (1.7961821901392556,-1.723033291984713) circle (2pt);
\draw [fill=uuuuuu] (-1.7961821901392556,1.723033291984713) circle (2pt);
\draw [fill=xdxdff] (1.3136826487660311,1.626876261442889) circle (2.5pt);
\draw[color=xdxdff] (1.1485258499452204,1.92) node {$E$};
\draw [fill=uuuuuu] (1.3136826487660311,-1.626876261442889) circle (2pt);
\draw[color=uuuuuu] (1.246864269949297,-2) node {$E_{2}$};
\draw [fill=uuuuuu] (-1.3136826487660311,-1.626876261442889) circle (2pt);
\draw[color=uuuuuu] (-1.156963774594796,-2) node {$E_{3}$};
\draw [fill=uuuuuu] (-1.3136826487660311,1.626876261442889) circle (2pt);
\draw[color=uuuuuu] (-1.1788167568179242,1.92) node {$E_{4}$};
\draw [fill=xdxdff] (0.5178133080004113,-1.1629023852174727) circle (2.5pt);
\draw [fill=uuuuuu] (-0.5178133080004113,-1.1629023852174727) circle (2pt);
\draw [fill=uuuuuu] (-0.5178133080004113,1.1629023852174727) circle (2pt);
\draw [fill=uuuuuu] (0.5178133080004113,1.1629023852174727) circle (2pt);
\draw [fill=xdxdff] (-0.07850064402872259,-0.4804999154541802) circle (2.5pt);

\draw[color=xdxdff] (-0.4,-0.35) node {$G$};
\draw [fill=uuuuuu] (-0.07850064402872259,0.4804999154541802) circle (2pt);
\draw[color=uuuuuu] (-0.4,0.35) node {$G_{2}$};
\draw [fill=uuuuuu] (0.07850064402872259,0.4804999154541802) circle (2pt);
\draw[color=uuuuuu] (0.4,0.35) node {$G_{3}$};
\draw [fill=uuuuuu] (0.07850064402872259,-0.4804999154541802) circle (2pt);
\draw[color=uuuuuu] (0.4,-0.35) node {$G_{4}$};

\draw [fill=uuuuuu] (-1.7961821901392556,-1.723033291984713) circle (2pt);

\draw [fill=xdxdff] (-3.9929980876658053,-0.1448066535775655) circle (2.5pt);
\draw [fill=uuuuuu] (3.9929980876658053,-0.1448066535775655) circle (2pt);
\draw [fill=uuuuuu] (-3.9929980876658053,0.1448066535775655) circle (2pt);
\draw [fill=uuuuuu] (3.9929980876658053,0.1448066535775655) circle (2pt);

\end{tikzpicture}}
         \caption{}
     \end{subfigure}
     \begin{subfigure}[]{0.45\textwidth}
         \scalebox{0.7}{\definecolor{uuuuuu}{rgb}{0.26666666666666666,0.26666666666666666,0.26666666666666666}
\definecolor{xdxdff}{rgb}{0.49019607843137253,0.49019607843137253,1}
\begin{tikzpicture}[line cap=round,line join=round,>=triangle 45,x=1cm,y=1cm]

\clip(-5,-3) rectangle (5,3);

\draw [line width=1pt] (-3.836930294054438,0.685029430819441)-- (-2.57156380627509,1.6598150835310685);
\draw [line width=1pt] (-2.57156380627509,1.6598150835310685)-- (-1.3136826487660311,1.626876261442889);
\draw [line width=1pt] (-1.3136826487660311,1.626876261442889)-- (-0.07850064402872259,0.4804999154541802);
\draw [line width=1pt] (-0.07850064402872259,0.4804999154541802)-- (-0.07850064402872259,-0.4804999154541802);
\draw [line width=1pt] (-0.07850064402872259,-0.4804999154541802)-- (-1.3136826487660311,-1.626876261442889);
\draw [line width=1pt] (-1.3136826487660311,-1.626876261442889)-- (-2.57156380627509,-1.6598150835310685);
\draw [line width=1pt] (-2.57156380627509,-1.6598150835310685)-- (-3.836930294054438,-0.685029430819441);
\draw [line width=1pt] (-3.836930294054438,-0.685029430819441)-- (-3.836930294054438,0.685029430819441);
\draw [line width=1pt] (-0.07850064402872259,-0.4804999154541802)-- (0.07850064402872259,0.4804999154541802);
\draw [line width=1pt] (-0.07850064402872259,0.4804999154541802)-- (0.07850064402872259,-0.4804999154541802);
\draw [line width=1pt] (-0.07850064402872259,-0.4804999154541802)-- (0.07850064402872259,-0.4804999154541802);
\draw [line width=1pt] (0.07850064402872259,0.4804999154541802)-- (-0.07850064402872259,0.4804999154541802);
\draw [line width=1pt] (0.07850064402872259,-0.4804999154541802)-- (-1.3136826487660311,-1.626876261442889);
\draw [line width=1pt] (0.07850064402872259,-0.4804999154541802)-- (1.3136826487660311,-1.626876261442889);
\draw [line width=1pt] (1.3136826487660311,-1.626876261442889)-- (2.57156380627509,-1.6598150835310685);
\draw [line width=1pt] (2.57156380627509,-1.6598150835310685)-- (3.836930294054438,-0.685029430819441);
\draw [line width=1pt] (3.836930294054438,-0.685029430819441)-- (3.836930294054438,0.685029430819441);
\draw [line width=1pt] (3.836930294054438,0.685029430819441)-- (2.57156380627509,1.6598150835310685);
\draw [line width=1pt] (2.57156380627509,1.6598150835310685)-- (1.3136826487660311,1.626876261442889);
\draw [line width=1pt] (1.3136826487660311,1.626876261442889)-- (0.07850064402872259,0.4804999154541802);
\draw [line width=1pt] (1.3136826487660311,1.626876261442889)-- (-0.07850064402872259,0.4804999154541802);
\draw [line width=1pt] (-1.3136826487660311,1.626876261442889)-- (0.07850064402872259,0.4804999154541802);
\draw [line width=1pt] (-0.07850064402872259,-0.4804999154541802)-- (1.3136826487660311,-1.626876261442889);

\draw [fill=xdxdff] (-3.836930294054438,-0.685029430819441) circle (2.5pt);
\draw[color=xdxdff] (-4.2,-0.75) node {$A$};
\draw [fill=uuuuuu] (-3.836930294054438,0.685029430819441) circle (2pt);
\draw[color=uuuuuu] (-4.2,0.75) node {$A_{2}$};
\draw [fill=uuuuuu] (3.836930294054438,0.685029430819441) circle (2pt);
\draw[color=uuuuuu] (4.2,0.75) node {$A_{3}$};
\draw [fill=uuuuuu] (3.836930294054438,-0.685029430819441) circle (2pt);
\draw[color=uuuuuu] (4.2,-0.75) node {$A_{4}$};

\draw [fill=xdxdff] (-2.57156380627509,1.6598150835310685) circle (2.5pt);
\draw[color=xdxdff] (-2.7136720084034573,2) node {$C$};
\draw [fill=uuuuuu] (2.57156380627509,1.6598150835310685) circle (2pt);
\draw[color=uuuuuu] (2.65638162334106,2) node {$C_{2}$};
\draw [fill=uuuuuu] (-2.57156380627509,-1.6598150835310685) circle (2pt);
\draw[color=uuuuuu] (-2.642966565767508,-2) node {$C_{4}$};
\draw [fill=uuuuuu] (2.57156380627509,-1.6598150835310685) circle (2pt);
\draw[color=uuuuuu] (2.65638162334106,-2) node {$C_{3}$};

\draw [fill=xdxdff] (1.3136826487660311,1.626876261442889) circle (2.5pt);
\draw[color=xdxdff] (1.1389290845362432,1.92) node {$E$};
\draw [fill=uuuuuu] (1.3136826487660311,-1.626876261442889) circle (2pt);
\draw[color=uuuuuu] (1.246864269949297,-2) node {$E_{2}$};
\draw [fill=uuuuuu] (-1.3136826487660311,-1.626876261442889) circle (2pt);
\draw[color=uuuuuu] (-1.156963774594796,-2) node {$E_{3}$};
\draw [fill=uuuuuu] (-1.3136826487660311,1.626876261442889) circle (2pt);
\draw[color=uuuuuu] (-1.1788167568179242,1.92) node {$E_{4}$};
\draw [fill=xdxdff] (-0.07850064402872259,-0.4804999154541802) circle (2.5pt);

\draw[color=xdxdff] (-0.4,-0.35) node {$G$};
\draw [fill=uuuuuu] (-0.07850064402872259,0.4804999154541802) circle (2pt);
\draw[color=uuuuuu] (-0.4,0.35) node {$G_{2}$};
\draw [fill=uuuuuu] (0.07850064402872259,0.4804999154541802) circle (2pt);
\draw[color=uuuuuu] (0.4,0.35) node {$G_{3}$};
\draw [fill=uuuuuu] (0.07850064402872259,-0.4804999154541802) circle (2pt);
\draw[color=uuuuuu] (0.4,-0.35) node {$G_{4}$};
\end{tikzpicture}}
         \caption{}
     \end{subfigure}

    \caption{Example of the Equivariant Ball Mapper construction on a point cloud $X$ sampled from a wedge of two circles (A)
    with a symmetry group determined by the reflections on the horizontal and vertical axis. 
    The point cloud enriched by each point's orbit is shown in (B). Selected landmarks $A, C, E$ and $G$, and all the points in their orbits that are selected as well are shown in (C) with a covering consisting of 16 balls. The resulting equivariant Ball Mapper graph is depicted in (D).}
    \label{fig:eq_bm_8_rotated}
\end{figure}

To ensure the "covering" condition described above is satisfied, the procedure of selection of $\epsilon$-net $L \subset X$ is adjusted. In the Ball Mapper implementations $L$ is chosen in the greedy way presented in Algorithm~\ref{alg:greedy_epsilon_net}.  
The  main idea is to add the whole orbit  $\Omega(x)=\{ g(x) \}_{g \in H}$ to the constructed set of landmark points together with the any added point $x$. This idea is formalized in the Algorithm~\ref{alg:greedy_equivariant_epsilon_net}, \NEW{which adjusts  Algorithm~\ref{alg:greedy_epsilon_net} so that the obtained $\epsilon$-net $L$ is invariant under the action of $H$.}
\smallskip
\begin{center}
\begin{algorithm}[H]
 \KwData{Point cloud $X$, $\epsilon > 0$, group $H$ acting on  $X$}
  \KwResult{$L \subset X$, a $H$-equivariant $\epsilon$-net}
  $L = \emptyset$\\
  \While{There exists a point $x \in X$ farther than $\epsilon$ from any point in $L$}
  {
     $L = L \cup \Omega(x)$ ;
  }
  \Return L 
  \label{alg:greedy_equivariant_epsilon_net}
  \caption{Equivariant greedy $\epsilon$-net.}
 \end{algorithm}
 \end{center}

Since each $g \in H$ is an isometry, for every $l \in L$, if $y \in B(l,\epsilon) \cap X$,  $g(y) \in B( g(l) , \epsilon )\cap X$. As each $B(l,\epsilon) \cap X$ corresponds to a vertex in the Ball Mapper graph, the action of $H$ on $X$ is induced to the vertices of the Ball Mapper graph as described above.  

\NEW{In Algorithm~\ref{alg:greedy_equivariant_epsilon_net}, the \emph{while} loop iterates once through the points in $X$. The time needed to locate all points at a distance less or equal than $\epsilon$ to a given point is bounded by $|X|$, hence the overall complexity of the presented procedure is bounded by $|X|^2$.}
 
This equivariant construction is Ball Mapper specific and it is not clear that the analogous construction for  Mapper is possible. For example, the equivariant Mapper construction requires the lens function  to be either invariant with respect to the group action (i.e. points in the same orbit should obtain the same value of the function), or to map points from the input point cloud to different cover elements in such a way that there is an induced group acting on the cover. Adjusting the lens function to satisfy either of this requirements is non-trivial and also requires the clustering algorithm to be performed in such a way that there is an appropriate group action induced on the obtained clusters. These obstacles to obtaining equivariant Mapper emphasize the advantage of Ball Mapper with respect to equivariance and justify our choice.

\subsection{Mapper on Ball Mapper}\label{sec:BonBM}
This section provides a new construction of a Mapper graph, as described in Section~\ref{sec:conventional_mapper}. It uses the  output of the Ball Mapper algorithm, described in Section~\ref{sec:ball_mapper}
to cover the range of the lens function.

Recall that the Mapper construction is based on a interval cover 
of the range of the lens function $f : X \rightarrow \mathbb{R}^n$. Typically $n = 1$ or another small positive integer, as the range of $f$ needs to be covered with a collection of overlapping cubes, defined as a product of intervals. There is antagonism between wanting large $n$ when the lens function is more likely to preserve essential information about the point cloud, and the fact that having $k$ intervals in each of $n$ directions requires  $k^n$ cover elements which is not computationally feasible for large values of $n$.
\ifdefined\showOldText 
\OLD{To overcome this obstacle we modify Mapper construction to use the overlapping, adaptive cover of the point cloud $f(X)$ obtained by the Ball Mapper. }\fi

\NEW{To overcome this obstacle we propose the following \emph{Mapper on Ball Mapper (MoBM)} construction to leverage the overlapping, adaptive cover of the point cloud obtained from Ball Mapper.}
\ifdefined\showOldText 
\OLD{More generally, suppose that given two point clouds denoted by  $X$ and $Y$ that are related via a function $p: X \rightarrow Y$.  Our aim is to analyze the first point cloud $X$. In this setting we propose the following pipeline, illustrated in Figure~\ref{fig:mapper_on_ball_mapper}: }
\fi
\NEW{This algorithm applies to the more general setting with two point clouds $X$, $Y$ and a relation  $p: X\rightarrow Y$.
MoBM is formalized by the following pseudocode:} \\

\NEW{
\begin{algorithm}[H]
\KwData{Point clouds $X, Y$, a relation $p : X \rightarrow Y$, $\epsilon > 0$, a clustering algorithm $C$}
\KwResult{MoBM graph representing $X$} 
Let $\{B( l,\epsilon )\}_{l \in L}$ be the balls in a Ball Mapper of $p(X)$ constructed for the radius $\epsilon$; \\
Apply $C$ to each $\{ p^{-1}( B( l,\epsilon )) \cap  X  \}_{l \in L}$ to obtain an overlapping cover of $X$;\\
MoBM = 1-dimensional nerve of such cover; \\
\Return MoBM 
\label{alg:MoBM}
\caption{Mapper on Ball Mapper.}
\end{algorithm}
\medskip
}

This construction, illustrated in Figure~\ref{fig:mapper_on_ball_mapper}, requires only two parameters: the radius $\epsilon$ and the choice of the clustering algorithm. There is no need to define the number of intervals or the overlapping percentage as in the conventional Mapper algorithm; the covering of $Y$, being the range of the lens, is completely determined by the Ball Mapper graph. However, varying the selection of landmark points can lead to potentially different Ball Mapper graphs. Selection of different landmarks can be obtained by permuting the points of $Y$ and re-running the algorithm. \NEW{The time complexity of MoBM is bounded by the time required to run both Mapper and Ball Mapper.}

\ifdefined\ifNotPictures
\begin{figure}[h!]
  \centering
\resizebox{0.85\textwidth}{!}{\input{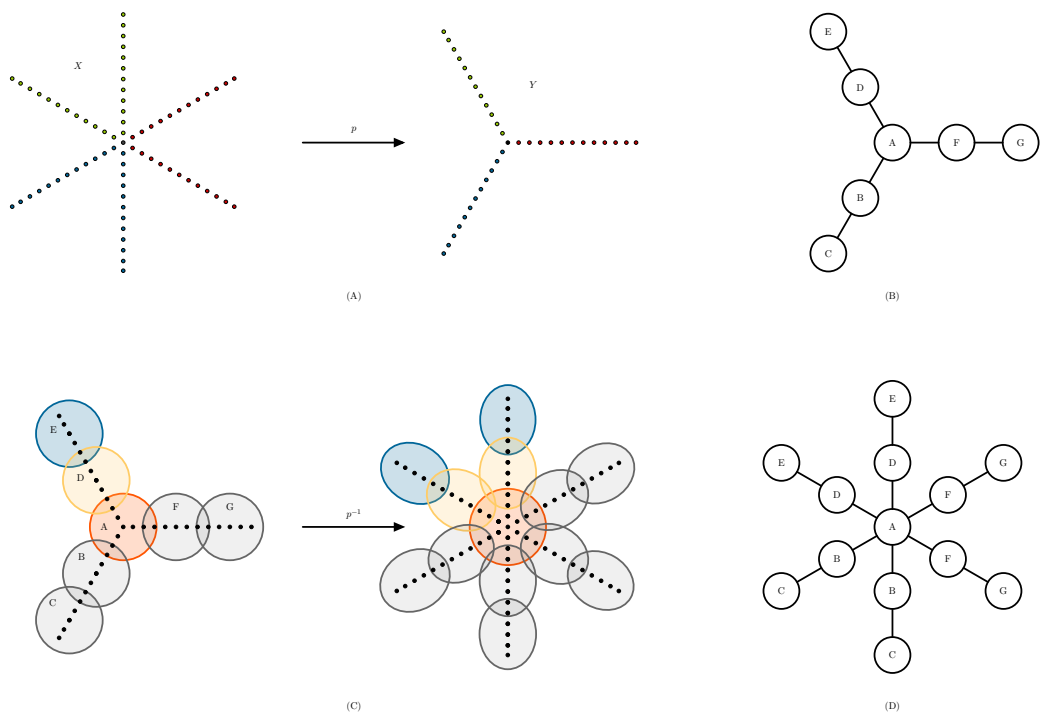}}
\caption{The Mapper on Ball Mapper construction: an illustration. The point cloud $X$ is mapped via $p$ to the point cloud $Y$ (A) where the map is indicated by corresponding colors.
     The Ball Mapper graph for $Y$ (B) is used to obtain a cover of $Y$ (C, left). Such cover is then pulled-back to obtain a covering of $X$  (C, right).
     The resulting Mapper on Ball Mapper graph for $X$ is shown in (D), where, similarly to Figure~\ref{fig:idea_of_conventional_mapper}, node labels indicate the originating covering ball in $Y$.}
     \label{fig:mapper_on_ball_mapper}
 \end{figure}
\fi
%
%
%
% %
% %
% %
% %
\subsection{MappingMappers}
\label{sec:relation_between_mapper}

The standard construction of Mapper or Ball Mapper for a point cloud $X$ with a metric $d : X \times X \rightarrow \mathbb{R}^{\geq 0}$ provides a model of the point cloud $X$. If the point cloud $X$ is equipped with a function $f : X \rightarrow \mathbb{R}$, an induced function $\hat{f}$ can be defined on the Mapper or Ball Mapper graph $G$ as explained in  Sections~\ref{sec:conventional_mapper} and~\ref{sec:ball_mapper}. This way, Mapper and Ball Mapper graphs can be used to visualize functions $f : \mathbb{R}^n \rightarrow \mathbb{R}$, where $f$ is defined on $X \subset \mathbb{R}^n$.

Consider a more general question of using Mapper and Ball Mapper to visualize functions $f : \mathbb{R}^n \rightarrow \mathbb{R}^m$ for larger values of $m$ and $n$. Assume that we are given $X \subset \mathbb{R}^n$ and $Y \subset \mathbb{R}^m$ and a relation $f \subset X \times Y$ (note that a function is a particular case of such a relation).
In this instance we focus on the Ball Mapper--based construction. A construction for Mapper is analogous. 
In the first step, let us build Ball Mapper graphs $G(X)$ and $G(Y)$ corresponding to point clouds $X$ and $Y$ and denote with $V(X)$ and $V(Y)$ the corresponding vertex sets. Given a relation $f \subset X \times Y$ assigning points from $X$ to the points from $Y$, define a map $\tilde{f} : V(X) \rightarrow [0,1]^{V(Y)}$ in the following way.
For every $v$ in $V(X)$, corresponding to a ball $B(l_X,\epsilon_X)$ compute $f( B(l_X,\epsilon_X) \cap X ) \subset Y$. For every vertex $w$ in $V(Y)$ corresponding to a ball $B(l_Y,\epsilon_Y)$, compute $\frac{| B(l_Y,\epsilon_Y) \cap f( B(l_X,\epsilon_X) \cap X ) |}{|B(l_Y,\epsilon_Y)  \cap Y|}$. This fraction indicate the percentage of points in $B(l_Y,\epsilon_Y) \cap Y$ that are in the image of the points covered by the vertex $v$ in $G(X)$. When computed for every vertex in $G(X)$, this fraction gives us $|V(X)|$ different coloring functions on $G(Y)$ indicating where the image of each vertex $v$ is mapped.

This construction works analogously for arbitrarily large unions of vertices of the graph $G(X)$. The procedure described above is automatized and a reliable interface can be found in  \cite{zenodo}. It allows to see which regions of $G(Y)$ corresponds to chosen region of $G(X)$. By doing so, a visualization of the map $f : X \rightarrow Y$ is obtained. A simplified example of the procedure is given at the Figure~\ref{fig:function}.

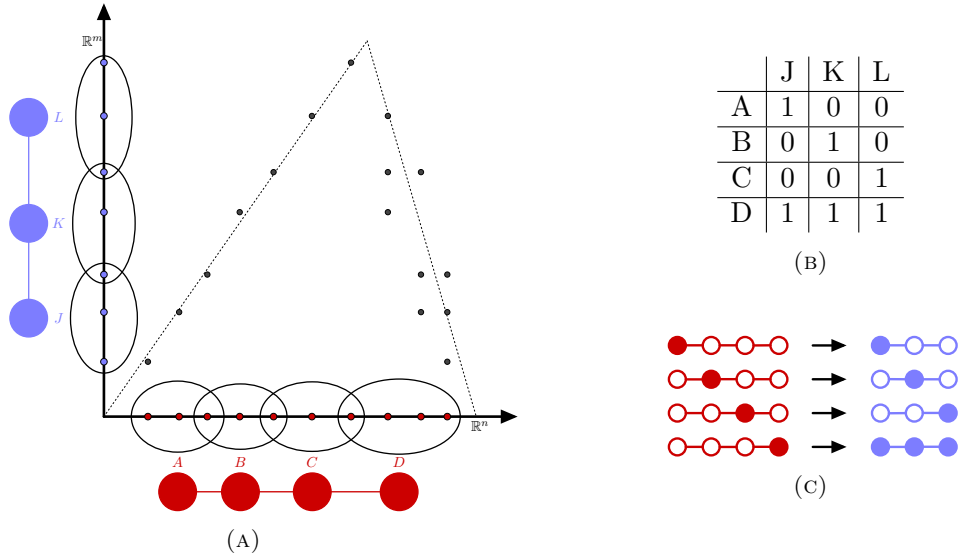
\begin{figure}[!htb]
\centering
\begin{tabular}[t]{cc}
\begin{subfigure}{0.4\textwidth}
    \centering
    \scalebox{0.5}{\definecolor{uuuuuu}{rgb}{0.26666666666666666,0.26666666666666666,0.26666666666666666}
\definecolor{ccqqqq}{rgb}{0.8,0,0}
\definecolor{xdxdff}{rgb}{0.49019607843137253,0.49019607843137253,1}
\begin{tikzpicture}[line cap=round,line join=round,>=triangle 45,x=1cm,y=1cm]

\draw [line width=1pt,color=ccqqqq,fill=ccqqqq,fill opacity=1] (1.96,-2) circle (0.5cm);
\draw [line width=1pt,color=ccqqqq,fill=ccqqqq,fill opacity=1] (3.63,-2) circle (0.5cm);
\draw [line width=1pt,color=ccqqqq,fill=ccqqqq,fill opacity=1] (5.54,-2) circle (0.5cm);
\draw [line width=1pt,color=ccqqqq,fill=ccqqqq,fill opacity=1] (7.85,-2) circle (0.5cm);
\draw [line width=1pt,color=xdxdff,fill=xdxdff,fill opacity=1] (-2,7.962727272727269) circle (0.5cm);
\draw [line width=1pt,color=xdxdff,fill=xdxdff,fill opacity=1] (-2,5.142727272727269) circle (0.5cm);
\draw [line width=1pt,color=xdxdff,fill=xdxdff,fill opacity=1] (-2,2.62) circle (0.5cm);
\draw [->,line width=2pt] (0,0) -- (0,11);
\draw (10,-0.25) node {$\mathbb{R}^n$};
\draw [->,line width=2pt] (0,0) -- (11,0);
\draw (-0.3,10) node {$\mathbb{R}^m$};
\draw [line width=0.4pt,dash pattern=on 1pt off 2pt] (0,0)-- (7,10);
\draw [line width=0.4pt,dash pattern=on 1pt off 2pt] (7,10)-- (9.9,0);
\draw [rotate around={90:(0,2.62)},line width=1pt] (0,2.62) ellipse (1.4616200120739025cm and 0.889231724408723cm);
\draw [rotate around={90:(0,5.142727272727266)},line width=1pt] (0,5.142727272727266) ellipse (1.5957801663129818cm and 0.8303545744820077cm);
\draw [rotate around={90:(0,7.962727272727269)},line width=1pt] (0,7.962727272727269) ellipse (1.6380357935431857cm and 0.7480089967879913cm);
\draw [rotate around={0:(1.96,0)},line width=1pt] (1.96,0) ellipse (1.2319294825127454cm and 0.94527786913908cm);
\draw [rotate around={0:(3.63,0)},line width=1pt] (3.63,0) ellipse (1.237170515508981cm and 0.8695923668275599cm);
\draw [rotate around={0:(5.54,0)},line width=1pt] (5.54,0) ellipse (1.3845619557684994cm and 0.9252631027775257cm);
\draw [rotate around={0:(7.85,0)},line width=1pt] (7.85,0) ellipse (1.6240664762153556cm and 0.9995958779259518cm);
\draw [line width=1pt,color=ccqqqq] (1.96,-2)-- (3.63,-2);
\draw [line width=1pt,color=ccqqqq] (3.63,-2)-- (5.54,-2);
\draw [line width=1pt,color=ccqqqq] (5.54,-2)-- (7.85,-2);
\draw [line width=1pt,color=xdxdff] (-2,7.962727272727269)-- (-2,5.142727272727269);
\draw [line width=1pt,color=xdxdff] (-2,5.142727272727269)-- (-2,2.62);

%WARNING: PGF/Tikz and Gnuplot don't support implicit curves
%Rather try PSTricks export
%Cannot draw null

\draw [fill=xdxdff] (0,3.78) circle (2.5pt);
\draw [fill=xdxdff] (0,1.46) circle (2.5pt);
\draw [fill=xdxdff] (0,6.505454545454538) circle (2.5pt);
\draw [fill=xdxdff] (0,9.42) circle (2.5pt);
\draw [fill=ccqqqq] (1.17,0) circle (2.5pt);
\draw [fill=ccqqqq] (2.75,0) circle (2.5pt);
\draw [fill=ccqqqq] (4.51,0) circle (2.5pt);
\draw [fill=ccqqqq] (6.57,0) circle (2.5pt);
\draw [fill=ccqqqq] (9.13,0) circle (2.5pt);
\draw [fill=uuuuuu] (9.13,1.46) circle (2pt);
\draw [fill=uuuuuu] (9.13,3.78) circle (2pt);
\draw [fill=uuuuuu] (1.17,1.46) circle (2pt);
\draw [fill=uuuuuu] (2.75,3.78) circle (2pt);
\draw [fill=uuuuuu] (6.57,9.42) circle (2pt);
\draw [fill=uuuuuu] (4.51,6.505454545454538) circle (2pt);
\draw[color=ccqqqq] (1.96,-1.2) node {$A$};
\draw[color=ccqqqq] (3.63,-1.2) node {$B$};
\draw[color=ccqqqq] (5.54,-1.2) node {$C$};
\draw[color=ccqqqq] (7.85,-1.2) node {$D$};
\draw[color=xdxdff] (-1.2,7.962727272727269) node {$L$};
\draw[color=xdxdff] (-1.2,5.142727272727269) node {$K$};
\draw[color=xdxdff] (-1.2,2.62) node {$J$};
\draw [fill=ccqqqq] (2,0) circle (2.5pt);
\draw [fill=xdxdff] (0,2.78) circle (2.5pt);
\draw [fill=uuuuuu] (2,2.78) circle (2pt);
\draw [fill=ccqqqq] (3.61,0) circle (2.5pt);
\draw [fill=xdxdff] (0,5.44) circle (2.5pt);
\draw [fill=uuuuuu] (3.61,5.44) circle (2pt);
\draw [fill=ccqqqq] (5.53,0) circle (2.5pt);
\draw [fill=xdxdff] (0,8) circle (2.5pt);
\draw [fill=uuuuuu] (5.53,8) circle (2pt);
\draw [fill=ccqqqq] (8.43,0) circle (2.5pt);
\draw [fill=uuuuuu] (9.13,2.78) circle (2pt);
\draw [fill=uuuuuu] (8.43,2.78) circle (2pt);
\draw [fill=uuuuuu] (8.43,3.78) circle (2pt);
\draw [fill=uuuuuu] (8.43,6.505454545454538) circle (2pt);
\draw [fill=ccqqqq] (7.55,0) circle (2.5pt);
\draw [fill=uuuuuu] (7.55,6.505454545454538) circle (2pt);
\draw [fill=uuuuuu] (7.55,5.44) circle (2pt);
\draw [fill=uuuuuu] (7.55,8) circle (2pt);
\end{tikzpicture}}
    \caption{}
\end{subfigure}
    &
\begin{tabular}{c}% if you add [t], than sub images are pushed down
\begin{subfigure}[t]{0.4\textwidth}
\centering
    \begin{tabular}{c|c|c|c}
        & J & K & L \\ 
        \hline
        A & 1 & 0 & 0 \\
        \hline
        B & 0 & 1 & 0 \\
        \hline
        C & 0 & 0 & 1 \\
        \hline
        D & 1 & 1 & 1
    \end{tabular}
    \caption{}
\end{subfigure}\\
\bigskip \\
\begin{subfigure}[t]{0.45\textwidth}
\centering
\scalebox{0.45}{\definecolor{xdxdff}{rgb}{0.49019607843137253,0.49019607843137253,1}
\definecolor{ccqqqq}{rgb}{0.8,0,0}

\begin{tikzpicture}[line cap=round,line join=round,>=triangle 45,x=1cm,y=1cm]
\draw [line width=2pt,color=ccqqqq,fill=ccqqqq,fill opacity=1] (0,0) circle (0.25cm);
\draw [line width=2pt,color=xdxdff,fill=xdxdff,fill opacity=1] (6,0) circle (0.25cm);
\draw [line width=2pt,color=ccqqqq,fill=ccqqqq,fill opacity=1] (1,-1) circle (0.25cm);
\draw [line width=2pt,color=xdxdff,fill=xdxdff,fill opacity=1] (7,-1) circle (0.25cm);
\draw [line width=2pt,color=ccqqqq,fill=ccqqqq,fill opacity=1] (2,-2) circle (0.25cm);
\draw [line width=2pt,color=xdxdff,fill=xdxdff,fill opacity=1] (8,-2) circle (0.25cm);
\draw [line width=2pt,color=ccqqqq,fill=ccqqqq,fill opacity=1] (3,-3) circle (0.25cm);
\draw [line width=2pt,color=xdxdff,fill=xdxdff,fill opacity=1] (6,-3) circle (0.25cm);
\draw [line width=2pt,color=xdxdff,fill=xdxdff,fill opacity=1] (7,-3) circle (0.25cm);
\draw [line width=2pt,color=xdxdff,fill=xdxdff,fill opacity=1] (8,-3) circle (0.25cm);
\draw [line width=2pt,color=ccqqqq] (1,0) circle (0.25cm);
\draw [line width=2pt,color=ccqqqq] (2,0) circle (0.25cm);
\draw [line width=2pt,color=ccqqqq] (3,0) circle (0.25cm);
\draw [->,line width=2pt] (4,0) -- (5,0);
\draw [line width=2pt,color=xdxdff] (7,0) circle (0.25cm);
\draw [line width=2pt,color=xdxdff] (8,0) circle (0.25cm);
\draw [line width=2pt,color=ccqqqq] (0,-1) circle (0.25cm);
\draw [line width=2pt,color=ccqqqq] (2,-1) circle (0.25cm);
\draw [line width=2pt,color=ccqqqq] (3,-1) circle (0.25cm);
\draw [line width=2pt,color=xdxdff] (6,-1) circle (0.25cm);
\draw [line width=2pt,color=xdxdff] (8,-1) circle (0.25cm);
\draw [line width=2pt,color=ccqqqq] (0,-2) circle (0.25cm);
\draw [line width=2pt,color=ccqqqq] (1,-2) circle (0.25cm);
\draw [line width=2pt,color=ccqqqq] (3,-2) circle (0.25cm);
\draw [line width=2pt,color=xdxdff] (6,-2) circle (0.25cm);
\draw [line width=2pt,color=xdxdff] (7,-2) circle (0.25cm);
\draw [line width=2pt,color=ccqqqq] (0,-3) circle (0.25cm);
\draw [line width=2pt,color=ccqqqq] (1,-3) circle (0.25cm);
\draw [line width=2pt,color=ccqqqq] (2,-3) circle (0.25cm);
\draw [->,line width=2pt] (4,-1) -- (5,-1);
\draw [->,line width=2pt] (4,-2) -- (5,-2);
\draw [->,line width=2pt] (4,-3) -- (5,-3);
\draw [line width=2pt,color=ccqqqq] (0.24997564645839412,0.0034894380217946673)-- (0.7500275036600406,0.0037082440007300136);
\draw [line width=2pt,color=ccqqqq] (0.24997648783673734,-1.0034286336651081)-- (0.7500280865427957,-1.0037473300553834);
\draw [line width=2pt,color=ccqqqq] (0.25,-2)-- (0.75,-2);
\draw [line width=2pt,color=ccqqqq] (0.249997969124202,-2.9989923126603206)-- (0.7500020774783458,-2.9989808165734444);
\draw [line width=2pt,color=ccqqqq] (1.2499733163106619,0.0036525515259597744)-- (1.7500277000194173,0.003721457028844713);
\draw [line width=2pt,color=ccqqqq] (1.2499853092305167,-0.9972898027889234)-- (1.7500152362589958,-0.9972399461319941);
\draw [line width=2pt,color=ccqqqq] (1.2499984015038506,-1.999105994675911)-- (1.7500018366043781,-1.9990417209091371);
\draw [line width=2pt,color=ccqqqq] (1.2499957300887856,-3.001461142489635)-- (1.7500047253717246,-3.0015370958113095);
\draw [line width=2pt,color=ccqqqq] (2.2499756773955046,-2.0034872210510106)-- (2.750026957020913,-2.003671210124113);
\draw [line width=2pt,color=ccqqqq] (2.249982267452754,-2.997022423942218)-- (2.7500182919189324,-2.9969758265803828);
\draw [line width=2pt,color=ccqqqq] (2.249997294521149,-0.9988369298792085)-- (2.7500029146861085,-0.9987927988738717);
\draw [line width=2pt,color=ccqqqq] (2.25,0)-- (2.750031510964379,0.003969192518482048);
\draw [line width=2pt,color=xdxdff] (6.249973634035112,0.0036307419737197084)-- (6.750026953113758,0.003670944075942525);
\draw [line width=2pt,color=xdxdff] (7.249973582661386,0.0036342772915005836)-- (7.750027959207548,0.0037388262940405603);
\draw [line width=2pt,color=xdxdff] (6.249969983938162,-0.9961260962898976)-- (6.75,-1);
\draw [line width=2pt,color=xdxdff] (7.249998520896956,-1.0008599705426424)-- (7.750001506140831,-1.0008677949912965);
\draw [line width=2pt,color=xdxdff] (6.2499448932454,-2.0052488418289993)-- (6.750000719404223,-2.0005997512769036);
\draw [line width=2pt,color=xdxdff] (7.249966031370834,-1.9958789369435093)-- (7.750034589447477,-1.9958414512978035);
\draw [line width=2pt,color=xdxdff] (6.249999732488965,-3.0003657259162466)-- (6.75000027335672,-3.000369700264054);
\draw [line width=2pt,color=xdxdff] (7.2499646318906175,-2.9957949073983716)-- (7.750034700689633,-2.995834770032957);
\end{tikzpicture}}
\caption{}
\end{subfigure}
\end{tabular}\\
\end{tabular}
\caption{Consider a relation $f$ between $X \subset \mathbb{R}^n$  located on the x-axis 
and $Y \subset \mathbb{R}^m$  on the y-axis,  represented by the black dots in the plane,  see (A). Next, create Ball Mapper graphs on both domain $X$ and co-domain $Y$. Set $X$ is covered by four balls corresponding to the red vertices in the Ball Mapper graph, while $Y$,  with the three balls corresponding to the blue vertices in the Ball Mapper graph. Based on the relation, the points in $X$ covered by $A$ are mapped into points covered by $J$ in $Y$. Similarly, $B$ into $K$, $C$ into $L$.  Lastly, the points in $D$ in $X$ are mapped to points in $J$, $K$ and $L$.
In addition, when we track the proportions of points in each of the ball $J$, $K$ and $L$ that are reached by points in $A$, $B$, $C$ or $D$. Such proportions are indicated in the matrix (B). Each row in this matrix provides a  coloring function on the image Ball Mapper graph (C). This idea generalizes and provides a way to visualize maps between point clouds in high dimensional spaces.}
\label{fig:function}
\end{figure}

\section{Applications: knot theory}
%Inspiration: knot theory}
\label{sec:KT_intro}

Analyzing and visualizing data from knot theory is the main motivation and inspiration for the development of algorithms presented in Sections~\ref{sec:equivarinat_ball_mapper}, \ref{sec:relation_between_mapper} and \ref{sec:BonBM}.  This section discusses the basis of knot theory and the data used in the our analysis, followed by the results obtained using Equivariant Ball Mapper,  Mapper on Ball Mapper, and MappingMappers.

% %
% %
% %
% %
% %======================
\subsection{Knot theory: a brief introduction}
\label{sec:KnotIntro}

A knot is a class of embeddings of $S^1$ into $\mathbb{R}^3$ up to ambient isotopy \cite{ livingston1993knot, lickorish2012introduction, jablan2007linknot}. Knots are hard to distinguish and their classification and tabulation \cite{hoste1998first, hoste2005enumeration,KA, knotinfo, BBurton18} solicits techniques from a range of mathematical disciplines.
Consequently, a number of \emph{knot invariants}  have been introduced in an attempt to compare and classify knots. A knot invariant should be thought of as a quantity assigned to each knot  such that if two knots are the same (isotopic), the values assigned to them are the same. The most common knot invariants are integers, one- or two-variable polynomials, groups, etc.  In this paper we focus on the following knot invariants and their relations:
\begin{itemize}
    \item Numerical: minimal crossing number, signature $\sigma(K)$ \cite{kauffman1976signature} defined as a signature of a matrix obtained using a Seifert surface.
    \item Polynomial: Alexander $\Delta(K)(t)$ \cite{alexander1928topological}, Jones $J(K)(q))$ \cite{jones1997polynomial}, $P(K)(a,z))$ HOMFLYPT \cite{freyd1985new}.
\end{itemize} 

Since knot invariants often rely on advanced algebraic, geometric,  and combinatorial topology, in lieu of definitions we provide references  \cite{livingston1993knot, rolfsen2003knots, lickorish2012introduction, jablan2007linknot}, and key insights sufficient for analyzing these knot invariants and their relations. 
Note that Alexander and Jones polynomials are $1$-variable polynomials, while HOMFLYPT polynomial is of two variables. Moreover, HOMFLYPT  is more general than both the Alexander and Jones polynomial. HOMFLYPT specializes to the Jones polynomial \cite{jones1997polynomial} by substituting $a=t^{-1}$ and $z=\sqrt{t}-\frac{1}{\sqrt{t}}$, and to the Alexander polynomial \cite{alexander1928topological} by substituting $a=1$ and $z=\sqrt{t}-\frac{1}{\sqrt{t}}$.

\begin{table}[h!]
 \centering
 \begin{tabular}{|c|c|c|c|}
     \hline
           Invariant & Unknot &  Trefoil & Data vector \\      
           \midrule
         Alexander  & $1$ &  $t^{-1}-1+t$ & (0,1,-1,1,0)\\
         \hline
         Jones & $1$ &  $q+q^3-q^4$& (0,0,0,0,0,1,0,1,-1)\\ \hline
         HOMFLYPT & $1$ & $-a^4 + 2a^2 + a^2z^2$& (0,0,0,0,0,0,2,0,-1,\textbf{0,0,0,0,0,0,1,0,0})\\ 
    \hline
     \end{tabular}
     \caption{Values of several knot polynomials for the unknot and the trefoil, and the corresponding data vector. Note that for the 2-variable HOMFLYPT polynomial the coefficient matrix is flattened into a vector: the variable $z$ corresponding to rows, $a$ to columns.}
     \label{tab:trefoil}
\end{table}

Note that it is common for knot tables to contain only one knot from each mirror pair $(K,mir(K))$, since if a knot $K$ is achiral it coincides with its mirror $mir(K)$, {and many invariants either do not distinguish between mirrors or have a straightforward relation between the two values.} The Alexander polynomial does not distinguish between the knot and its mirror, and the Jones polynomial \cite{lickorish2012introduction,jones1997polynomial} satisfies  the following relation: 
\begin{equation}\label{MirJ}
 J(mir(K))(q) = J(K)(q^{-1}). 
\end{equation}

The signature of a knot and its mirror have opposite signs  ~\cite{lickorish2012introduction}, \begin{equation}\label{MirS}\sigma (mir(K))=-\sigma(K).
\end{equation} 
 and for the HOMFLYPT polynomials of mirror knots the following relation holds:
\begin{equation}\label{MirHOM} P(mir(K))(a,z)=
P(K)(a^{-1},z)
\end{equation}

%$$$$$$$$$$$$$$$$$$$$$$$$$$$$$$$$$$$$$$$$$$$$$$$$$

\subsection{Knot invariants as  point clouds}\label{sec:Data} The use of big data techniques \cite{hughes2016neural,jejjala2019deep,ward2018using,levitt2019big} is warranted by the result of Ernst and Summers showing that 
 the number of knots with a given number of crossings grows exponentially~\cite{ernst1987growth}. The point clouds are obtained in the way introduced and described in \cite{levitt2019big}, where each knot is  represented by a vector of coefficients of a knot polynomial such as the Alexander, Jones or HOMFLYPT polynomial. 
The datasets we consider were created by J.S. Levitt \cite{levitt2019big}, and preprocessed by D. Gurnari. The data is freely available at \cite{zenodo} and includes Alexander and Jones polynomials, together with numerical knot invariants  like minimal crossing number and signature for all \NEW{$9755329$ knots up to $17$ crossings. HOMFLYPT polynomials are provided for all $313231$ knots up to $15$ crossings}.

Following \cite{levitt2019big}, 
given a finite collection of knots $\mathcal{K}$, we construct a point cloud $\mathcal{I(K)}$ corresponding to the coefficients of the one-variable  polynomial invariant $\mathcal{I}$, in the following way:
\begin{itemize}
     \item[Step 1] Given a knot $K \in \mathcal{K}$ and its single variable polynomial $I(K)$, extract a vector of the coefficients.
       \item[Step 2] Compute the minimal and maximal powers $min_{t}$, $max_{t}$ of the variable denoted by $t$  among all knots in $\mathcal{K}$. Then the maximal length of all considered vectors is $d=max_{t}-min_{t}+1$.
       \item[Step 3] Add zeros on both sides of each vector of coefficients to obtain a vector $I(K)_v \in \mathbb{R}^{d}$ to ensure a correct alignment of corresponding powers.
 \end{itemize}

Note that in this way all vectors are of the same length determined by the overall minimum and maximum exponent, and,  the coefficients of a given power are in the same position in the vector for all the considered knots. 
In case of a two-variable polynomial, such as the HOMFLYPT polynomial, we apply Steps 1-3, as described above in the case of one-variable polynomial,  to both variables. In this way we first obtain a matrix padded with zeros, and then create a corresponding vector by linearizing this matrix (concatenating its rows). HOMFLYPT data belongs to $\mathbb{R}^d$ where  $d=(max_{a}-min_{a}+1)(max_{z}-min_{z}+1)$ where $a$ and $z$ stand for the two variables in the corresponding polynomial. \NEW{Examples of the coefficient vectors obtained from both the unknot and the trefoil are presented in Table~\ref{tab:trefoil}.}

Unlike some databases, we choose to consider knots and their mirrors although that increases the dimension of the point cloud coefficient vectors. The Jones coefficient vector of its mirror is obtained by reversing the original vector, see  \eqref{MirJ}, and therefore the point cloud admits a symmetry given by the exchange matrix. The HOMFLYPT coefficient matrix of the mirror knot can be obtained by reversing the columns of the original matrix  \eqref{MirHOM}. Hence, in the case of the Jones polynomial $min_t(mir(K))=-max_t(K)$ and $max_t(mir(K))=-min_t(mir(K))$ as a consequence of relation \eqref{MirJ} and the point cloud belongs to $\mathbb{R}^d$ where $d=2max(|max_t|, |min_t|) +1$. 

\NEW{The size of the obtained tables of polynomial coefficients are the following: $9755329$ rows $\times$ $17$ columns for Alexander, $19510658$ rows $\times$ $51$ columns for Jones and $626462$ rows $\times$ $152$ columns for HOMFLYPT. Note that the number of rows in Jones and HOMFLYPT data is double the number of prime knots since  mirrors are also included.}

\ifdefined\ifNotPictures
\begin{figure}[h!]
     \centering
     \begin{subfigure}[b]{0.3\textwidth}
     \includegraphics[width=\textwidth]{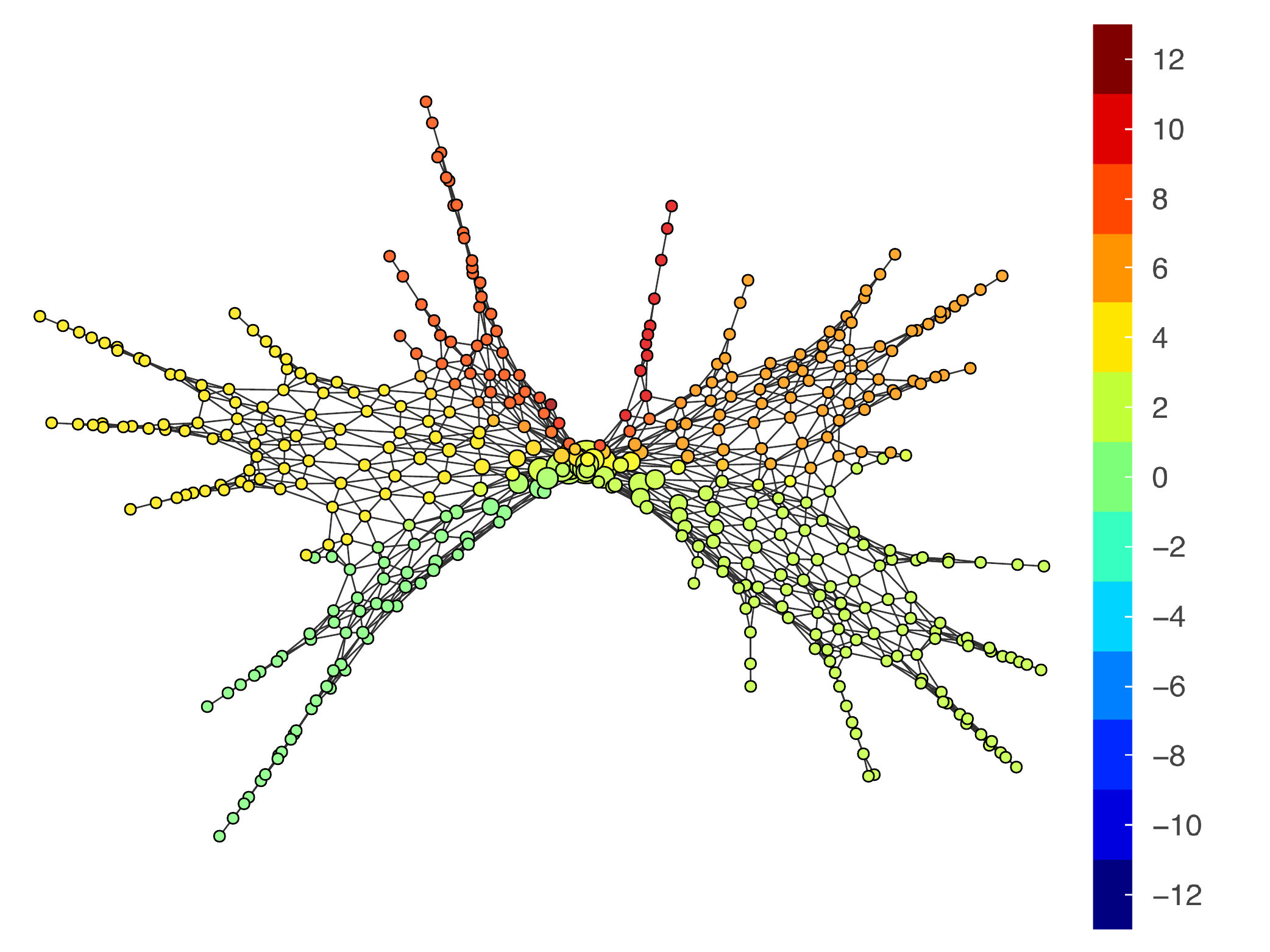}
     \caption{}%BM with no mirrors}
     \end{subfigure}
     \begin{subfigure}[b]{0.3\textwidth}
     \includegraphics[width=\textwidth]{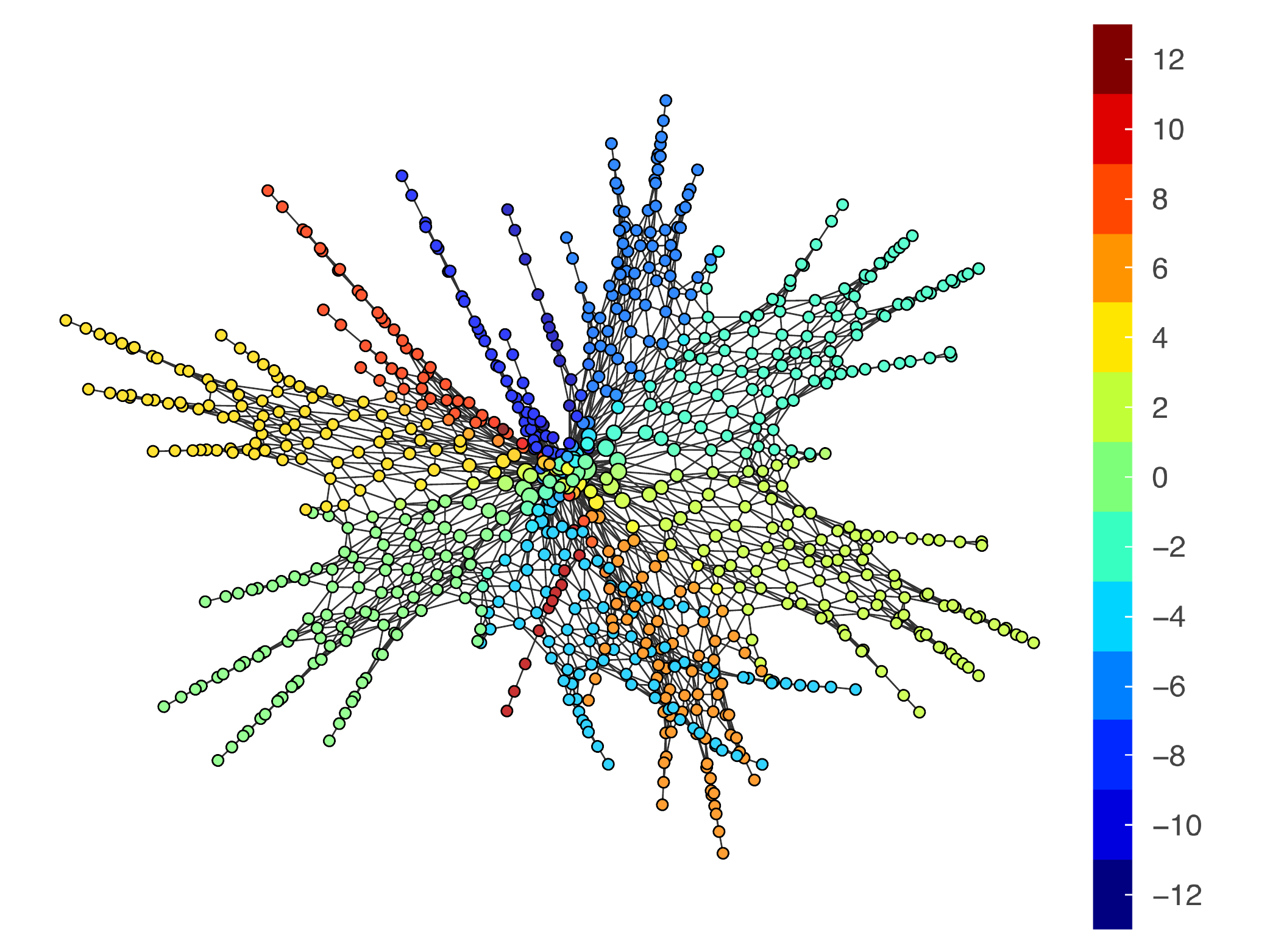}
     \caption{}%BM with mirrors}
     \end{subfigure}
     \begin{subfigure}[b]{0.3\textwidth}
     \includegraphics[width=\textwidth]{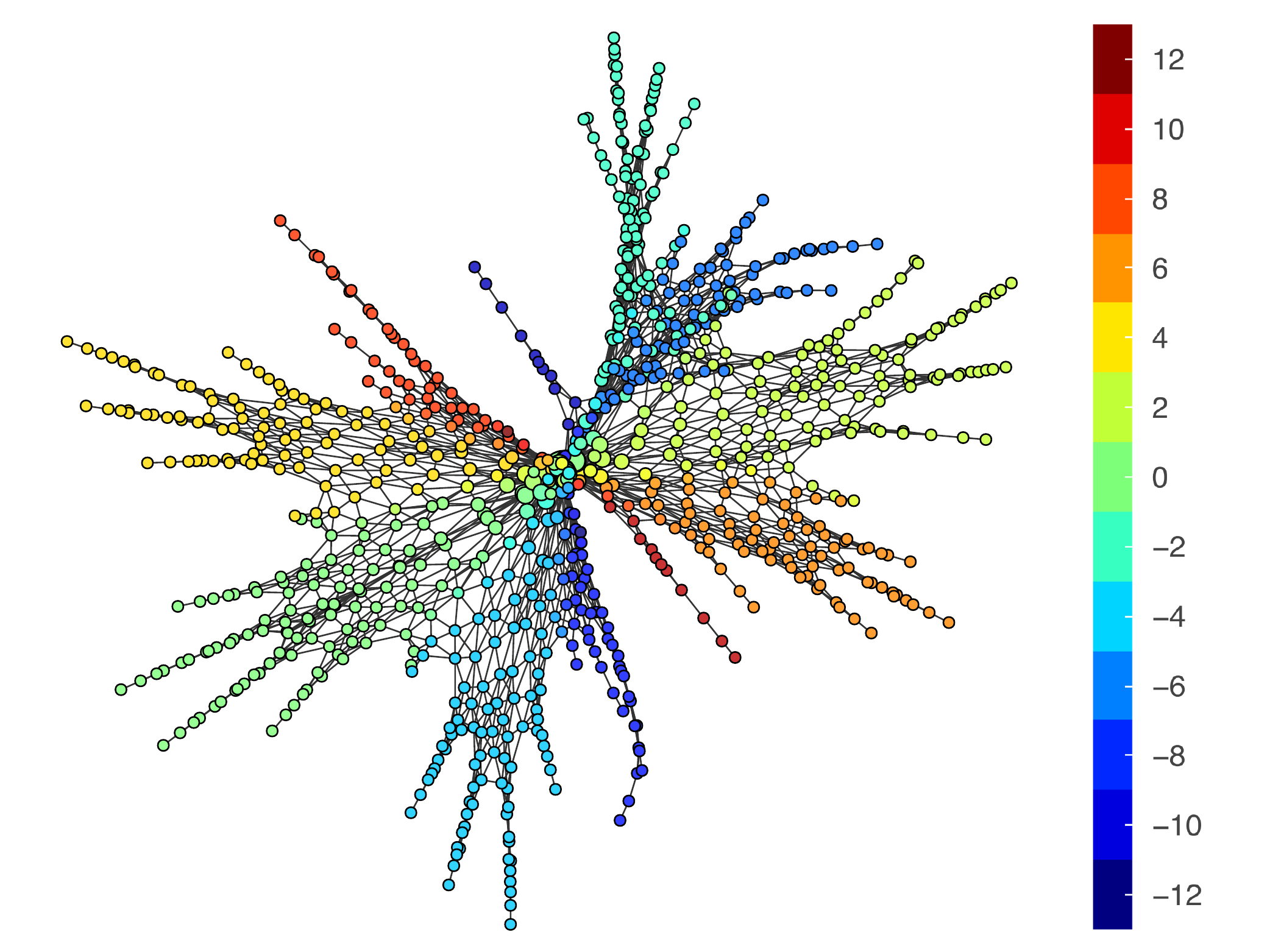}
     \caption{}
     \end{subfigure}
     \caption{Ball Mapper applied to Jones polynomial data  of knots up to 17 crossings with just one choice of a mirror (A), knots and their mirrors with standard Ball Mapper (B) and Equivariant Ball Mapper (C).  Color reflects  the average signature of knots in each cluster. Note that the graph in (C) is symmetric, although this fact is not accurately represented in this image  due to the chosen graph plotting subroutine.}
     \label{fig:Jones_mirors_and_symmetries}
\end{figure}
\fi
%
%
%$$$$$$$$$$$$$$$$$$$$$$$$$$$$$$$$$$$$$$$$$$$$$$$$$

\subsection{Ball Mapper: structure of knot polynomial data}\label{sec:KnotDataBM}

In this section we apply standard and Equivariant Ball Mapper to data obtained from Jones, Alexander and HOMFLYPT polynomials for all knots up to $17$ crossings. The choice to use Ball Mapper is natural, as there is no obvious lens function for the Mapper construction. 

\ifdefined\ifNotPictures
\begin{figure}[h!]
      \centering
      \begin{subfigure}[b]{0.3\textwidth}
          \centering
          \includegraphics[width=\textwidth]{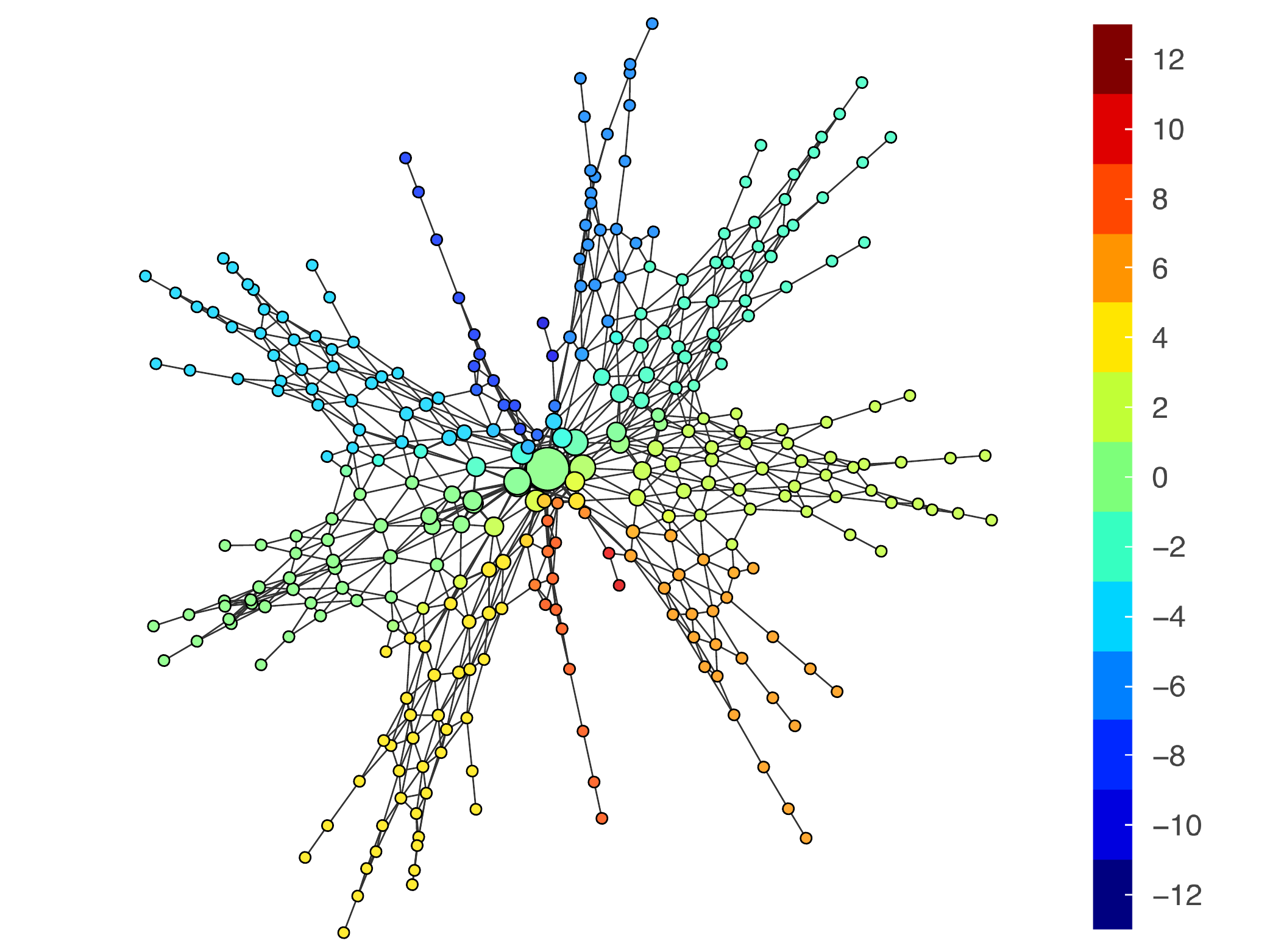}
          \caption{}
          %Up to 15,  $\epsilon=30$ with $826$ nodes}
      \end{subfigure}
      \hfill
      \begin{subfigure}[b]{0.3\textwidth}
          \centering
          \includegraphics[width=\textwidth]{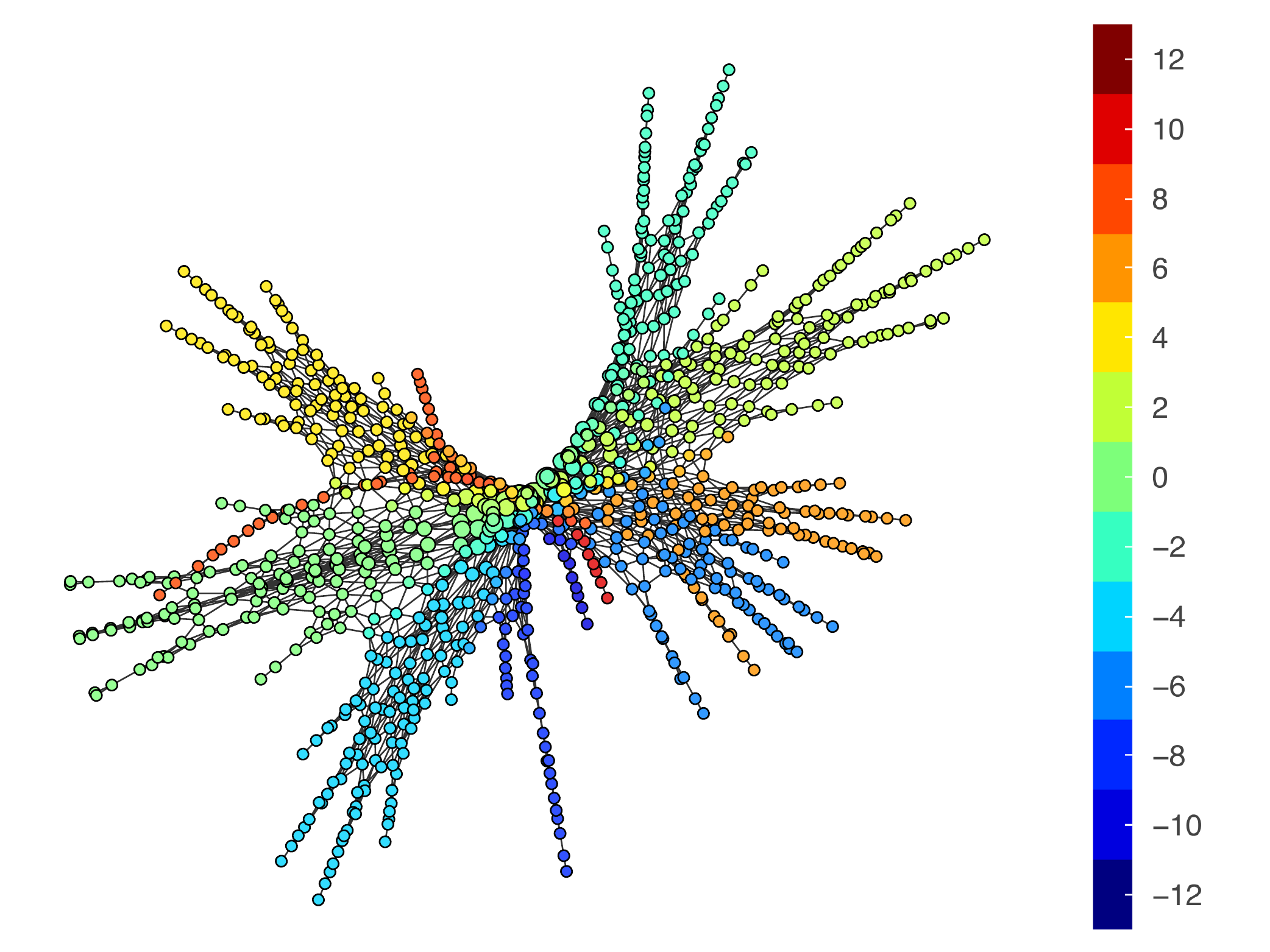}
          \caption{}
      \end{subfigure}
      \hfill
      \begin{subfigure}[b]{0.3\textwidth}
          \centering
          \includegraphics[width=\textwidth]{jones_upto_17/jones_upto_17_100_signature-eps-converted-to.pdf}
          \caption{}
          %up to 17 $\epsilon=100$ with $890$ nodes.}
      \end{subfigure}
     \caption{Stability of  \NEW{Equivariant} Ball Mapper for Jones polynomial data with respect to the crossing number filtration: EqBM graphs of knots up to $15$ crossings $\epsilon=30$ with $826$ nodes (A), $16$ with $\epsilon=50$ with $1008$ nodes (B), $17$ crossings with $\epsilon=100$ with $890$ nodes (C).}
     \label{fig:stability_Jones_construction}
 \end{figure}
 \fi

\NEW{Since our data contains knots and their mirror images, the Jones polynomial data cloud admits a symmetry generated by the permutation of the coordinates for those knots which are not identical to their mirrors (see \eqref{MirJ}). }
Figure \ref{fig:Jones_mirors_and_symmetries} shows Ball Mapper graphs with just one of the mirrors included (A), knots and their mirrors with the standard Ball Mapper (B) and the equivariant Ball Mapper construction from Section \ref{sec:equivarinat_ball_mapper} in \ref{fig:Jones_mirors_and_symmetries}(C). The symmetry of the data is preserved by the equivariant Ball Mapper (C): for each flare there is an identical one with opposite signature, as a consequence of relation \eqref{MirS}.  

The structure of the BM graph is stable with respect to the filtration by the number of crossings, as illustrated in Figure \ref{fig:stability_Jones_construction}, in addition to stability across the choice of parameter/radius shown as shown in Figure \ref{fig:stability_Jones_construction_epsilon}. This observation is important since sampling knot data is a hard problem and it is known that knots with lower crossing number do not provide a sample representative of the space of knots \cite{levitt2019big}.

\ifdefined\ifNotPictures
\begin{figure}[h!]
      \centering
      \begin{subfigure}[b]{0.3\textwidth}
          \centering
          \includegraphics[width=\textwidth]{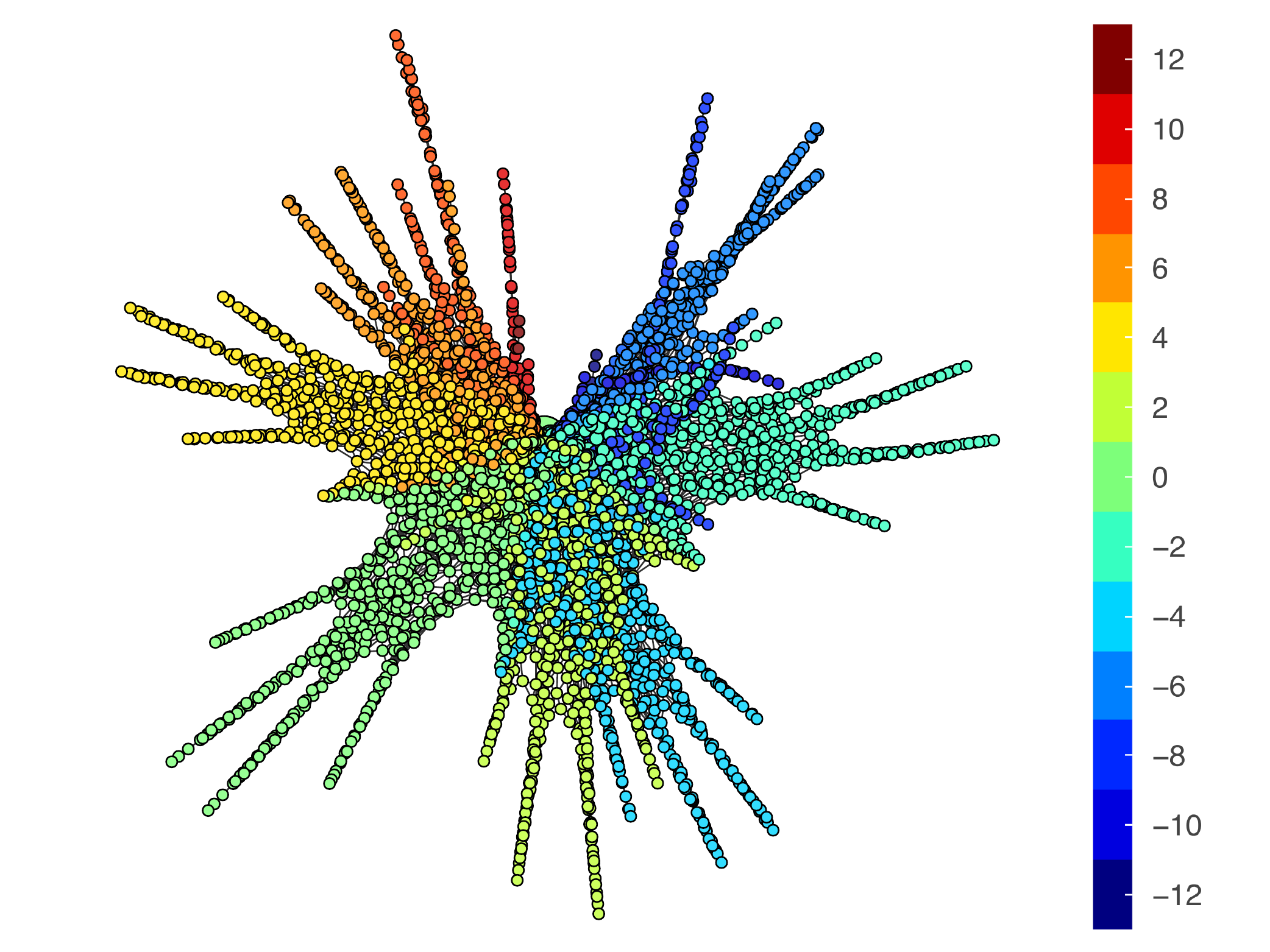}
          \caption{}
          %Up to 15,  $\epsilon=30$ with $826$ nodes}
      \end{subfigure}
     \hfill
      \begin{subfigure}[b]{0.3\textwidth}
          \centering
          \includegraphics[width=\textwidth]{jones_upto_17/jones_upto_17_100_signature-eps-converted-to.pdf}
          \caption{}
      \end{subfigure}
     \hfill
      \begin{subfigure}[b]{0.3\textwidth}
          \centering
          \includegraphics[width=\textwidth]{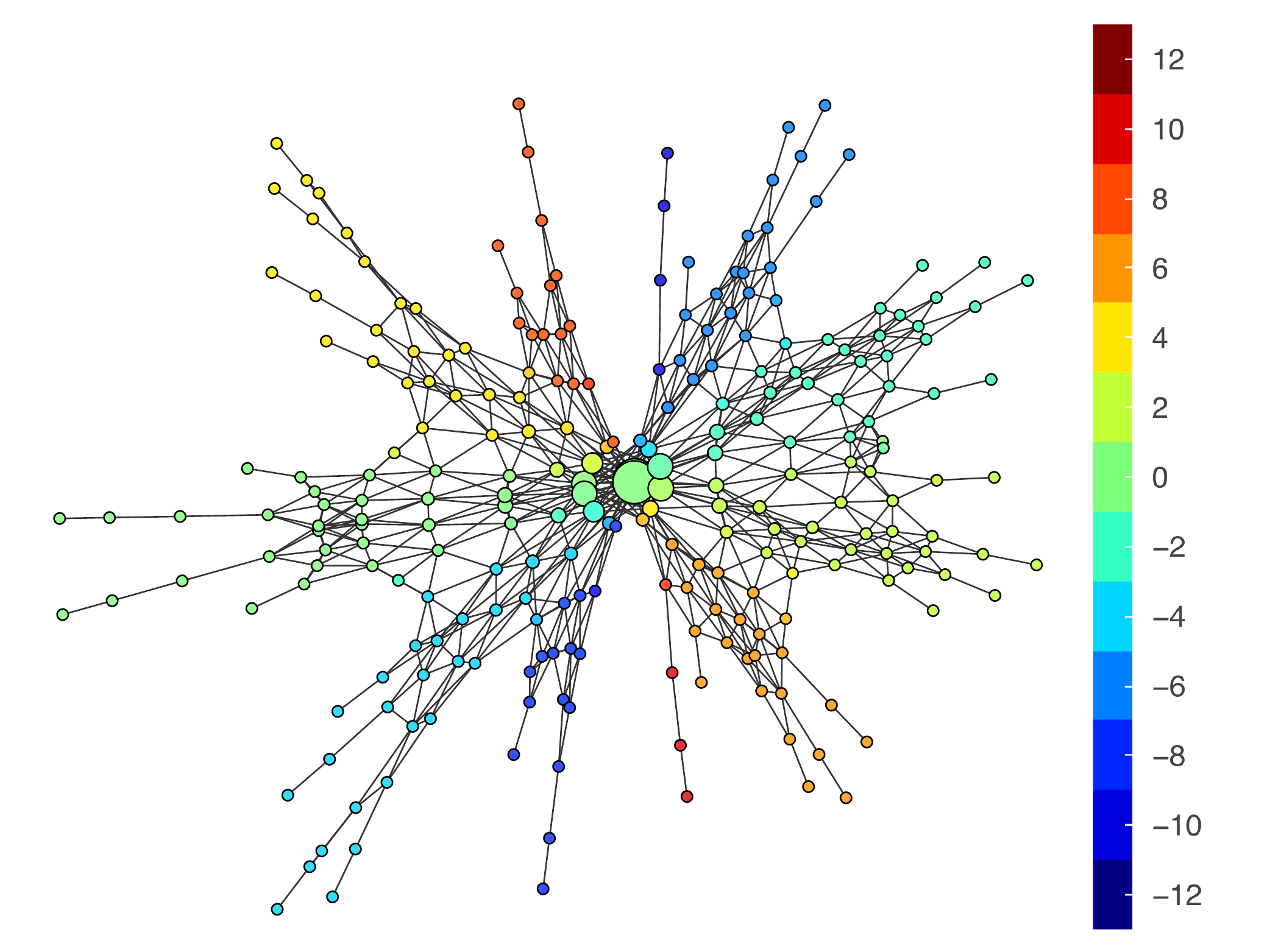}
          \caption{}
          %up to 17 $\epsilon=100$ with $890$ nodes.}
      \end{subfigure}
     \caption{Stability of \NEW{Equivariant} Ball Mapper with respect to the choice of parameter/radius $\epsilon$:  EqBM graphs of knots up to $17$ crossings $\epsilon=50$ with $3840$ nodes (A), $\epsilon=100$ with $890$ nodes (B), $\epsilon=200$ with $254$ nodes (C). Color is determined by the average signature of knots in that node/cluser.}
     \label{fig:stability_Jones_construction_epsilon}
\end{figure}
\fi

The Ball Mapper graph for the Alexander polynomial data has linear structure, see Figure \ref{fig:other_BM}(B). 
\NEW{The two flares contain clusters of knots whose  signature modulo $4$ is equal to zero or two, respectively.}
Similarly, the Ball Mapper graphs for HOMFLYPT data exhibit a star-like structure whose flares contain knots with the same signature, Figure \ref{fig:other_BM}(C). As in the case of the Jones polynomial data, all Ball Mapper graphs are stable with respect to the crossing number filtration and the choice of Ball Mapper parameter $\epsilon$.

\ifdefined\ifNotPictures
\begin{figure}[h!]
    \centering
       \begin{subfigure}[t]{0.3\textwidth}
         \centering
         \includegraphics[width=\textwidth]{jones_upto_17/jones_upto_17_100_signature-eps-converted-to.pdf}
         \caption{}%Jones up to 17 crossings.}
     \end{subfigure}
    \begin{subfigure}[t]{0.3\textwidth}
    \includegraphics[width=\textwidth]{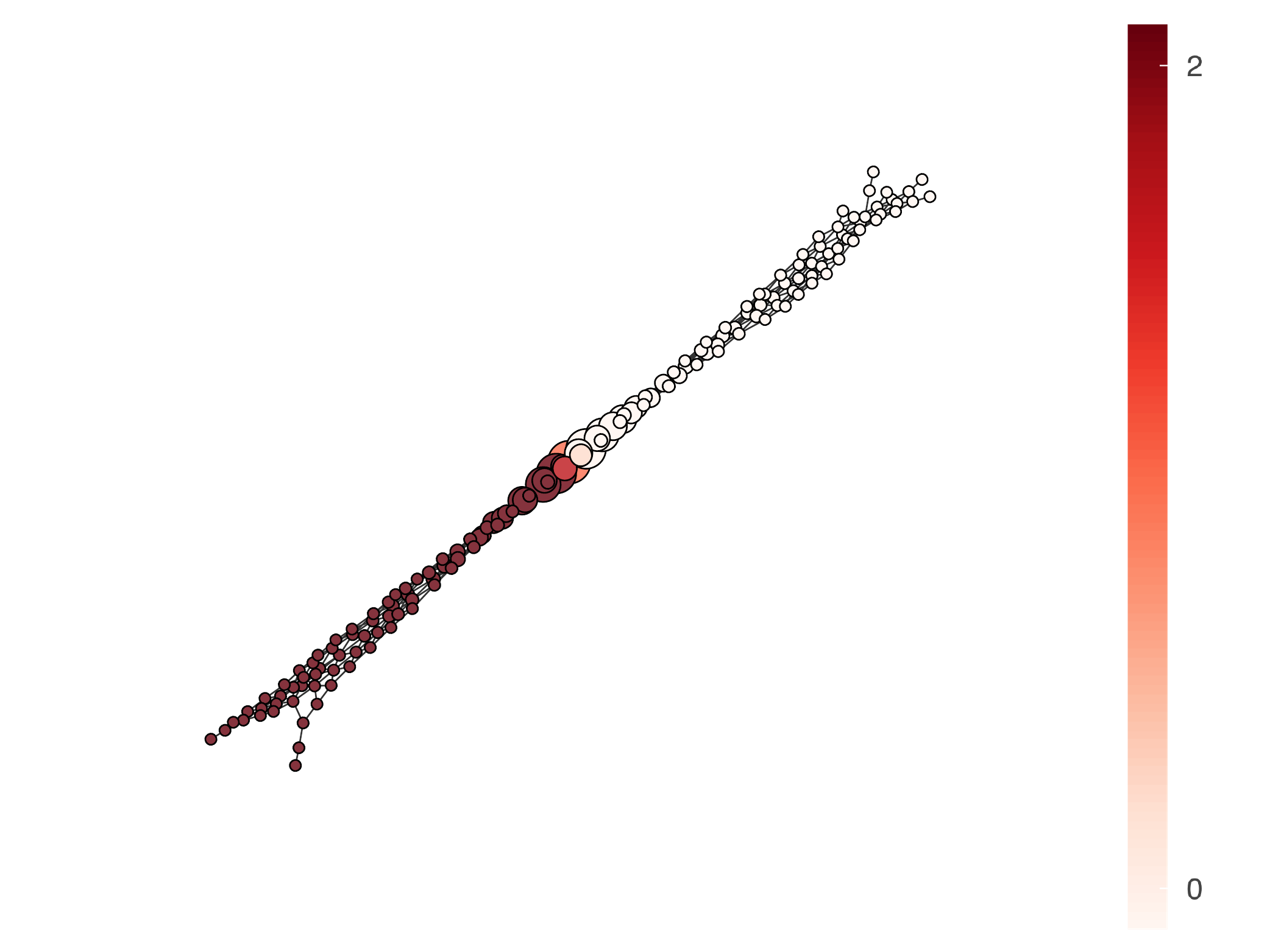}
    \caption{}%Alexander up to 17  crossings}
    \end{subfigure}
    \begin{subfigure}[t]{0.3\textwidth}
    \includegraphics[width=\textwidth]{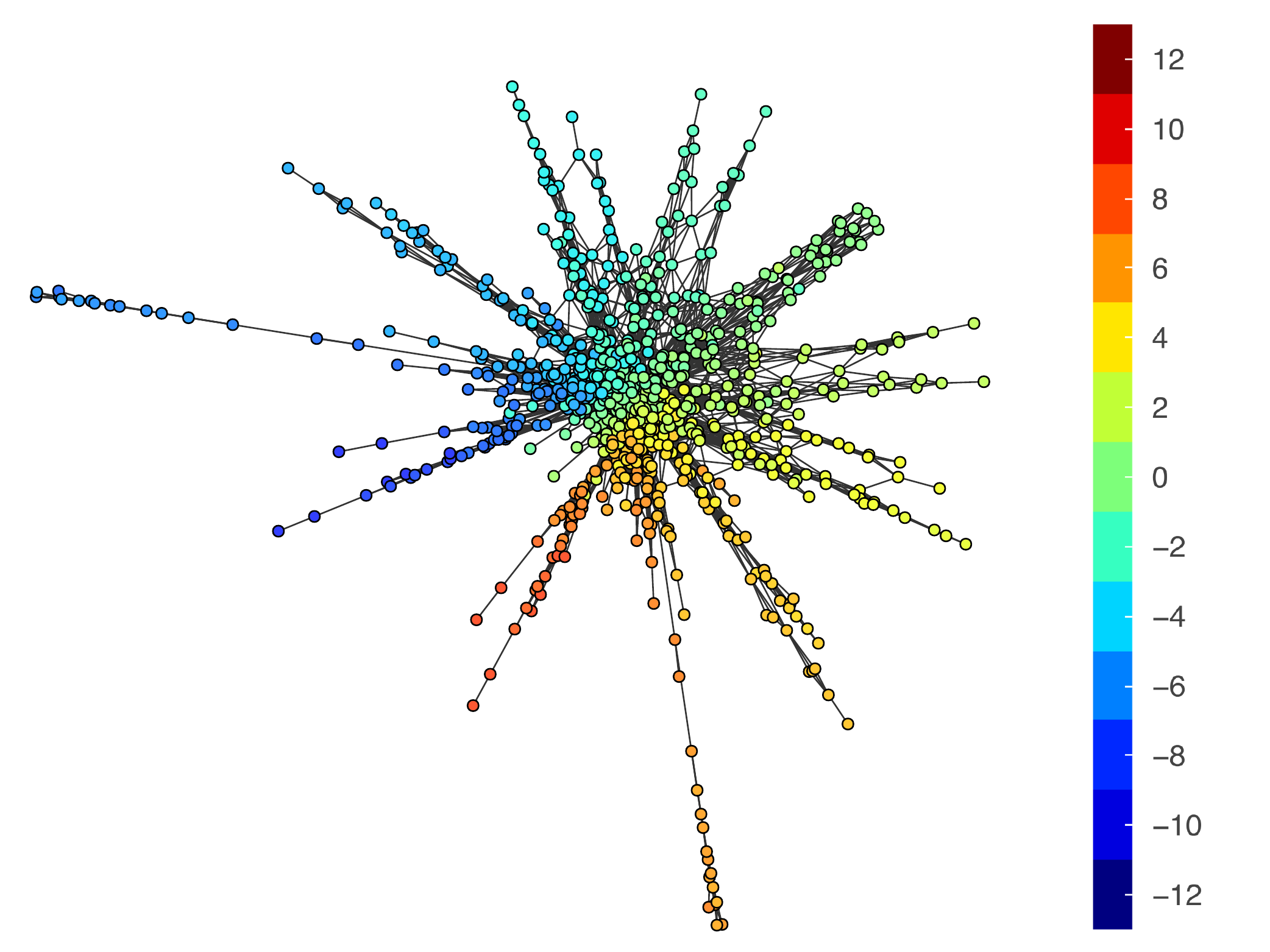}
    \caption{}%HOMFLYPT up to 15  crossings}
    \end{subfigure}
    %\begin{subfigure}[b]{0.3\textwidth}
    %\includegraphics[width=\textwidth]{khovanov_upto_15n/khov_upto_15n_20_svalue-eps-converted-to.pdf}
    %\caption{Khovanov non-alternating up to 15 crossings}
    %\end{subfigure}
    \caption{Equivariant Ball Mapper graphs for the 
    \ifdefined\showOldText 
    \OLD{Alexander (A), Jones(B)} 
    \fi
    \NEW{Jones (A), Alexander (B)}, and HOMFLYPT  (C) polynomial data of knots up to 17 crossings colored by the average signature \NEW{of knots}  in each cluster \NEW{(A), (C)} or  signature mod 4 in (B). }
     \label{fig:other_BM}
\end{figure}
\fi

%$$$$$$$$$$$$$$$$$$$$$$$$$$$$$$$$$$$$$$$$$$$$$$$$$

\subsection{MappingMappers and Mapper on Ball Mapper\NEW{: comparing knot invariants}}
\label{sec:comparisons}
In this section, we further investigate data corresponding to the Jones, Alexander, and HOMFLYPT polynomials whose Ball Mapper graphs are provided in Section \ref{sec:KnotDataBM}.  The main goal is to compare these spaces using Ball Mapper based tools developed in Sections \ref{sec:relation_between_mapper} and \ref{sec:BonBM}, in order to, implicitly, compare the invariants. First, spaces of two invariants can be compared by constructing their Ball Mapper graphs and visualizing the maps between them as described in Section \ref{sec:relation_between_mapper}. 
Next,  Mapper on Ball Mapper construction from Section \ref{sec:BonBM}  is used to emphasize relative strengths of two invariants with respect to distinguishing knots. To be more specific, given two invariants $A$ and $B$, in general data descriptors, of a dataset $\mathcal{K}$ we can think of them as maps $A : \mathcal{K} \rightarrow M_A$ and $B : \mathcal{K} \rightarrow M_B$, where $M_A$ and $M_B$ are metric spaces. 
Most commonly, $M_A$ and $M_B$ are finite point clouds in Euclidean spaces of different dimensions.  Inspired by comparison of knot invariants,  \NEW{ that invariant $A$ is considered to be \emph{stronger} than invariant $B$ if the elements covered by a single vertex or several closely-connected vertices in the Mapper on Ball Mapper graph of $M_B$ are spread across different regions of the Mapper on Ball Mapper graph of $M_A$.} MappingMappers and Mapper on Ball Mapper can be used for this type of analysis. 

\ifdefined\showOldText 
\OLD{Visualizing functions}
\fi
\NEW{MappingMappers, defined}  in  Section \ref{sec:relation_between_mapper}, uses two point clouds $X  \subset \mathbb{R}^n$ and $Y \subset \mathbb{R}^m$ and a relation $f \subset X\times Y$ as inputs. To illustrate this technique we use the 
 collection $\mathcal{K}$ of knots up to $17$ crossings along with the Jones  $J(\mathcal{K}) \subset \mathbb{R}^{51}$ 
 and the Alexander $A(\mathcal{K})\subset \mathbb{R}^{17}$ point clouds 
 obtained in Section \ref{sec:Data}. 

\begin{figure}[h!]
     \centering
     \begin{subfigure}[b]{0.4\textwidth}
         \centering
         \includegraphics[width=\textwidth]{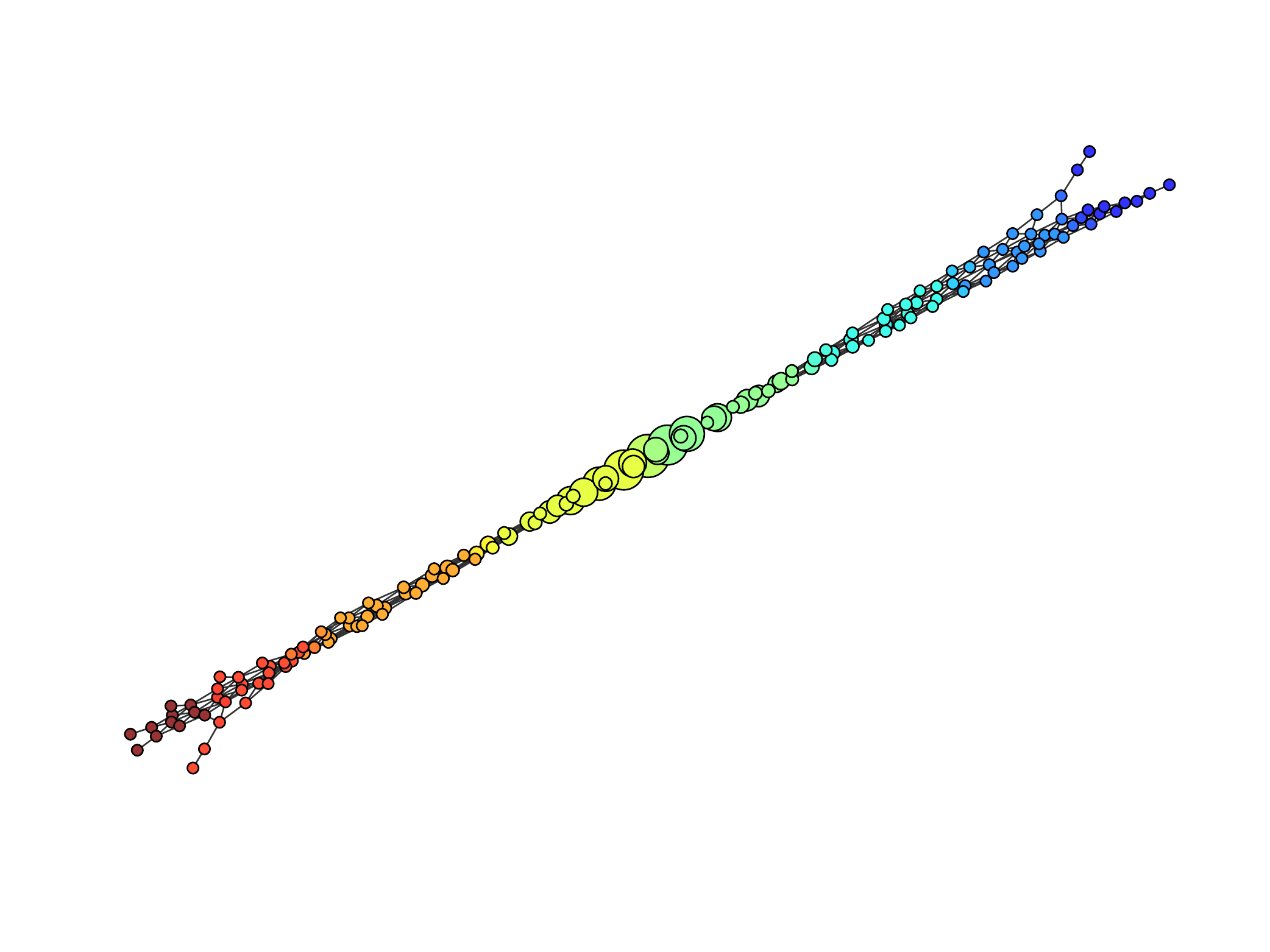}
         \caption{}
     \end{subfigure}
     \hspace{1cm}
     \begin{subfigure}[b]{0.4\textwidth}
     \includegraphics[width=0.9\textwidth]{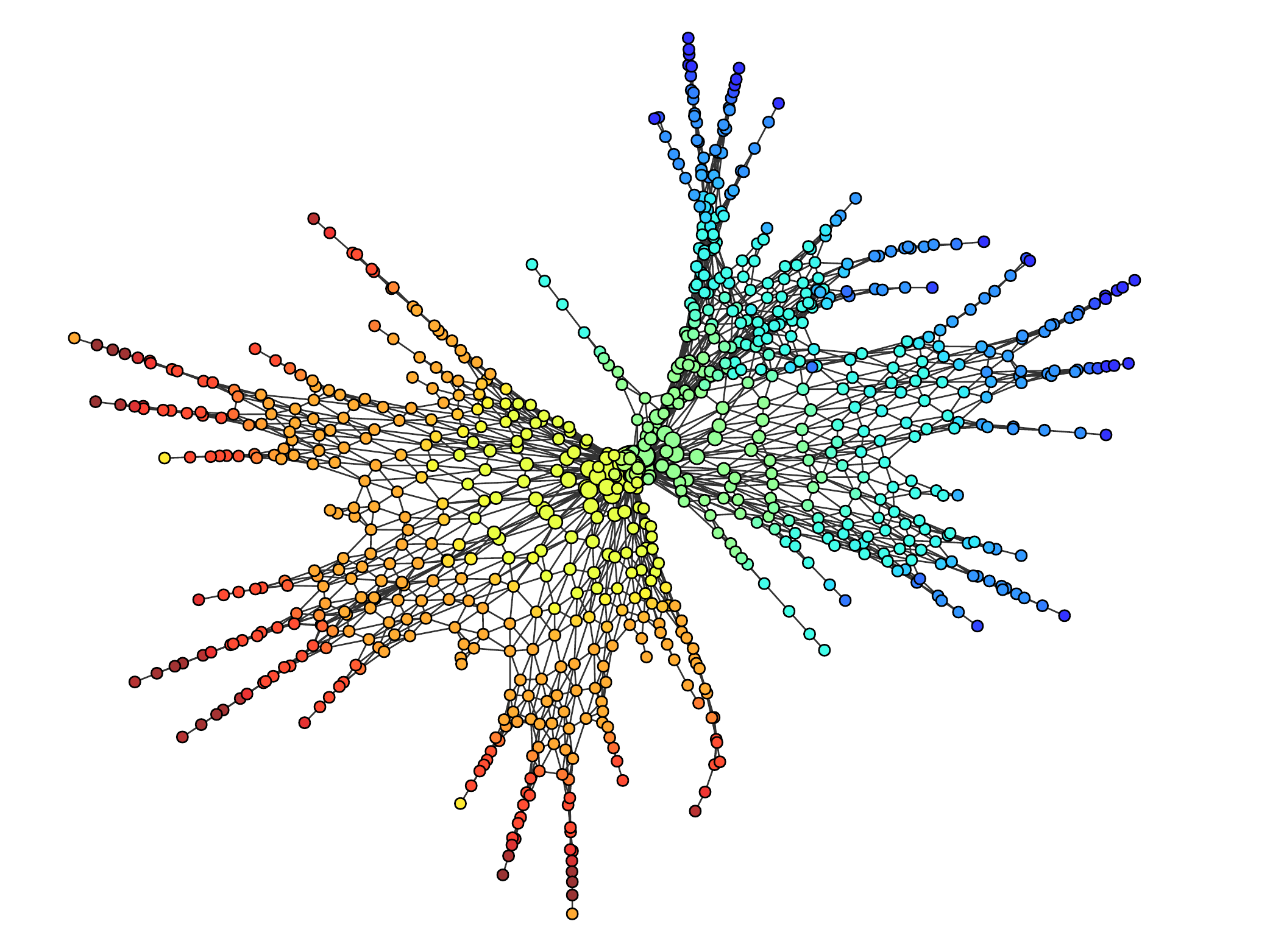}
      \caption{}
     \end{subfigure}
    % \hfill
        \caption{MappingMappers: Representation of  a map from the space of the  Alexander (A) and the space of the Jones polynomials of  knots up to 17 crossings (B). \NEW{Rainbow coloring of consecutive clusters of the linear embedding of the Equivariant Ball Mapper graph for the Alexander data and used to color the corresponding regions of the Equivariant Ball Mapper graph for the Jones data.} 
        \ifdefined\showOldText 
        \OLD{ This Ball Mappers are custom colored to show region correspondence, i.e. regions of the same color are mapped to each other.}
        \fi
        }
        \label{fig:mapping_alex_to_jones}
\end{figure}

\begin{figure}[h!]
     \centering
     \begin{subfigure}[b]{0.4\textwidth}
         \centering
         \includegraphics[width=\textwidth]{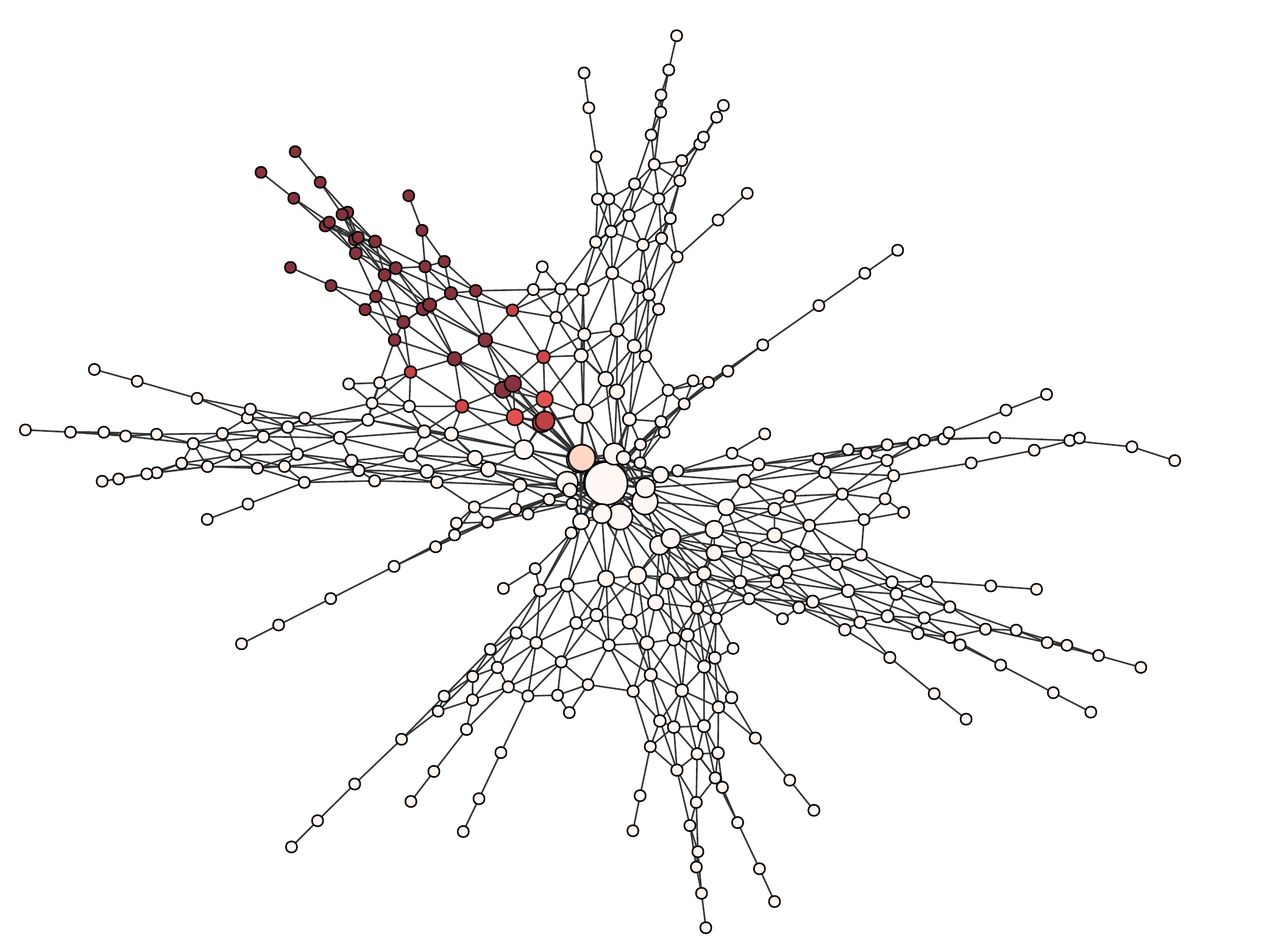}
         \caption{}%{Ball Mapper graph (containing $326$ nodes) of Jones polynomial of all knots up to 15 crossings with $\epsilon=50$.}
     \end{subfigure}
    % \hfill
     \begin{subfigure}[b]{0.4\textwidth}
         \centering
         \includegraphics[width=\textwidth]{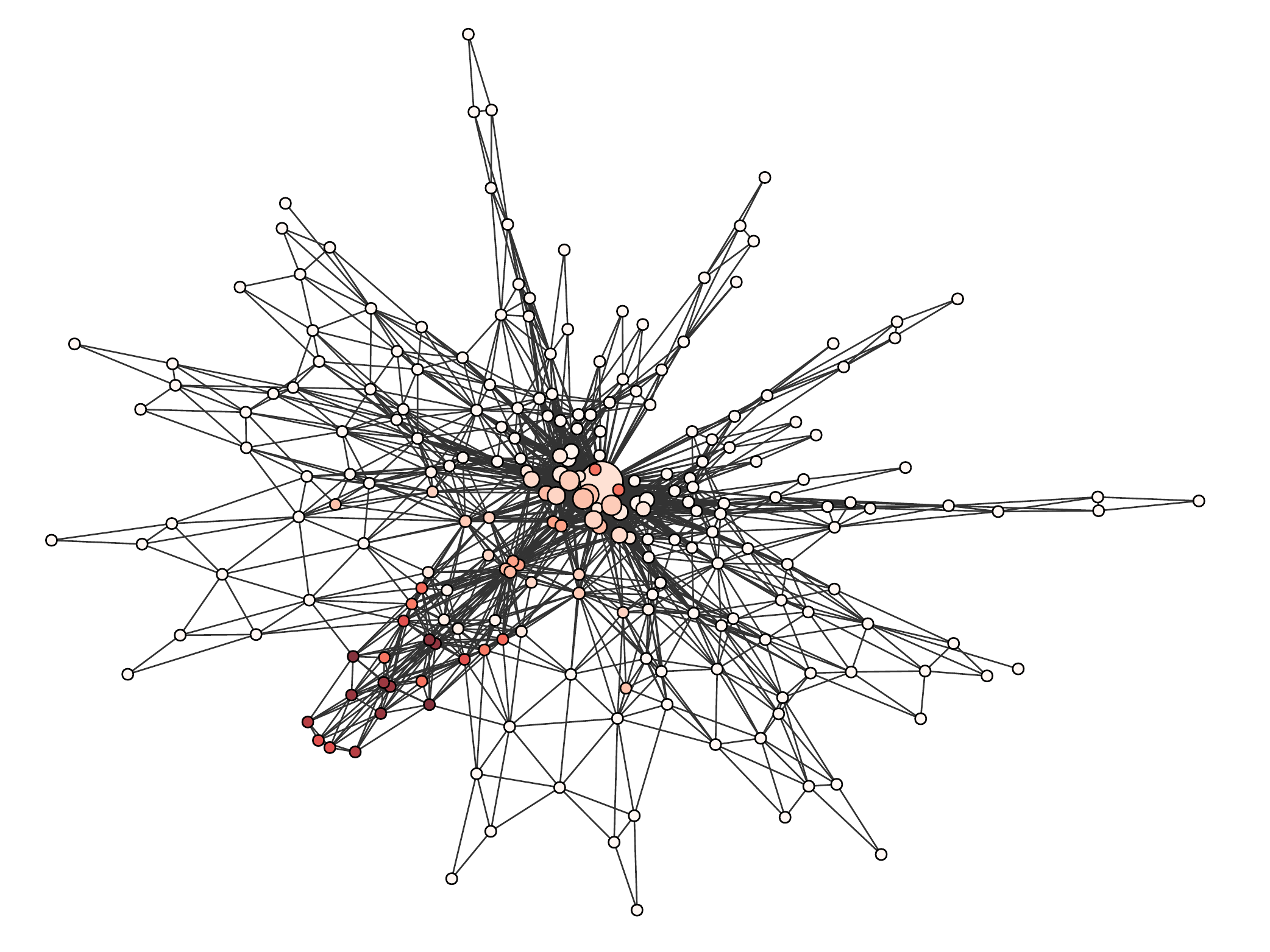}
         \caption{}
     \end{subfigure}
        \caption{MappingMappers: Visualizing the map from the space of Jones polynomials to the space of HOMFLYPT polynomials of knots up to 15 crossings using the Ball Mapper graph (containing $326$ nodes) of Jones polynomial with $\epsilon=50$ (A) and the Ball Mapper graph (containing $258$ nodes) of HOMFLYPT polynomial of 
 with $\epsilon=50$ shown in (B). All clusters containing knots with signature equal to zero in both Ball Mapper graphs are shown as a shade of red.}
        \label{fig:mapping_jones_to_homfly}
\end{figure}

Using the set of knots $\mathcal{K}$ as the common indexing set, we define the relation $f \subset A(\mathcal{K}) \times J(\mathcal{K})$ in the following way: for a given knot $K \in \mathcal{K}$ its Alexander polynomial $A(K)$ in $A(\mathcal{K})$ is related to its Jones polynomial $J(K)$ in $J(\mathcal{K})$. Note that this is a relation rather than a function, since some knots from $\mathcal{K}$ may have the same Alexander polynomial, but different Jones polynomials, and vice versa. In this case, a single point from $A(\mathcal{K})$ can be related to multiple points of $J(\mathcal{K})$ or the other way around.
Figure~\ref{fig:mapping_alex_to_jones} illustrates this relation: colors indicate matching regions in the Alexander and Jones Ball Mapper graphs. Roughly speaking, the linear structure of the Alexander Ball Mapper induces the linear structure among the flares in the star-like Ball Mapper graph of the Jones data. In the opposite direction, flares of the Jones Ball Mapper merge according to the signature module $4$ which is consistent with the fact that one end of the Alexander Ball Mapper contains knots whose signature mod 4 is zero and the other one is two.  Figure~\ref{fig:mapping_jones_to_homfly} illustrates the  analogous relation between the Jones and HOMFLYPT Ball Mapper graphs, respectively. \NEW{The non-linear nature of the Jones Ball Mapper graph prohibits using a gradient-like coloring; instead, the cluster color in this graph reflects the percentage of knots with signature equal to zero and the corresponding clusters in the HOMFLYPT Ball Mapper graph.}

\NEW{To illustrate how methods introduced in this paper can be used to discover potentially new results in knot theory, let us reconsider the Ball Mapper graph in Figure \ref{fig:other_BM}(B). The distribution of signature values suggests that signature mod 4 can be determined as a function of the coefficients of the Alexander polynomial. This hypothesis was tested by training a Support Vector Machine classifier~\cite{bishop2006pattern}.  
According to it, the perfect separation between two classes of knots, those whose signature mod 4 is zero and those with it equal to 2, is achieved
with an "anti-diagonal" hyperplane with normal vector $[1, -1, 1, \ldots]$. This observation indicates that the sign of the alternating sum of the coefficients of the Alexander polynomial determines the signature  mod 4. 
The correspondence between clusters of Ball Mapper graphs of the Alexander and Jones polynomials shown in Figure~\ref{fig:mapping_alex_to_jones} suggests that analogous property should hold for the Jones polynomial. Turns out that both statements are true according to the well-known theorem \cite{conway1970enumeration}. Our approach provides a new  way of obtaining 
such theorem  and paves the way for using mapper-type algorithms to aid  discovery in knot theory in particular and  theoretical mathematics and sciences in general.
}
% 
%
%
%
%$$$$$$$$$$$$$$$$$$$$$$$$$$$$$$$$$$$$$$$$$$$$$$$$$

\NEW{Next, we employ MoBM on the} algebraic relation between the HOMFLYPT polynomial and both the Alexander and Jones polynomials, see Section \ref{sec:KnotIntro}. 
These specializations tie in with the framework introduced in  Section \ref{sec:BonBM} as they can be used as lenses for  Mapper algorithm when using the Mapper on Ball Mapper, or MoBM, construction. 

As a clustering method in the MoBM construction, we use the DBSCAN algorithm~\cite{dbscan}. DBSCAN requires a new parameter $\epsilon_{DB}$  in addition to the $\epsilon$ denoting the radius of balls used in the Ball Mapper construction.

 The MoBM construction using knots data is illustrated in  Figure \ref{fig:pullback_HOMFLYPT_JONES_small}:  Ball Mapper on Jones data Figure \ref{fig:pullback_HOMFLYPT_JONES_small} (A) is used as the input covering to the MoBM graph \ref{fig:pullback_HOMFLYPT_JONES_small} (B) the HOMFLYPT data  whose Ball Mapper is in Figure \ref{fig:mapping_jones_to_homfly}(B). The coloring in Figure  \ref{fig:pullback_HOMFLYPT_JONES_small} (A) represents the number of clusters into which the points in each node split when they are pulled back from the space of Jones to HOMFLYPT. \NEW{This pullback is not trivial (clusters split into more than one cluster in the pre-image) in the center region and in between the flares.} The obtained MoBM graph thus achieves better separation in those regions.

\ifdefined\ifNotPictures
\begin{figure}[h!]
     \centering
     \begin{subfigure}[b]{0.3\textwidth}
         \centering
         \includegraphics[width=\textwidth]{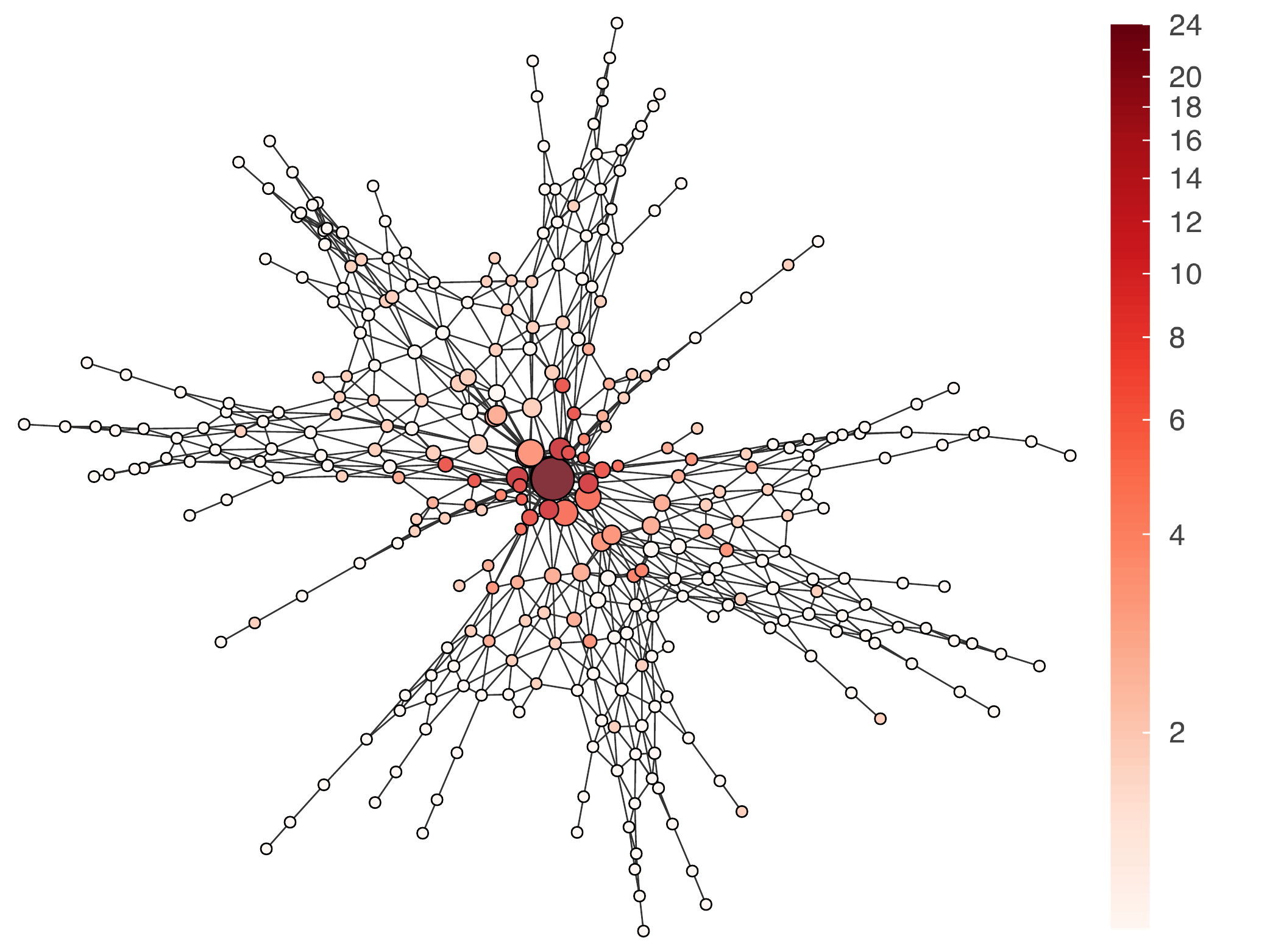}
         \caption{}
        %  Jones BM $\epsilon=50$ with $326$ nodes colored by the number of DBSCAN clusters .}
     \end{subfigure}
     \hfill
     \begin{subfigure}[b]{0.3\textwidth}
         \centering
         \includegraphics[width=\textwidth]{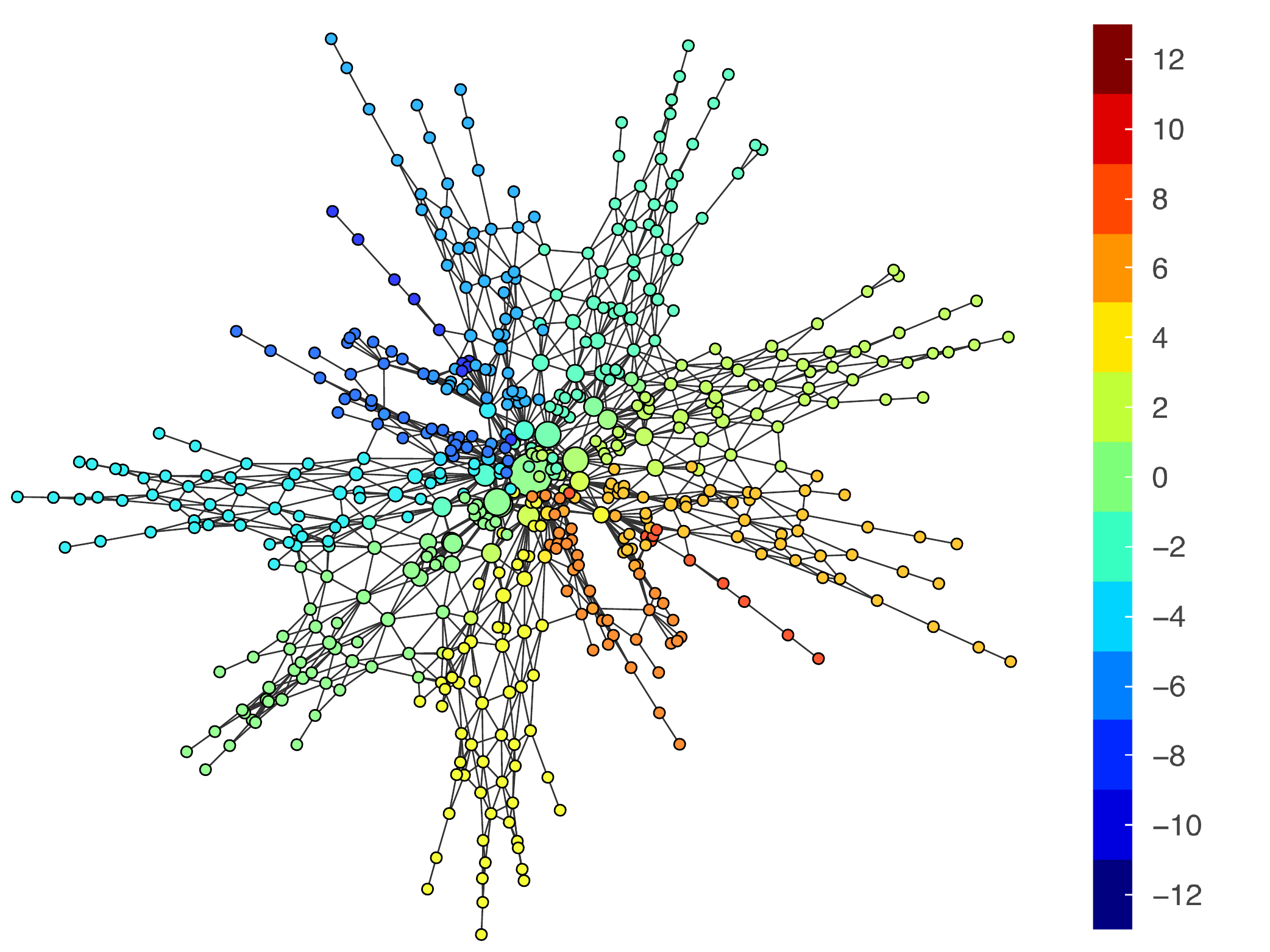}
         \caption{}
        %  Pullback to HOMFLY $\epsilon_{DB}=40$ with $644$ nodes.}
     \end{subfigure}
    \hfill
     \begin{subfigure}[b]{0.3\textwidth}
         \centering
         \includegraphics[width=\textwidth]{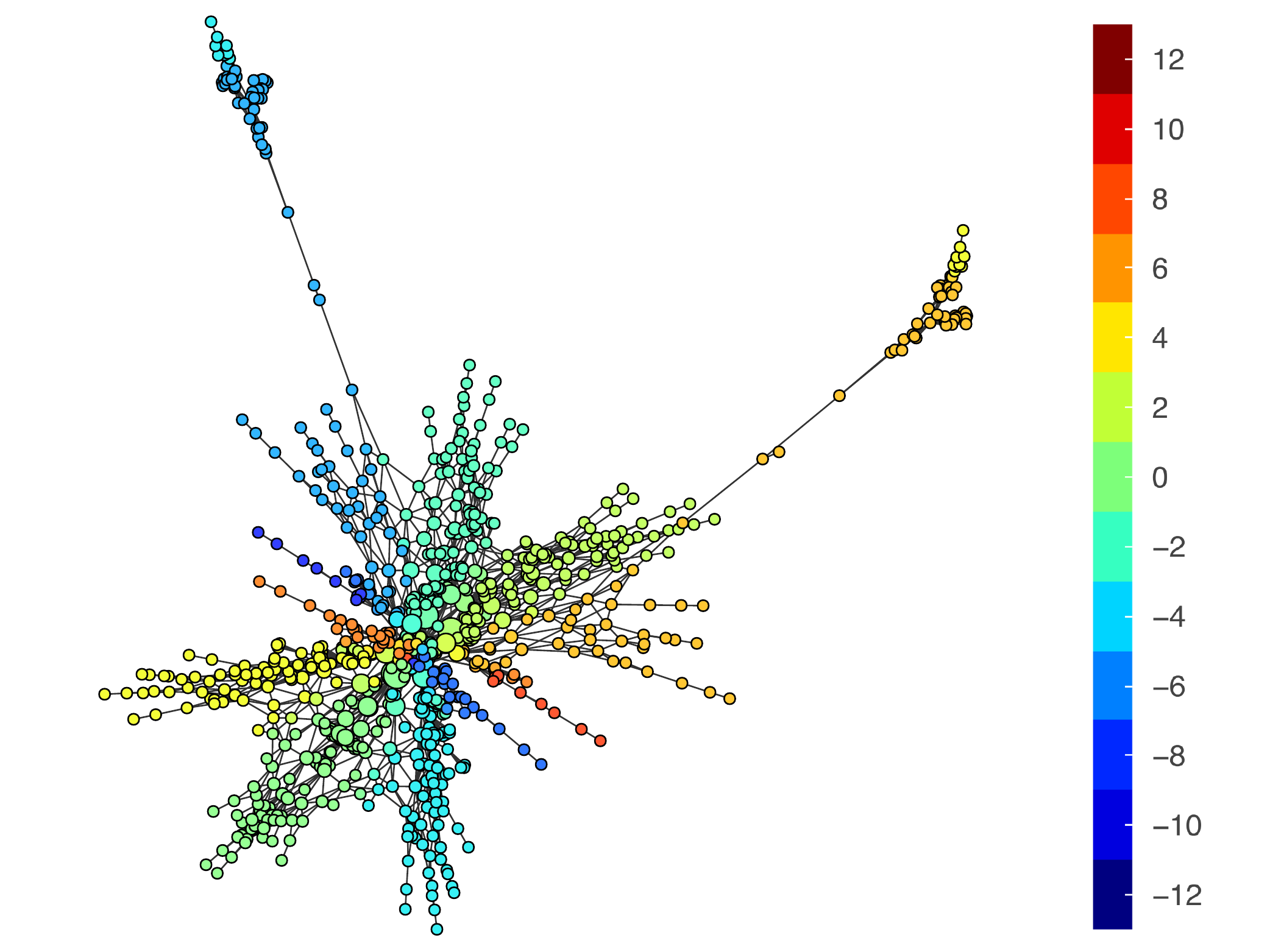}
         \caption{}
        %  Pullback to HOMFLY $\epsilon_{DB}=40$ with $644$ nodes.}
     \end{subfigure}
        \caption{The Mapper on Ball Mapper graph on Jones-HOMFLYPT data pair. Image (A) shows the Ball Mapper graph for Jones data for knots up to 15 crossings at  $\epsilon_{BM}=50$ with a total of $326$ nodes colored by the number of clusters found in each node. The Mapper on Ball Mapper construction in (B)  is obtained from the Jones Ball Mapper (A) and the HOMFLYPT data for the DBSCAN clustering algorithm parameter equal $\epsilon_{DB}=40$ with $644$ nodes, and colored by signature. \NEW{  Analogous MoBM construction for $\epsilon_{DB}=30$ reveals a different structure on HOMFLYPT data (C): the two long flares emerging in (C) consist knots with the same s-invariant \cite{rasmussen2010khovanov}, but different values of signature.} }
        % \RS{great}\PD{Deeper analysis indicates that the two flares emerging in (C) are build up with knots of the same s-invariant, but different values of signature.}        
        % }
        \label{fig:pullback_HOMFLYPT_JONES_small}
\end{figure}
\fi

\section{Applications: Game theory, materials science and cancer research}
\label{sec:further_examples}
The techniques discussed in this paper apply to a wide range of datasets. In this section we present additional sample applications of the proposed techniques to artificial datasets, game theory, materials science and cancer research.

\subsection{Equivariant Ball Mapper: Tic-Tac-Toe data}
\label{sec:additional_equivariant_bm}
Any dataset with a nontrivial isometric group action can be visualized with preservation of the action using Equivariant Ball Mapper. As an example we analyze the Tic-Tac-Toe endgame dataset~\cite{Dua:2019} that consists of all possible board configurations at the end of Tic-Tac-Toe games. \NEW{These configurations, numbering 958 in total, are represented by 3 by 3 matrices (interpreted as vectors in $\mathbb{R}^9$) since the game is played between two players on a 3 by 3 grid where the first player places noughts and the  second places crosses.} The winner is the first player who places three noughts or crosses in a vertical, horizontal or a diagonal line.

The input for our analysis consists of $3$ by $3$ grids, interpreted as vectors with $9$ entries, with values $-1$ (corresponding to a nought), $0$ (corresponding to an empty slot) or $1$ (corresponding to a cross). 
The symmetries of the $3 \times 3$ configurations are given by a dihedral group consisting of four rotations and four reflections. It is straightforward to see that all configurations in one orbit are all wins, loses or ties, since rotation and reflection of the board does not change the outcome of the game. These symmetries induce relations between vectors in $\mathbb{R}^9$ and the Euclidean distance between any two configurations and their images via one of the actions will be the same.
Hence, Equivariant Ball Mapper is the natural choice for this data; see the resulting Ball Mapper graph in the Figure~\ref{fig:tictactoe}(A). 

\ifdefined\ifNotPictures
\begin{figure}[h!]
     \centering
     \begin{subfigure}[t]{0.3\textwidth}
         \centering
         \includegraphics[width=\textwidth]{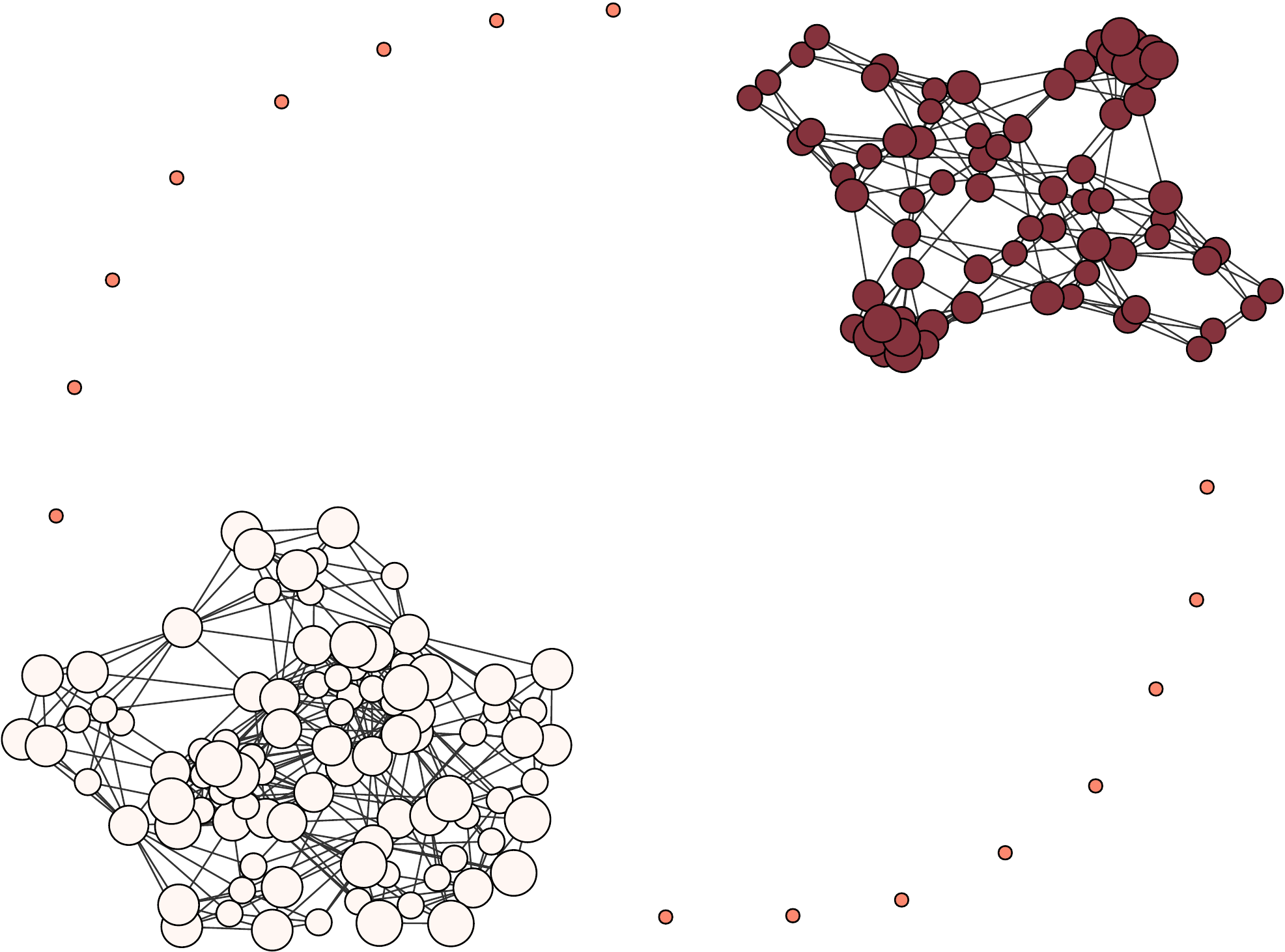}
         \caption{}%Equivariant Ball Mapper colored by outcome.}
         \label{subfig:tic_A}
     \end{subfigure}
     \hfill
     \begin{subfigure}[t]{0.3\textwidth}
         \centering
         \includegraphics[width=\textwidth]{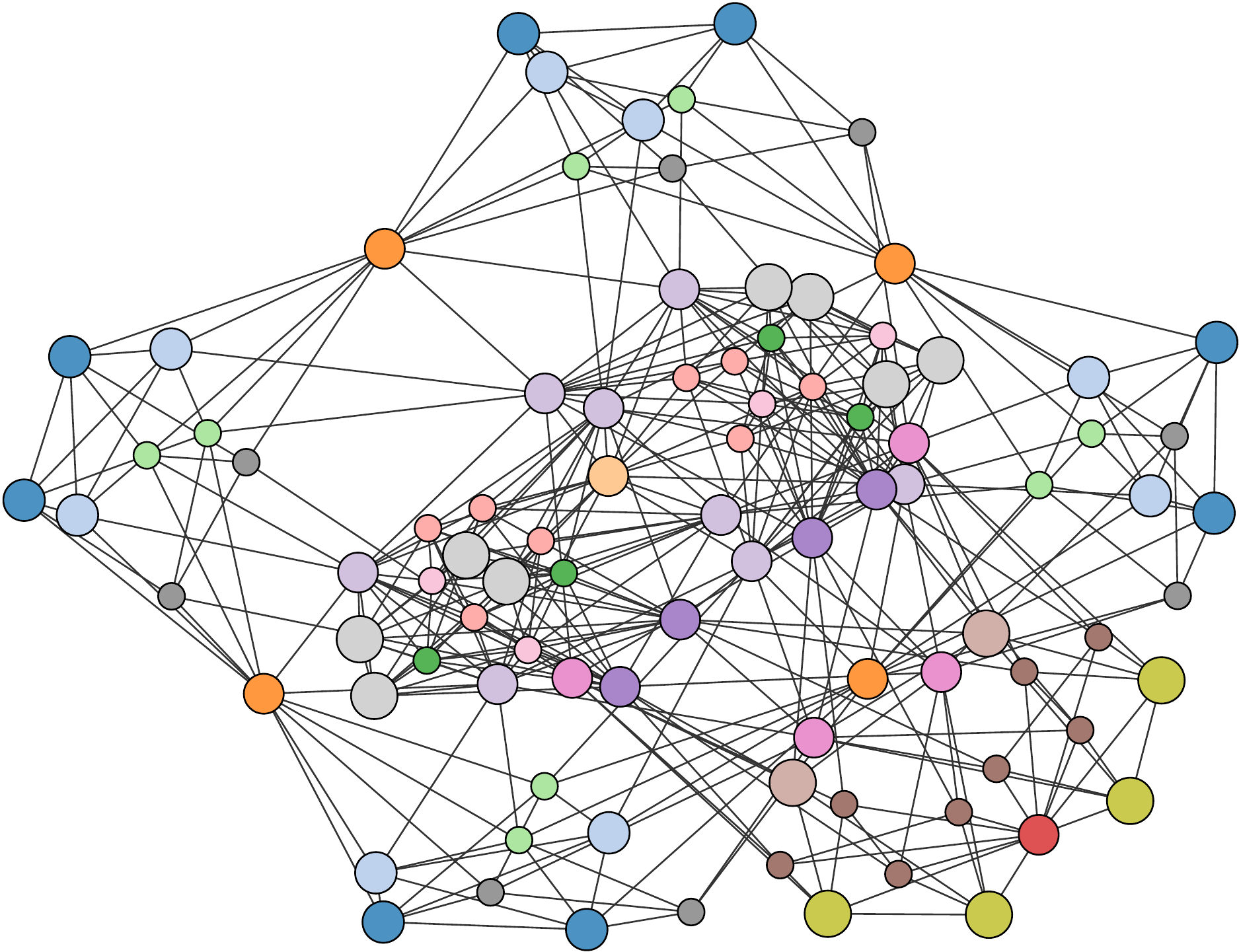}
         \caption{}%Orbits in the wins cluster.}
         \label{subfig:tic_B}
     \end{subfigure}
     \hfill
     \begin{subfigure}[t]{0.3\textwidth}
         \centering
         \includegraphics[width=\textwidth]{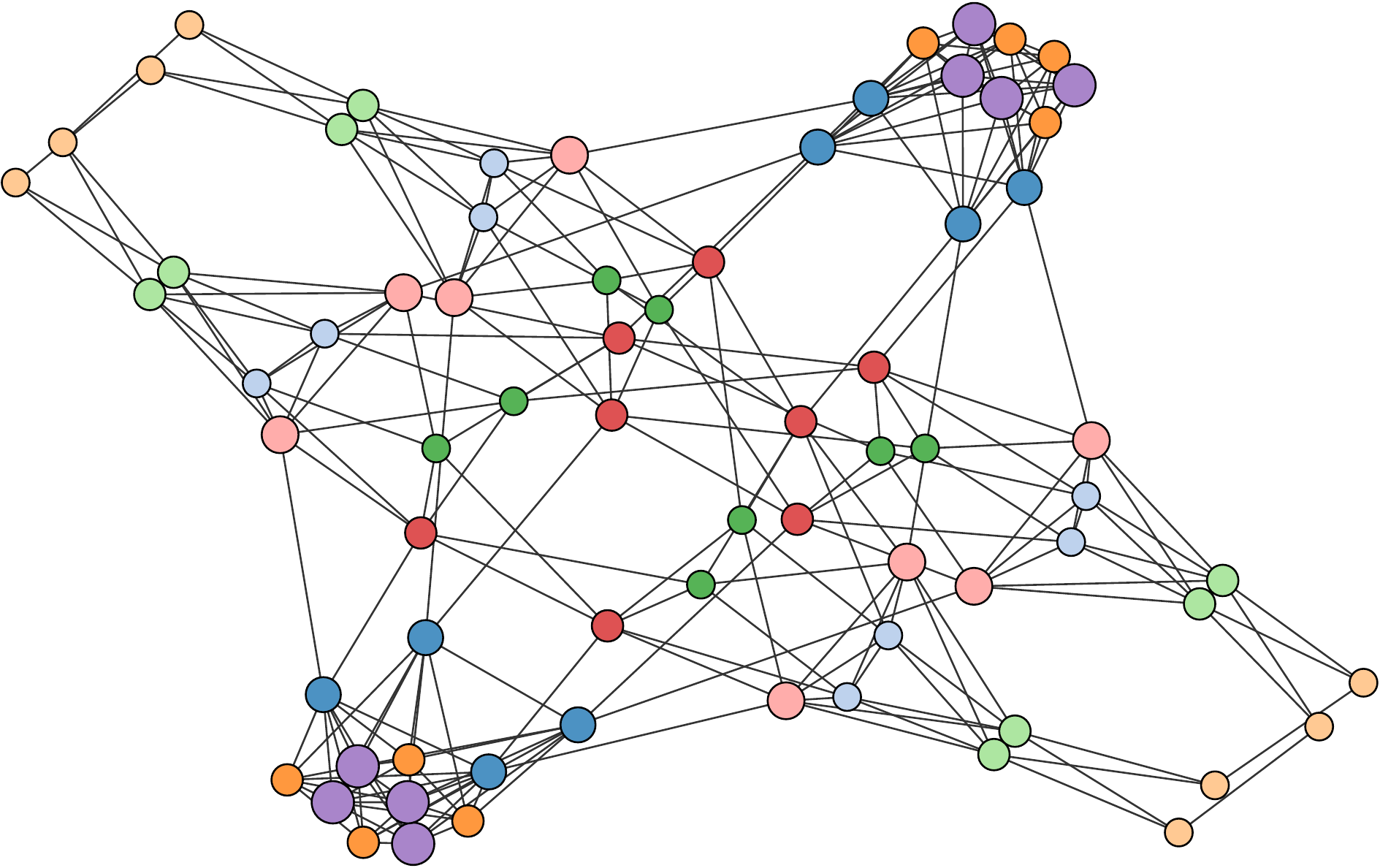}
         \caption{}%Orbits in the losses cluster.}
         \label{subfig:tic_C}
     \end{subfigure}
        \caption{The Equivariant Ball Mapper graph for the Tic-Tac-Toe dataset. \NEW{Figure \ref{subfig:tic_C} shows the game outcome: the wins clusters are colored white, the losses clusters red.} The sixteen isolated orange clusters correspond to all possible ties. The wins and losses clusters are shown in \ref{subfig:tic_B} and \ref{subfig:tic_C} respectively, with color denoting orbits. 
        Note that even if the same color palette is used, there is no relation between the nodes in \ref{subfig:tic_B} and \ref{subfig:tic_C} as the dihedral symmetry does not change the outcome of a game.}
        \label{fig:tictactoe}
\end{figure}
\fi

While our purpose is to showcase the technique rather than to draw a conclusion about the output, we include several observations. 
\NEW{The Equivariant Ball Mapper graph for radius $\epsilon = 2.5$ (using $l_1$ distance) provides a perfect separation between the configurations} in which the first player wins the game (white), loses the game (red) or when there is a tie (orange), see Figure~\ref{fig:tictactoe}(A). We attribute the clear separation of win-loss-tie clusters, as well as the symmetries of the win clusters (Figure~\ref{fig:tictactoe}(B)) and loss clusters (Figure~\ref{fig:tictactoe}(C)), to the combinatorial properties of the game. In panels (B) and (C) nodes belonging to the same orbit are colored with the same color. Note that different orbits might have different lengths. A configuration with no rotation or reflection symmetry will lead to a length 8 orbit. On the other hand, a configuration that already has some symmetries will have a shorter orbit. The maximally symmetric configuration has an orbit of length 1. This corresponds to the only red node in the bottom right of Figure~\ref{fig:tictactoe}(B). 
\ifdefined\showOldText 
\OLD{Further analysis of this data is a part of the  forthcoming knot theory paper.} 
\fi
\NEW{Intuition on the perfect win-loss-tie separation can be built by considering the smallest $l_1$ distance between configurations. It is not difficult to show that the distance between winning and losing configurations is at least $3$. Similarly, the minimum distance between two ties is $4$ and the minimum distance between winning (resp. loosing) configuration and a tie is $4$ (resp. $3$). On the other hand, any pair of winning (resp. losing) configurations can be connected by a sequence of winning (resp. loosing) configurations spaced by at most $2$. 
In lieu of the proof, consider a Ball Mapper graph with radius $2 < \epsilon < 3$, like the one depicted in Figure~\ref{fig:tictactoe}(A), where 
all the winning (resp. losing) configurations are in the same connected component evidencing the existence of such a path. 
}

%
%
%
%
%
%
%$$$$$$$$$$$$$$$$$$$$$$$$$$$$$$$$$$$$$$$$$$$$$$$$$
\subsection{MappingMappers: superconductors}
\label{sec:additional_relational_mapper}
As an example of a  map between two set of
descriptors of the same data, we  consider the  superconductor dataset available from the UC Irvine Machine Learning Repository~\cite{superconductor_dataset}. 

The dataset contains two descriptors of a family of $21263$ superconductors. The first one, the \emph{characteristics dataset}, is composed of $81$ features extracted from considered superconductors. The other, the \emph{composition dataset}, contains sparse vectors in $\mathbb{R}^{86}$  that describe  the chemical composition of the superconductor.

\NEW{A Ball Mapper graph for the characteristic dataset, with all features scaled to the interval $[0, 1]$ and with radius $\epsilon=2$, is provided in  Figure~\ref{fig:superconducting_dataset_BM}(A).
For the composition dataset, the Ball Mapper for any fixed radii  is either composed of a single vertex, or a large group of small, disconnected clusters. }
We attribute this behavior to a concentration of measure--type phenomena. In order to construct a meaningful Ball Mapper graph, we use a cosine similarity measure~\cite{cosine_similarity}  instead of the standard Euclidean metric. The corresponding Ball Mapper graph obtained for $\epsilon=0.25$ is shown in Figure~\ref{fig:superconducting_dataset_BM}(B). 
The red color in both Ball Mapper graphs on Figure~\ref{fig:superconducting_dataset_BM} denotes superconductors that work in the highest temperature and should be easiest to use in practice.

\ifdefined\ifNotPictures
\begin{figure}[h!]
     \centering
     \begin{subfigure}[t]{0.4\textwidth}
         \centering
         \includegraphics[width=\textwidth]{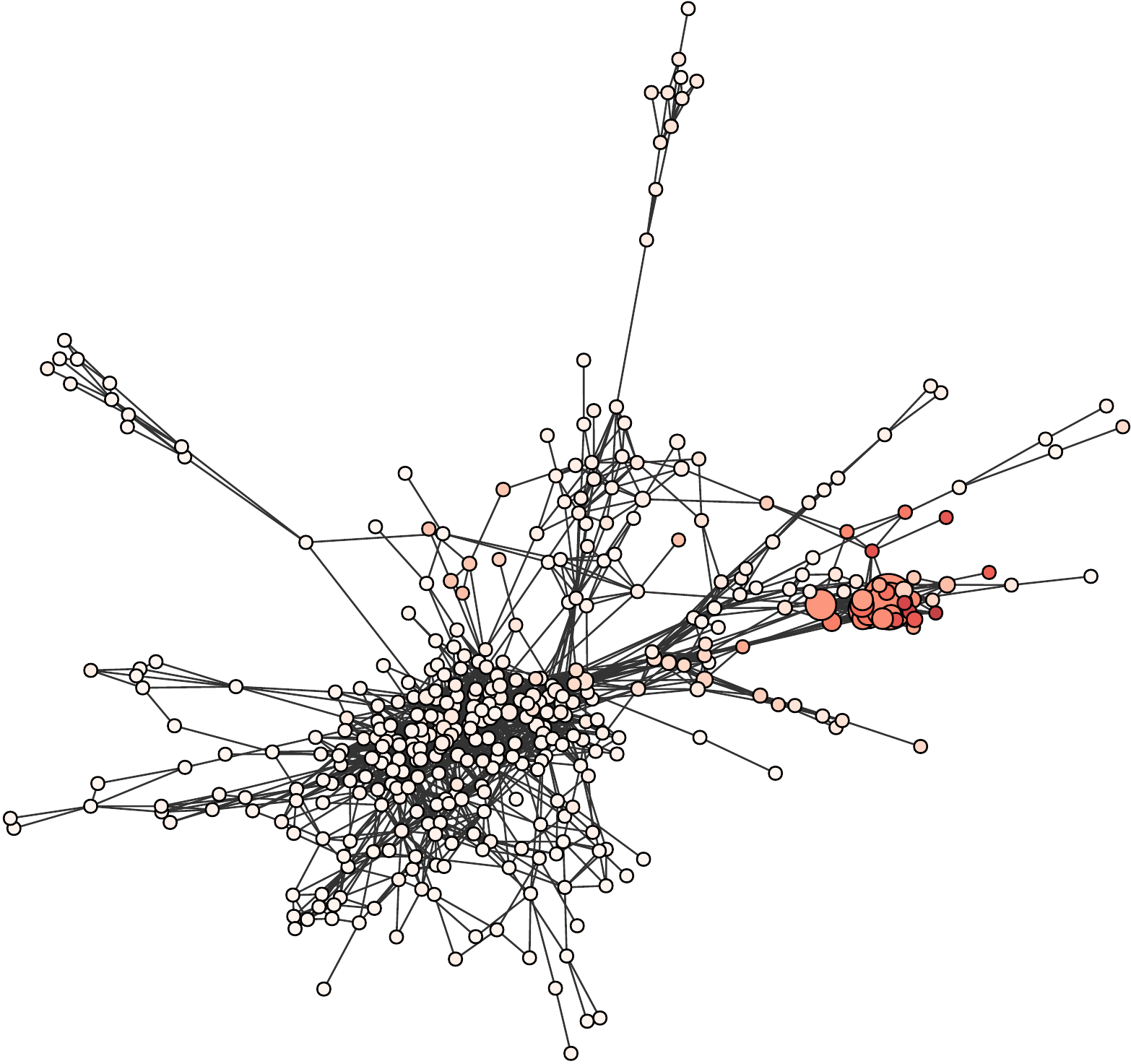}
         \caption{}%BM graph for the characteristic dataset.}
     \end{subfigure}
     \hspace{1cm}
     \begin{subfigure}[t]{0.4\textwidth}
         \centering
         \includegraphics[width=\textwidth]{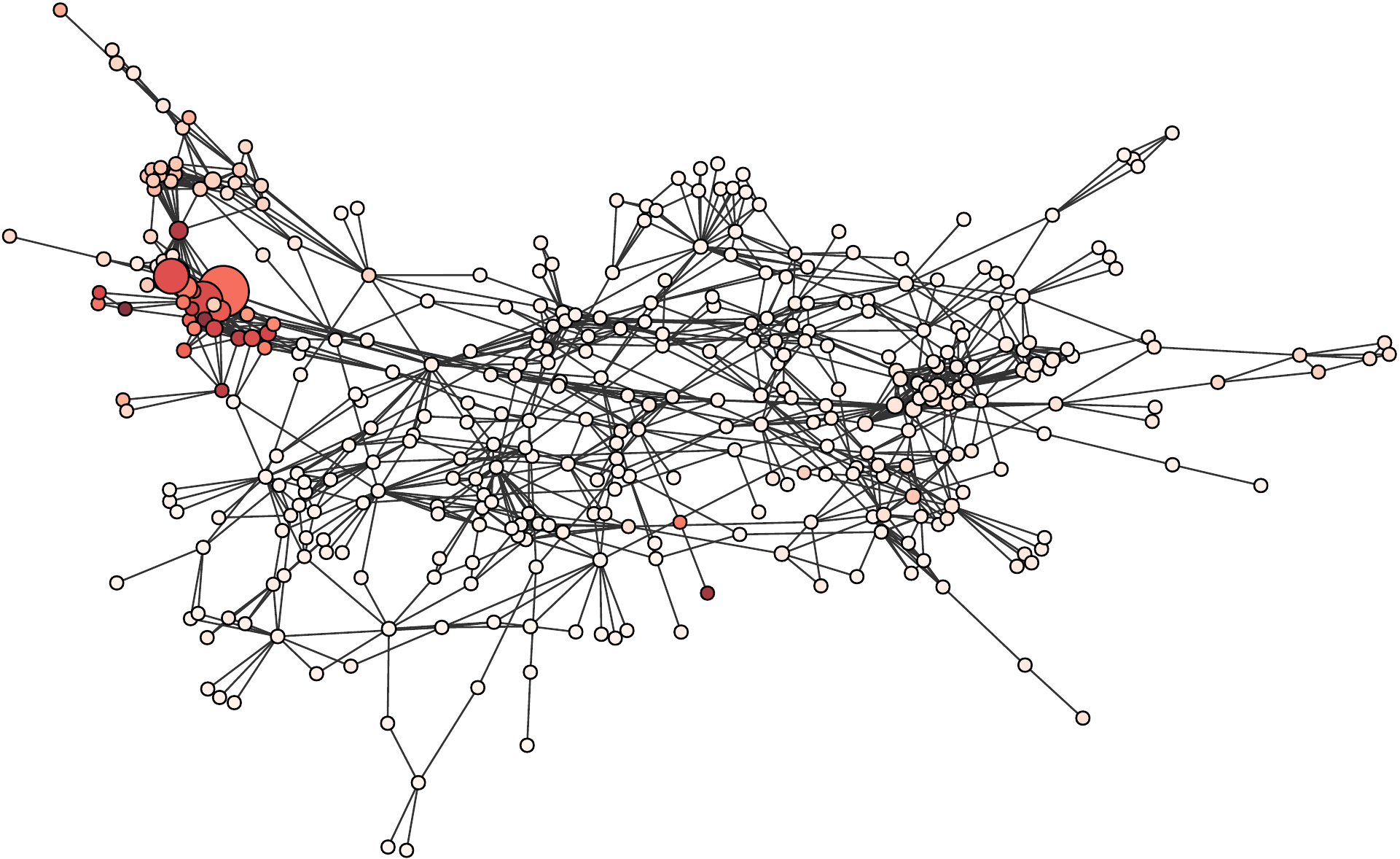}
         \caption{}%BM graph for the composition dataset.}
     \end{subfigure}
     \caption{Ball Mapper graphs of the superconductor characteristic (A) and composition dataset (B) for superconductor datasets available from the UC Irvine 
Machine Learning Repository,  colored by critical temperature. }
    \label{fig:superconducting_dataset_BM}
\end{figure}

Next we can compare structures of the two datasets characterizing superconductors using MappingMappers since the data is related via the map described in
~\cite{superconductor_dataset}. MappingMappers (figure not presented for brevity) shows that the image of a single vertex in the Ball Mapper graph of one of the datasets contains multiple disconnected regions in the second one. These results suggest that the characteristic and composition data provide distinct and likely unrelated information; hence, when analyzed together, these should provide more insights into properties of superconductors. 

\subsection{Mapper on Ball Mapper}
\label{sec:additional_mobm}
Last but not least, we illustrate the performance and properties of Mapper on Ball Mapper on both synthetic and real-word datasets. 

\subsubsection{Nonlinear datasets}

Our first artificial dataset is based on a set $C$ obtained by sampling $500$ points from the nonlinear embedding $F:\mathbb{R}^2 \rightarrow \mathbb{R}^7 $ defined by $F(x,y) = (xy , x^2 , y^2 , x^2y \ y^2x , x^3 , y^3) $ 
 of the unit circle in $\mathbb{R}^2.$ The Ball Mapper graph of $C$ is a cyclic graph, recovering the topology of the original unit circle. Then we take the Cartesian product of $C \times L$ where $L$ is a uniform sample of $100$ points from the unit interval $[0,10]$. This data set consists of $500 \cdot 100 = 50000$ points in  $\mathbb{R}^8$
and its shape comes from a cylinder with a unit base circle embedded in $\mathbb{R}^8$ via the product of $F$ with the identity map. 

The Ball Mapper graph of  $C \times L$, shown on Figure \ref{fig:MoBM_circle_A}, recovers the cylindrical shape for a  sufficiently small radius. 
However, the resulting Ball Mapper graph contains a large number of nodes and could be difficult to interpret. More importantly, if the radius is increased a bit, the base cycle in the Ball Mapper graph would collapse to a point and the information about the central hole would be lost. 

\ifdefined\ifNotPictures
\begin{figure}[h!]
     \centering
     \begin{subfigure}[t]{0.4\textwidth}
         \centering
         \includegraphics[height=0.6\textwidth]{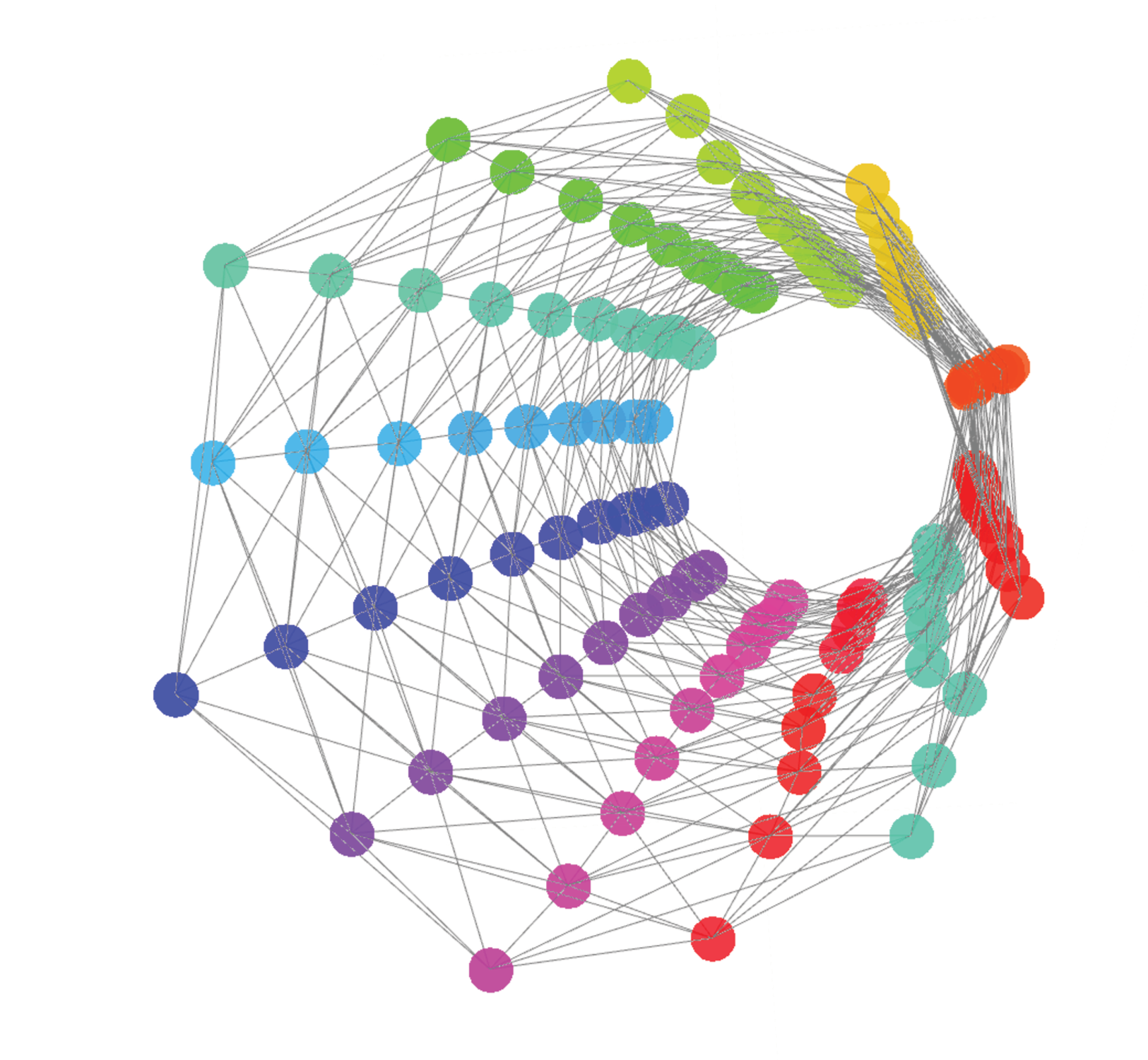}
         \caption{}%3d visualization of the BM graph on $C \times L$.}
         \label{fig:MoBM_circle_A}
     \end{subfigure}
     \hspace{1cm}
     \begin{subfigure}[t]{0.4\textwidth}
         \centering
         \includegraphics[height=0.55\textwidth]{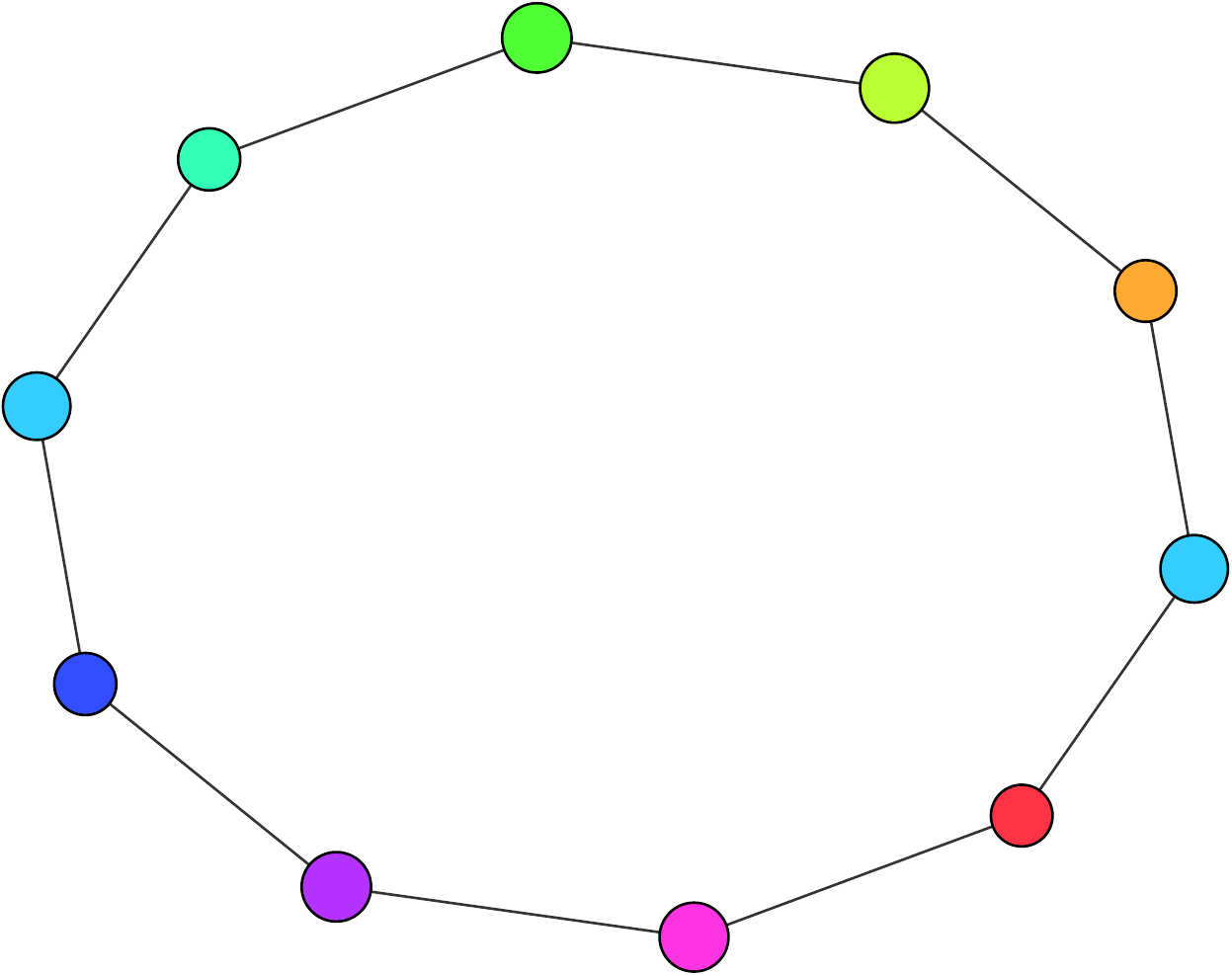}
         \caption{}%Mapper on BM graph from $C \times L$ to $C$.}
         \label{fig:MoBM_circle_B}
     \end{subfigure}
    \caption{The Ball Mapper graph obtained from $C \times L$ data is shown in (A). A sufficient small $\epsilon$ has to be chosen in order to observe the central hole, which leads to a lot of nodes being created. On the other hand, the Mapper on Ball Mapper graph obtained by pulling back the points covered by the nodes in the Ball Mapper graph on $C$ to $C \times L$ is presented in %Figure~{fig:MoBM_circle}
    (B). The coloring on both graphs is a function of the angle $\alpha=cos^{-1} (x)$ \NEW{and it used to show the correspondence between the Ball Mapper of the space $C$ and the MoBM for the space $C \times L$}. } 
    \label{fig:MoBM_circle}
\end{figure}

Since these two datasets are related via the projection $\pi : C \times L \subseteq \mathbb{R}^8 \rightarrow C \subseteq \mathbb{R}^7$, which takes each point in $C \times L$ to its component in $C$, we can apply the
Mapper on Ball Mapper construction on $C$ and $C \times L$ using the projection map $\pi$ as a lens.
The Mapper on Ball Mapper graph, shown on Figure \ref{fig:MoBM_circle}(B), captures the shape of the base space $C$ which is a cycle. Moreover, each copy of $L$ in $C \times L$ contains only one cluster, hence can be considered to be a single connected component. The Mapper on Ball Mapper graph, shown on Figure \ref{fig:MoBM_circle}(B), recovers the correct homotopy type of $C \times L$, which is the same as that of the base space $C$. Moreover, each fiber is a copy of $L$ corresponding to a single vertex/cluster in the Mapper on Ball Mapper graph, consistent with contractibility of the interval $L$ was sampled from. Note that this is not always the case, as connected fibers may have nontrivial topology.

\ifdefined\showOldText 
\OLD{The second example focuses on the standard two dimensional torus, the product of two unit circles $S^1 \times S^1$ with the standard parametrization: }
\begin{equation}
    \begin{cases}
    x = (R_1 + R_2\ \cos(v)) \cos(u)
    \\
    y = (R_1 + R_2\ \cos(v)) \sin(u)
    \\
    z = R_2\ \sin(v)
    \end{cases}.
\end{equation}
\OLD{
where $u,v \in [0,2\pi)$ and $R_1 > R_2$. Fixing $v = 0$ and letting $u \in [0,2\pi)$ determines the so-called base cycle 
$ ( (R_1+R_2)\cos(u) , (R_1+R_2)\sin(u) , 0 )$. }

\OLD{Consider the Mapper on Ball Mapper construction performed for the projection of torus to the base circle. As in the previous example, the Ball Mapper graph of the circle is a cyclic graph. Moreover, the inverse image of each vertex from the Ball Mapper graph of the basic circle has one connected component, hence, the Mapper on the Ball Mapper graph is isomorphic to the Ball Mapper graph of the basic circle. To summarize, both lens functions have connected fibers, corresponding to a single cluster in their Ball Mapper graphs, which are isomorphic to the Reeb graphs of the base space. Unlike the first example, where the fiber is contractible, the fiber in the second case is a circle, but the Mapper on Ball Mapper construction detects only the number of connected components. }
\fi

\NEW{
\subsubsection{Breast cancer dataset}
As a real-word application of the MoBM technique, we analyze the NKI breast cancer gene-expression dataset~\cite{van_de_vijver_gene-expression_2002}. Such dataset has already been analyzed in a seminal Mapper paper~\cite{nicolau} and can be downloaded from \url{https://data.world/deviramanan2016/nki-breast-cancer-data}. The data consist of the $1553$ top varying genes for $272$ breast cancer patients along with clinical metadata such as patient's age, type of treatment and survival. Our goal is to build a Mapper graph that differentiates between subgroups of patients with different survival rates. 
The first attempt to construct a Ball Mapper over the whole data set using the correlation distance (since gene expressions are given as z-scores) fails to provide any information, see Figure~\ref{fig:mobm_nki}.
The same holds for BM graphs constructed using other metrics. In the next attempt we project the data to a $200$-dimensional space using UMAP~\cite{mcinnes2018umap-software, 2018arXivUMAP}. The Ball Mapper graph obtained from the projected data using Euclidean distance  can be found in Figure~\ref{fig:mobm_nki}B. Finally, we pull this graph back to the original space using Mapper on Ball Mapper, thus obtaining the graph in Figure~\ref{fig:mobm_nki}C. In particular, we are able to identify a distinct subgroup of $44$ patients with 100\% survival rate. While the detailed medical interpretation of such a finding is beyond the scope of this paper, the presented pipeline shows how methods introduced in this paper can be used in combination with other techniques to analyze noisy, high dimensional data.

\begin{figure}[h!]
     \centering
     \begin{subfigure}[b]{0.3\textwidth}
         \centering
         \includegraphics[width=\textwidth]{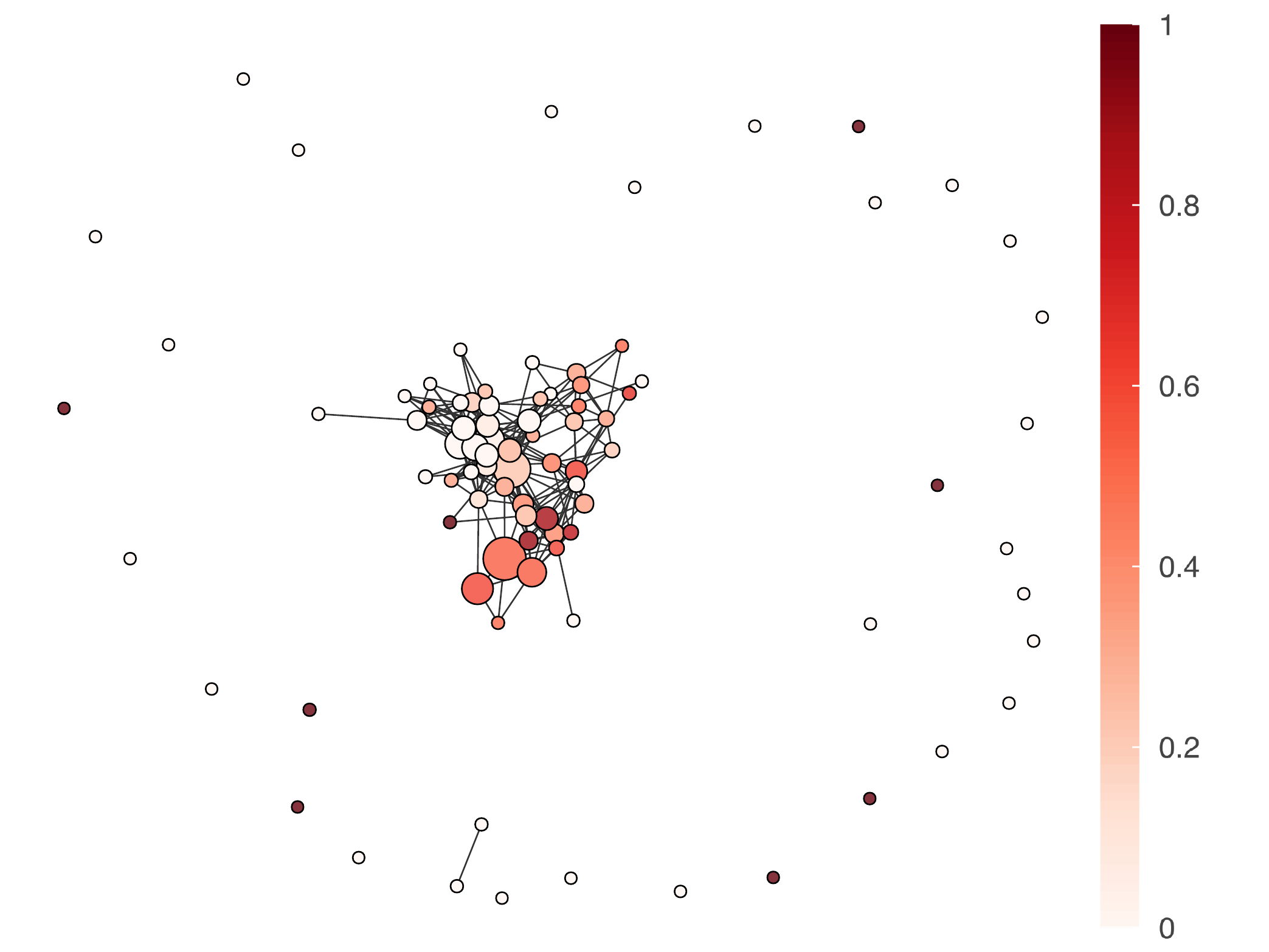}
         \caption{}
         \label{fig:mobm_nki_A}
     \end{subfigure}
     \hfill
     \begin{subfigure}[b]{0.3\textwidth}
         \centering
         \includegraphics[width=\textwidth]{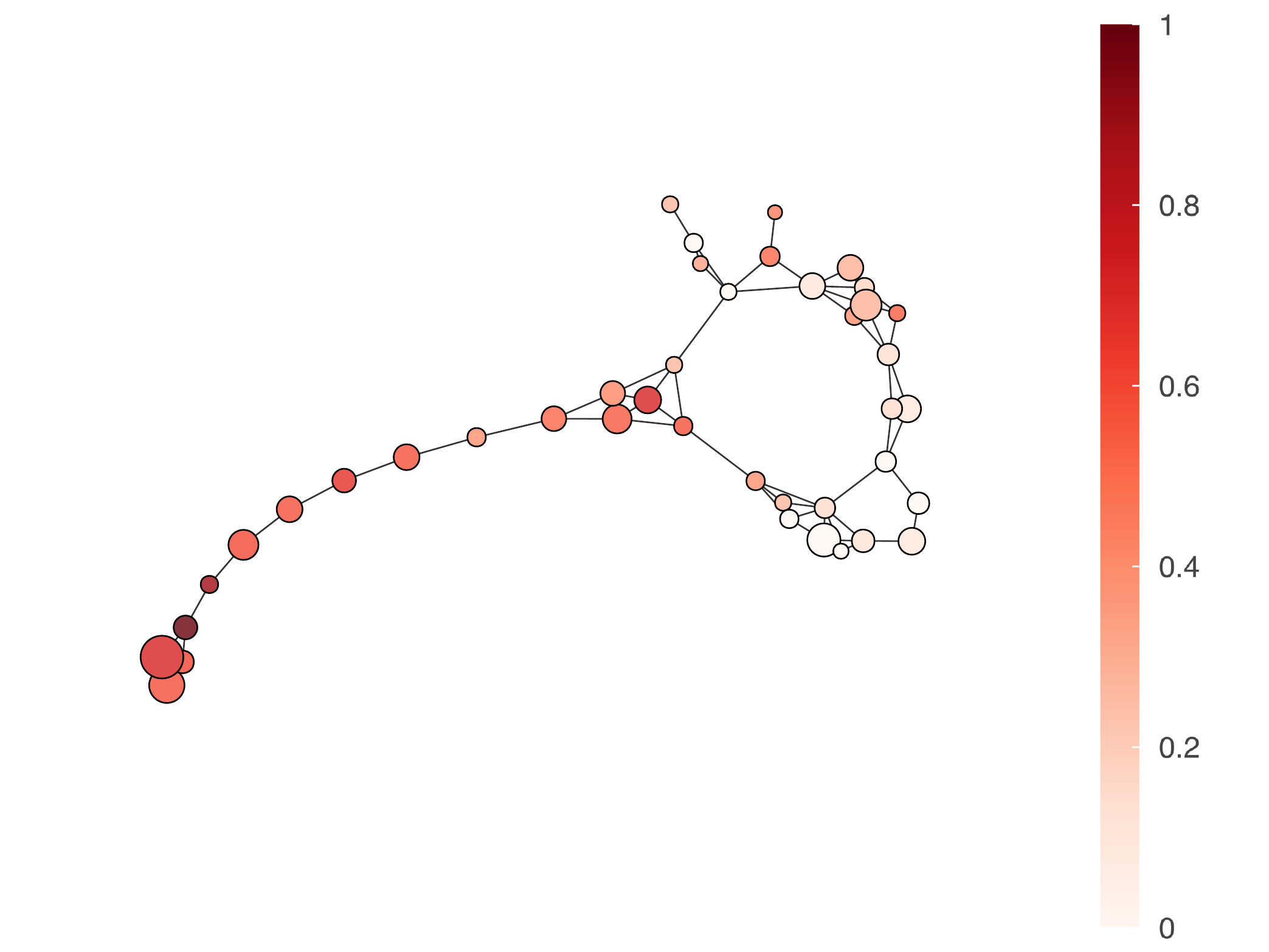}
         \caption{}
         \label{fig:mobm_nki_B}
     \end{subfigure}
    \hfill
     \begin{subfigure}[b]{0.3\textwidth}
         \centering
         \includegraphics[width=\textwidth]{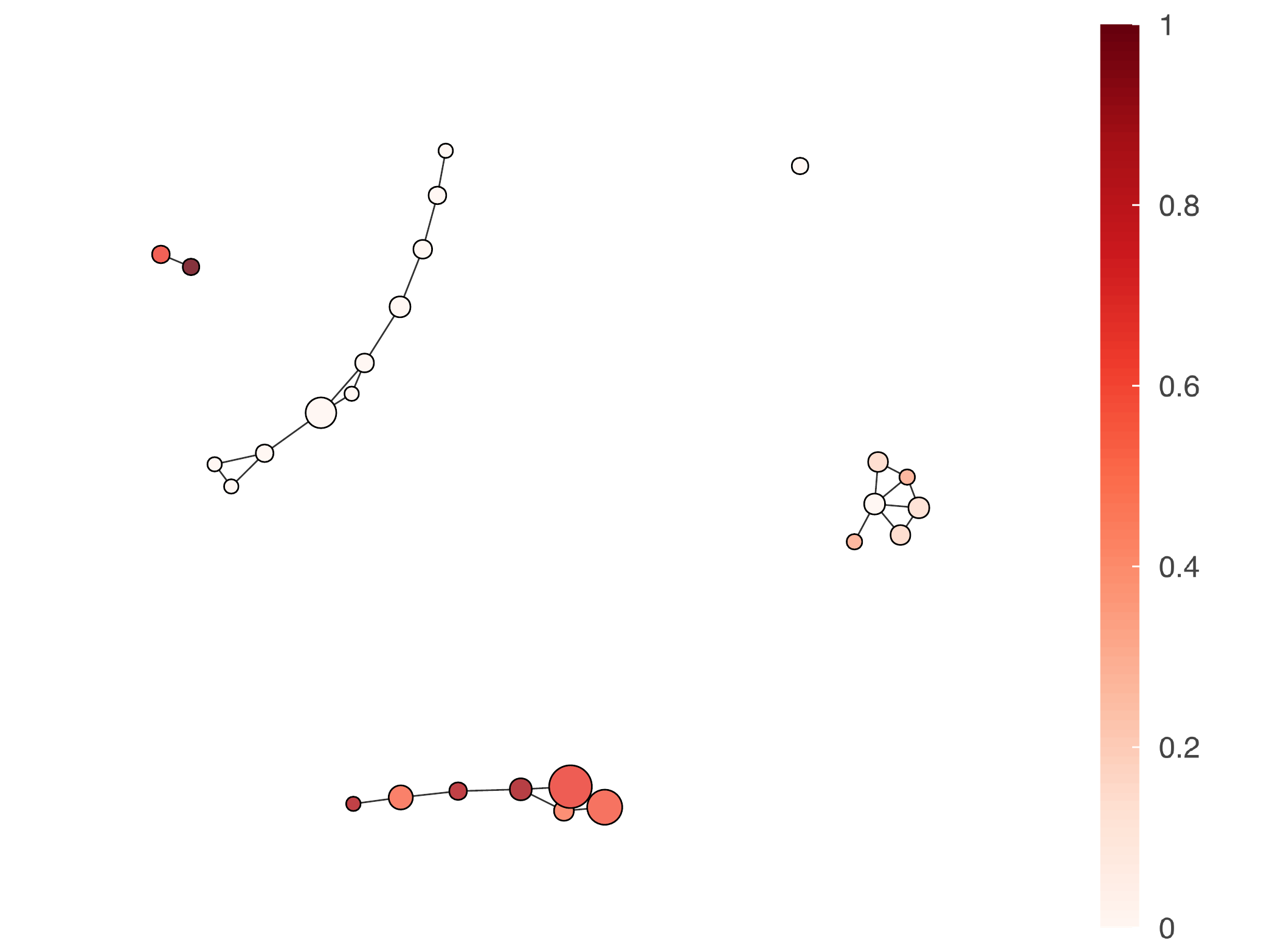}
         \caption{}
         \label{fig:mobm_nki_C}
     \end{subfigure}
        \caption{Three different Mapper-type graphs on the NKI dataset, nodes are colored by the fraction of patients who did not survive. Figure (A) depicts the Ball Mapper graph over the entire space, and it is not informative. Figure (B) depicts the Ball Mapper graph over a $200$-dimensional projection obtained using UMAP. Figure (C) depicts the Mapper on Ball Mapper from Figure (B) to the original $1553$-dimensional space. A clear separation between subgroups with different survival rates can be observed in Figure (C). In the last picture clusters containing only one or two elements have been removed for clarity.}
        \label{fig:mobm_nki}
\end{figure}

}

\ifdefined\showOldText 
\section{Summary}
\OLD{
This paper focuses on extensions of Mapper and Ball Mapper algorithms:  Equivariant Ball Mapper takes into account the structure of the data,  MappingMappers allows us to visualize maps between high dimensional datasets and enable comparisons, while Mapper on Ball Mapper extends the Mapper algorithm from 1-dimensional to high dimensional lens functions which are more likely to preserve and reveal information about the input data. To demonstrate the utility and applicability of these new methods across sciences we include a variety of examples accompanied by a public domain implementation.
}
\fi
%-----------------------------------------------

%
%
%

\bibliographystyle{plain}
\bibliography{references}
\end{document}